\documentclass[singlespace]{easychithesis}

\usepackage{amsmath, amsfonts, graphics, color, amsthm, amssymb, amsopn,
pstricks, multido, pst-node, rotating, epsfig}

\theoremstyle{definition}
 \newtheorem{thm}{Theorem}[section]
 \newtheorem{cor}[thm]{Corollary}
 \newtheorem{lem}[thm]{Lemma}
 \newtheorem{prop}[thm]{Proposition}
 \newtheorem{defn}{Definition}[section]
 \newtheorem{conj}{Conjecture}[section]

\newcommand{\inv}{\mathop{inv}}

\begin{document}

\title{Determining Whether Certain Affine Deligne-Lusztig Sets
are Empty}
\author{Daniel C. Reuman}
\date{August 2002}
\department{Mathematics}
\division{Physical Sciences} 
\degree{Doctor of Philosophy} 
\maketitle





\begin{abstract}
Let $F$ be a non-archimedean local field, let $L$ be the maximal unramified
extension of $F$, and let $\sigma$ be the Frobenius automorphism.  
Let $G$ be a split connected reductive group over $F$, and 
let $\mathcal{B}_{\infty}$ be the Bruhat-Tits building associated to $G(L)$.
Let $\mathcal{B}_1$ be the building associated to $G(F)$. 
Then $\sigma$ acts on $G(L)$ with fixed points $G(F)$. 
Let $I$ be the Iwahori associated to a chamber in $\mathcal{B}_1$. 
We have the relative position map
$\inv : G(L)/I \times G(L)/I \rightarrow \tilde{W}$, where $\tilde{W}$ is the
extended affine Weyl group of $G$. If $w \in \tilde{W}$ and $b \in G(L)$, 
then the affine Deligne-Lusztig set $X_w(b \sigma)$ is 
$\{ x \in G(L)/I : \inv(x,b \sigma (x)) = w$.  This dissertation
answers the question of which $X_w(b \sigma)$ are non-empty for certain
$G$ and $b$.

\end{abstract}

\topmatter{Acknowledgements}
Above all, I would like to thank my advisor, Robert Kottwitz, for being 
extremely helpful at every turn.  Professor Kottwitz's mathematical 
guidance has helped with every aspect of this dissertation.  His support
and understanding regarding non-research issues encountered during 
graduate education is also remarkable.  I am very 
grateful to have been his student.  I would
also like to thank Chris Degni and Sudheer Shukla for helpful discussions
and lectures.  Finally, I thank Moon Duchin, Kaj Gartz, Mark Behrens, and
Moses Hohman for moral support. 

%
%

\tableofcontents

%
%

\listoffigures

%
%

\listoftables

%
%

\mainmatter

%
%

\chapter{Introduction}

Let $F$ be a non-archimedean local field with ring of integers 
$\mathcal{O}_F$, and let $G$ be a split connected reductive group
over $F$.  Let $L$ be the completion of the maximal unramified extension
of $F$.  Let $\sigma$ be the Frobenius automorphism of $L$ over $F$.  Let
$\mathcal{B}_n$ be the affine building for $G(E)$ where $E/F$ is the 
unramified extension of degree $n$ in $L$, and let
$\mathcal{B}_{\infty}$ be the affine building for $G(L)$.  The theory of
buildings was developed by Bruhat and Tits in \cite{BT1} and \cite{BT2}, and 
is reviewed by Garrett in \cite{G}.  We know that $\sigma$
acts on $\mathcal{B}_{\infty}$, and the  
fixed points of $\sigma^n$ are $\mathcal{B}_n$. 
Let $T$ be a split torus in $G$ and let $I$ be an Iwahori in $G(L)$
containing $T(\mathcal{O}_L)$, where $\mathcal{O}_L$ is the ring of integers
of $L$.  Let $A_M$ and $C_M$ be the correspondingly specified
apartment and chamber, which we assume are in $\mathcal{B}_1$. 
We will call these the {\em main
apartment} and the {\em main chamber,} respectively.

If $b \in G(L)$ then the $\sigma$-conjugacy class of $b$ is $\{x^{-1}b\sigma
(x) : x \in G(L) \}$.  If $\tilde{W}$ is the extended affine Weyl
group, and $w \in \tilde{W}$ then we define, after Rapoport and Kottwitz,
the affine Deligne-Lusztig set $X_w(b\sigma) = \{x \in G(L)/I :
\inv(x, b\sigma(x)) = w \}$.  Here $\inv : G(L)/I \times
G(L)/I \longrightarrow \tilde{W}$ is the relative position map.  Note that
if $b_1$ and $b_2$ are $\sigma$-conjugate then the sets $X_w(b_1 \sigma)$
and $X_w(b_2 \sigma)$ are in bijective correspondence.  When 
$F = \mathbb{F}_q ((t))$, the affine Deligne-Lusztig set
 can be given variety structure over $\mathbb{F}_q$ (locally of finite
type).

Sometimes the $X_w(b \sigma)$ are empty.  
Our goal is to determine when the $X_w(b
\sigma)$ are non-empty.  This paper concerns itself mainly with 
the question of which
$X_w(b \sigma)$ are non-empty for groups of semisimple rank $1$ and $2$.

The question of which $X_w(b \sigma)$ are non-empty for 
$\{ x^{-1} b \sigma (x) : x \in G(L) \}$ a fixed $\sigma$-conjugacy 
class can be 
rephrased by saying that $X_w (b \sigma)$ is non-empty if and only if $\{ 
x^{-1} b \sigma (x) : x \in G(L) \}$ intersects the double-$I$-coset
corresponding to $w$ in a non-trivial way.  Here we are using the standard 
correspondence between $\tilde{W}$ and $I\backslash G(L) /I$.  We will use
this correspondence implicitly for the rest of the paper.  To see that the 
above rephrasing is accurate, note first that $\inv(x, b \sigma (x)) = 
\inv(1, x^{-1} b \sigma (x))$ is the relative position of $1$ and 
$x^{-1}b \sigma (x)$.  By definition, $\inv(1, x^{-1} b \sigma (x)) = 
Ix^{-1}b\sigma(x)I$.  We see that
$X_w (b \sigma)$ is non-empty if and only if there exists $x$ such that 
$Ix^{-1}b \sigma(x) I = w$, which is what was to be proven.

If $G$ is simply connected, then $\tilde{W} \simeq W_a$ 
is in bijective correspondence
with chambers in the main apartment $A_{M}$.  Here $W_a$ is the affine
Weyl group of $G$. So for a fixed 
$\sigma$-conjugacy class $\{x^{-1}b \sigma(x) : x \in G(L) \}$, the set of 
$w$ for which $X_w (b\sigma)$ is non-empty can be represented geometrically as
a collection of chambers in $A_M$.  For groups of rank $1$ and $2$, one can 
therefore draw pictures of the solution set.  
This is the form in which we present results in this paper.

Complete information can be obtained without difficulty for groups
of semisimple rank $1$.  This is done in Chapter~\ref{RankOne}.  The solution
sets for each $\sigma$-conjugacy class of $SL_2$ are pictured in 
Figure~\ref{Summary}.  The solution sets for each $\sigma$-conjugacy class of
$GL_2$ or $PGL_2$ are pictured in Figures~\ref{GL2Summary1} and~\ref{Summ}.

The group of semisimple rank $2$ for which the most information 
has been obtained is $G = SL_3$.  We have worked only with 
$\sigma$-conjugacy classes that are basic in the maximal torus.  
These are listed at the beginning of Section~\ref{SL3}.  
The general approach by which we have described 
$\{\inv(x,b \sigma(x)) : x \in SL_3(L) \} \subseteq 
\tilde{W} \simeq W_a$ in these cases has two main parts.  The first part,
done in Sections~\ref{SMGandComp},~\ref{Superset}, and~\ref{Relationship3}, 
seeks to produce a superset 
of the solution set.  The general methodology here is to establish a 
way of choosing, for each $x \in SL_3(L)$, a gallery $\Gamma_x$ 
between $xC_M$ and $b\sigma(x)C_M$.  Possible values of 
$\inv(x,b\sigma(x))$ can then be enumerated by listing the possible
``foldings'' of $x^{-1}\Gamma_x$.  Optimizations can be implemented so that 
this is a finite, but still prohibitively lengthy computation.
It has only been carried out to completion for two $\sigma$-conjugacy
classes.  The results are pictured in Figure~\ref{I1I2Whole} and 
Figure~\ref{TotalResults0_0}.
In the process of doing these two complete computations, we observe 
that in fact one need only do a much smaller part of the complete
computation to obtain all elements of $W_a$ that arise in the superset.  
We do not have a proof that this is the case in general, 
but we are very strongly 
confident that it is.  Section~\ref{Relationship3} 
develops an efficient way of 
doing the smaller, partial computation that works for 
any $\sigma$-conjugacy class.  Results of the smaller computation 
that are (confidently) conjectured to the be results of the entire
computation are listed for various $\sigma$-conjugacy classes in 
Figures~\ref{SomeResults2_0},~\ref{SomeResults2_1},~\ref{SomeResults4_0}, and~\ref{SomeResults4_minus2}.

The other half of the work for $SL_3$ is to develop a subset of the solution
set.  This is done in two different ways. The first method, which is discussed
in Sections~\ref{Rapoport} and~\ref{Subset}, is simple 
and potentially generalizable 
to higher rank groups, but only applies to the case $b = 1$.  The second method
(discussed in Section~\ref{SubsetGeometric}) applies for any $b$ in 
the collection of $\sigma$-conjugacy classes
under consideration, but is 
very complicated and somewhat aesthetically unpleasing.  Also, it 
probably could not be generalized to higher 
rank groups.  In all cases, the subset 
results turn out to be equal to the superset results discussed above.  
So it is known that Figures~\ref{I1I2Whole} and~\ref{TotalResults0_0} 
are solution
sets for their respective $\sigma$-conjugacy classes, and we are
very confident (although it is not proved) that 
Figures~\ref{SomeResults2_0},~\ref{SomeResults2_1},~\ref{SomeResults4_0}, and~\ref{SomeResults4_minus2} are solution sets for their
respective $\sigma$-conjugacy classes. It is known that these last four
figures are at least subsets of their respective solution sets.

The same methods can be applied to $Sp_4$, although this has been carried out
to a smaller extent.  Again, we only consider $\sigma$-conjugacy classes
$b$ which are basic in the maximal torus.  
The complete computation of the superset of the solution set has only
been done for $b = 1$.  The result is pictured in 
Figure~\ref{SupersetResultalpha0beta0Sp4}.  Again,
there is a partial computation which gives all results obtained from the 
complete computation in this case.  Carrying out
the partial computation for some other $\sigma$-conjugacy classes (done 
in Section~\ref{SupersetSp4}) yields
Figures~\ref{SupersetResultalpha3beta1Sp4},~\ref{SupersetResultalpha6beta3Sp4},~\ref{SupersetResultalpha6beta5Sp4}, and~\ref{SupersetResultalpha7beta1Sp4}.  We conjecture that these
are complete supersets, but we are not completely confident
of this because the chambers marked with a $*$ seem to be holes in 
the pattern.

The subset methods of Section~\ref{Rapoport} and~\ref{Subset} 
have been applied 
to $Sp_4$, $b = 1$
in Sections~\ref{RapoportSp4} and~\ref{SubsetSp4}.  
The result is the same as the 
superset pictured in
Figure~\ref{SupersetResultalpha0beta0Sp4}, so 
this figure gives the solution set for $b = 1$.
The complicated methods of Section~\ref{SubsetGeometric} 
have not been adapted to
$Sp_4$, although such an adaptation probably could be done.  
We conjecture that one
would obtain subsets equal to the sets pictured in 
Figures~\ref{SupersetResultalpha3beta1Sp4},~\ref{SupersetResultalpha6beta3Sp4},~\ref{SupersetResultalpha6beta5Sp4}, and~\ref{SupersetResultalpha7beta1Sp4}, making these figures actual 
solution sets.  So for $Sp_4$ we have definite
results for $b = 1$, and conjectural results of which we are not 
completely confident for other $b$.

The results obtained for $SL_3$ can easily be adapted to $GL_3$ and 
$PGL_3$.  This is done in Section~\ref{GL3}.  Similarly, the results for $Sp_4$
can be adapted to $GSp_4$ and $PSp_4$.  This is done in Section~\ref{GSP4}.
The methods used could be applied to $G_2$ as well.  This is also discussed 
in Section~\ref{GSP4}.

We also say a few words describing the basic ideas behind the methods used.  
The backbone of the method by which supersets of solution sets are
produced is the ``folding'' of galleries.  If $G$ is a gallery starting
at $C_M$, then the {\em folding results} of $G$ are the possible
values of the last chamber of $\rho(\tilde{G})$, where 
$\rho : \mathcal{B}_{\infty} \rightarrow A_M$ is the retraction
centered at $C_M$, and $\tilde{G}$ is any gallery of the same
``shape'' as $G$ that starts at $C_M$.  To produce a superset, it turns out
that the folding results have to be calculated for a gallery associated to
each chamber in $\mathcal{B}_{\infty}$.  Perhaps the real content of the
method is how these computations can be grouped so that they can be done
efficiently.  This will be discussed in Sections~\ref{Superset} 
and~\ref{Relationship3} for 
$SL_3$, and Sections~\ref{SupersetSp4} and~\ref{RelationshipSp4} for $Sp_4$.

As mentioned previously, there are two methods used to produce subsets
of solution sets.  The first, which only applies for $b = 1$, is
carried out in Sections~\ref{Rapoport} and~\ref{Subset} for $SL_3$, and 
Sections~\ref{RapoportSp4} and~\ref{SubsetSp4} for $Sp_4$.  
This method gets started using an
observation of Kottwitz and Rapoport that provides, with very little work,
a large and well distributed collection of chambers which must be contained
in the solution set.  This collection can then be enlarged as follows.
If $\inv(x,\sigma(x)) = w$ is one of the chambers provided by Kottwitz's
and Rapoport's observation, and if $\Gamma$ is a minimal 
gallery between $xC_M$ 
and $\sigma(x)C_M$, then $\Gamma$ describes the relative position $w$.
One can show that $\sigma$-conjugating $x$ by certain
very simple elements of $G(L)$ corresponds to adding ``appendages''
to the ends of $\Gamma$.  One can then understand
how the addition of such appendages affects relative position.  Some of the 
``appendage'' methods used are stated in generality in 
Section~\ref{Invariance}.

The second method of producing subsets of solution sets is lengthy
and unappealing, but applies for $b \neq 1$.  The basic idea is as 
follows.  In Section~\ref{SMGandComp}, we produce a gallery $\Gamma_x$ between
$xC_M$ and $b\sigma(x)C_M$.  We compute the folding results of 
$x^{-1}\Gamma_x$ to get possible candidates for the relative
position of $x$ and $b\sigma(x)$. In Section~\ref{SubsetGeometric}, we discover
which of these candidates can actually occur by defining some invariants 
which control how $x^{-1}\Gamma_x$ actually folds (i.e., these invariants 
control which of the folding results of $x^{-1}\Gamma_x$ is the actual
value of the last chamber of $\rho(x^{-1}\Gamma_x)$).  We then show
that one can choose $x$ such that these invariants take any pre-specified 
value.  One of the complications of this method is that one must proceed 
somewhat differently for each of the three cases $b = 1$, $b$ 
``degenerate'' but $b \neq 1$, and $b$ ``non-degenerate.''  The 
``non-degenerate'' $b$ are exactly
$$b = \left(
\begin{matrix}
	\pi^{\alpha} & 0 & 0 \\
	0 & \pi^{\beta} & 0 \\
	0 & 0 & \pi^{\gamma}
\end{matrix} \right)$$
with $\alpha + \beta + \gamma = 0$, $\alpha > \beta > \gamma$.

\chapter{Groups of Semisimple Rank $1$}\label{RankOne} 

If the split connected reductive group $G$ is of semisimple rank $1$, 
then it is possible to obtain a complete
characterization of which $X_w(b\sigma)$ are non-empty.  Some of the methods
which will be used on higher rank groups are similar to those used here.
Considering the rank $1$ case first allows one to develop 
a good intuition
about the geometric nature of the problem before proceeding to the more
complicated rank $2$ and higher rank cases.

Section~\ref{SL2} of this chapter will address the specific group $SL_2$, and 
Section~\ref{GL2} will adapt the $SL_2$ solution to other, non-simply connected
groups of semisimple rank $1$.

Everything in this chapter was previously known to Kottwitz and Rapoport. 

\section{$SL_2$}\label{SL2}

Results from a paper of Kottwitz \cite{K1} allow us to list 
the following elements of $SL_2$ as a complete collection of representatives
of the $\sigma$-conjugacy classes of $SL_2$:
$$\left\{ \left( 
\begin{matrix}
	\pi^n & 0 \\
	0 & \pi^{-n}
\end{matrix}
\right) : n \in \mathbb{Z}, n \geq 0 \right\}. $$
Let $\rho : \mathcal{B}_{\infty} \rightarrow A_M$ be the 
retraction of $\mathcal{B}_{\infty}$
onto $A_M$ relative to $C_M$.

In Section~\ref{Identity}, we answer the question of which $w$ have non-empty
$X_w (b \sigma)$ for $b$ the identity matrix.  In Section~\ref{Non-identity}, 
we consider the cases where $b$ is one of the above non-identity 
representatives.  In Section~\ref{Symmetry}, we prove that the results obtained
in Sections~\ref{Identity} and~\ref{Non-identity} are invariant under a certain
$\mathbb{Z}/2$-action.  The proof is independent of the results from 
Sections~\ref{Identity} and~\ref{Non-identity}.  Section~\ref{Relationship}
describes the relationship between $X_w (1\sigma)$ and $X_w (b \sigma)$ for 
$b \neq 1$.

\subsection{Identity $\sigma$-conjugacy Class}\label{Identity}

As mentioned previously, to find out which $w$ have non-empty $X_w (b \sigma)$,
it suffices to determine the set $\{ \inv(1,x^{-1}b \sigma(x)) : x \in SL_2(L) 
\}$.  However, it is easy to see that $\inv(1, x^{-1}b \sigma (x) ) = 
\rho (x^{-1}b \sigma(x) C_M)$, where again it is understood that we 
are making use of the correspondence between $\tilde{W}$ and the set
of chambers in $A_M$.  
So we seek to describe the set $\{ \rho ( x^{-1} b 
\sigma (x) C_M ) : x \in SL_2(L) \}$.

Let $D$ be any chamber in $\mathcal{B}_{\infty}$.  

\begin{defn} Let the distance $\tilde{d}_1 (D)$ 
between $D$ and $\mathcal{B}_1$ be $n-1$, where $n$ is the length of 
the shortest gallery containing $D$ and some chamber $E$ in $\mathcal{B}_1$.
\end{defn}

\noindent Note that $D$ is in $\mathcal{B}_1$ if and only if $\tilde{d}_1 (D) = 0$.

\begin{prop} If $\tilde{d}_1 (D) > 0$ then there 
is a minimal gallery between $D$ and $\sigma(D)$ of length $2 \tilde{d}_1 (D)$.
\end{prop}

\begin{proof} Let $G$ be a gallery of minimal length containing 
$D$ and some chamber $E$ in $\mathcal{B}_1$.  Let $E = G_1, G_2,G_3, \ldots,
G_n = D$ be the successive chambers of $G$.  Then $G_2$ is not contained 
in $\mathcal{B}_1$, or $G$ would not be of minimal length.  So the gallery
$G_n, G_{n-1}, \ldots , G_2, \sigma (G_2), \ldots , \sigma(G_{n-1}), 
\sigma (G_n)$ is non-stuttering and therefore minimal.
\end{proof}

Now let $D = xC_M$.  By the nature of the $SL_2$ action on the building,
$x$ can be chosen such that $\tilde{d}_1 (xC_M)$ is any non-negative integer.  
Therefore, $x$ can be chosen such that the minimal gallery between $xC_M$ and
$\sigma (x C_M) = \sigma(x) C_M$ is any even integer greater than or equal 
to $2$.  In the case where $\tilde{d}_1 (x C_M) = 0$ (i.e., $xC_M$ is 
contained in $\mathcal{B}_1$), we have $xC_M = \sigma (x C_M)$.

We are interested in $\rho (x^{-1} \sigma (x) C_M)$, so the minimal
gallery between $xC_M$ and $\sigma(x) C_M$ constructed above is useful.
We will need the following terminology to describe the set $\{ \rho
(x^{-1} \sigma (x)C_M) : x \in SL_2(L) \}$.

\begin{figure}
\centerline{\input{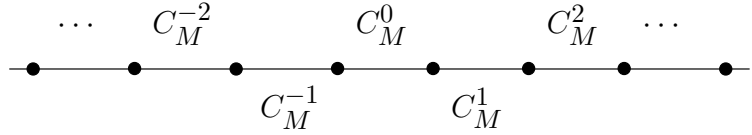}}
\caption{Chamber labels}
\label{ChamberLabels}
\end{figure}

\begin{defn} Let $C_M^i$ denote the chamber in $A_M$ as 
described in Figure~\ref{ChamberLabels}, where $C_M^0 = C_M$.
\end{defn}

\begin{figure}
\centerline{\input{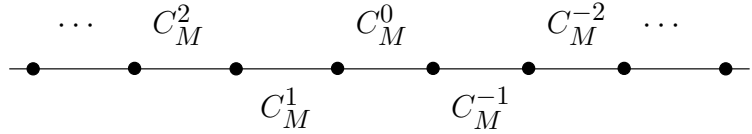}}
\caption{Alternate chamber labels}
\label{AlternateChamberLabels}
\end{figure}

\noindent Note that the choice between this labeling and the one in
Figure~\ref{AlternateChamberLabels} is arbitrary, but we choose the former.

\begin{prop} 
The $w \in \tilde{W}$ with non-empty
$X_{w}(1 \sigma)$ are $C_M$ and $C_M^i$ for $i$ odd.
\label{ResultId}
\end{prop}
\begin{proof} First, note that $C_M = \rho
(x^{-1}\sigma(x) C_M)$ for $x$ such that $xC_M$ is contained completely 
in $\mathcal{B}_1$.  Further, since $x$ can be chosen such that the minimal
gallery between $xC_M$ and $\sigma (xC_M)$ is of any even integral length, we
see that for each odd $i$, either $C_M^i$ or $C_M^{-i}$ is obtained.  We must
show that in fact both are obtained.  For this, note that if $\rho
(x^{-1}\sigma(x) C_M) = C_M^i$ then $\rho(y^{-1} \sigma(y) C_M) = 
C_M^{-i}$ for $y = gxq^{-1}$ where
$$g = \left(
\begin{matrix}
	\pi & 0 \\
	0 & 1 
\end{matrix}
\right) \hspace{.25in} ; \hspace{.25in} q = \left(
\begin{matrix}
	0 & - \pi \\
	1 & 0
\end{matrix} \right).$$  See Figure \ref{ProofHelp}.
\begin{figure}
\centerline{\input{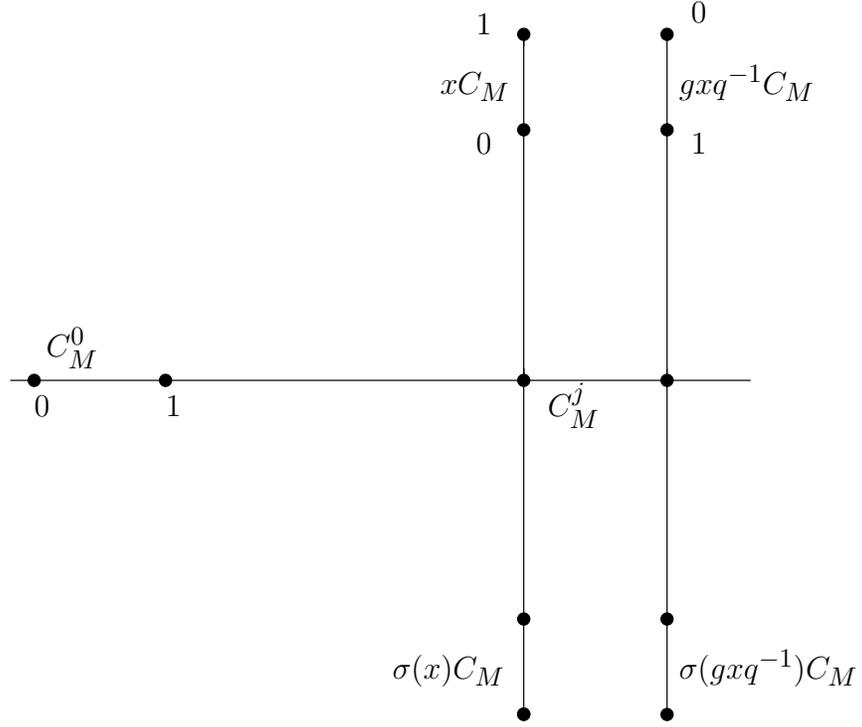}}
\caption{Proof of Proposition~\ref{ResultId}}
\label{ProofHelp}
\end{figure}
In this figure, the vertex types are labeled.  The 
key point is that the
vertices on the chambers $xC_M$ are of types opposite to the corresponding
vertices on $gxq^{-1}C_M$.  Note that we still have $\det (y) = 1$.  See
Section~\ref{Symmetry} for a reminder of how $q$ acts on $C_M$ and $A_M$.
\end{proof}

\subsection{Non-identity $\sigma$-conjugacy Classes}\label{Non-identity}

Again, let $D$ be any chamber in $\mathcal{B}_{\infty}$.  Let $b = \left(
\begin{matrix}
	\pi^s & 0 \\
	0 & \pi^{-s}
\end{matrix}
\right)$, $s > 0$.  We will enumerate the $w$ for which $X_w (b\sigma)$
is non-empty.

\begin{defn} The distance $\tilde{d}_{A_M}(D)$ between $D$
and $A_M$ is $n-1$, where $n$ is the length of the shortest gallery containing
$D$ and some chamber $E$ in $A_M$.
\end{defn}

\noindent Note that $D$ is in $A_M$ if and only if $\tilde{d}_{A_M}(D) = 0$.

\begin{prop} If $\tilde{d}_{A_M}(D) > 0$ then there exists
a minimal gallery between $D$ and $b\sigma(D)$ of length $2s +
2\tilde{d}_{A_M}(D)$.
\end{prop}

\begin{proof} Let $G$ be a gallery of minimal length containing
$D$ and some chamber $E$ in $A_M$.  Let $E = G_1, G_2, G_3, \ldots, G_n = D$
be the successive chambers of $G$.  Here $n = \tilde{d}_{A_M}(D)+1$.  Then 
$G_2$ is not contained in $A_M$, or $G$ would not be of minimal length.  Let
$v$ be the vertex that $E$ and $G_2$ share.  There is a unique minimal
gallery contained in $A_M$ between $v$ and $b\sigma (v) = bv$.  Call this
gallery $H$, and its chambers $H_1, H_2, \ldots , H_{2s}$, where $v$ is the
vertex of $H_1$ that is not also a vertex of $H_2$.  Consider the gallery
$G_n, G_{n-1}, \ldots , G_2, H_1, H_2, \ldots , H_{2s}, b\sigma (G_2),
b \sigma (G_3), \ldots , b \sigma (G_n)$. This is  non-stuttering and therefore
minimal.  It has the desired length.
\end{proof}

Now consider $D = xC_M$.  Clearly, $x$ can be chosen such that $\tilde{d}_{A_M}
(xC_M)$ is any non-negative integer.  Therefore, $x$ can be chosen such that
the minimal gallery between $xC_M$ and $b\sigma (xC_M)$ is of length
$2s + 2n$ for any $n \geq 1$.  If $xC_M$ is contained in $A_M$ then the 
minimal gallery between $xC_M$ and $b\sigma(xC_M)$ is of length $2s + 1$.

\begin{prop}
The $w \in \tilde{W}$ with non-empty
$X_w (b\sigma )$ are $C_M^{2s}$, $C_M^{-2s}$, $C_M^{2s+i}$ and $C_M^{-2s-i}$
for $i \geq 1$ any odd integer.
\end{prop}

\begin{proof} Similar to Proposition~\ref{ResultId} in 
Section~\ref{Identity}.
\end{proof}

Figure~\ref{Summary} summarizes results from this and the previous 
section.
\begin{figure}
\centerline{\input{fig4_2.pstex_t}}
\caption{Summary of $SL_2$ results}
\label{Summary}
\end{figure}

\subsection{Symmetry Under a $\mathbb{Z}/2$-action}\label{Symmetry}

\begin{figure}
\centerline{\input{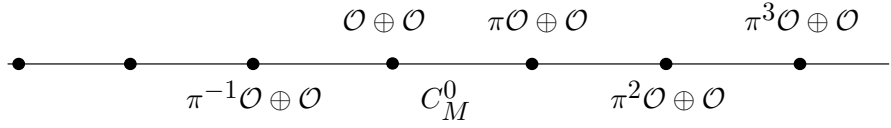}}
\caption{Lattice labels for vertices}
\label{Lattices}
\end{figure}

Let $q = 
\left(
\begin{matrix}
	0 & -\pi \\
	1 & 0
\end{matrix}
\right) \in GL_2 (L)$.
Since $\mathcal{B}_{\infty}$ is also the building for $GL_2(L)$, we 
can ask how $q$ acts on $A_M$.  The main apartment can be represented
using lattice classes as in Figure \ref{Lattices}.
The matrix $q$ takes 
$$\begin{matrix}
	\mathcal{O}\oplus\mathcal{O} & \rightarrow & 
	\pi\mathcal{O} \oplus \mathcal{O} \\
	\pi \mathcal{O} \oplus \mathcal{O} & \rightarrow &
	\mathcal{O} \oplus \mathcal{O} \\ 
	\pi^{-1} \mathcal{O} \oplus \mathcal{O} & \rightarrow &
	\pi^2 \mathcal{O} \oplus \mathcal{O} \\
	\pi^2 \mathcal{O} \oplus \mathcal{O} & \rightarrow &
	\pi^{-1} \mathcal{O} \oplus \mathcal{O}
\end{matrix}$$
and therefore represents a flip about the midpoint of $C_M$. One can see
by looking at the results of the previous two sections that $\{ w \in 
\tilde{W} : X_w (b \sigma) \neq \emptyset \}$ is invariant under this 
$q$ action.  It is instructive to develop a proof 
of this fact that is independent of the previous sections.

Since $\inv$ is a map from $SL_2(L)/I \times SL_2(L)/I \rightarrow
\tilde{W}$, and since $SL_2(L)/I$ is in one-to-one correspondence
with elements of the set of chambers in $\mathcal{B}_{\infty}$, we
can view $\inv$ as a map that takes two chambers in $\mathcal{B}_{\infty}$
as arguments.
We will use both this and the standard method of referring to $\inv$
throughout the rest of the paper.

We know that $\inv(xC_M, b\sigma(x)C_M) = \inv (gxC_M ,gb\sigma(x)C_M )$ for
any $g \in SL_2(L)$.  If $g \in SL_2 (F)$ is diagonal then
$\inv(xC_M, b\sigma (x)C_M ) = \inv(gxC_M,gb\sigma(x)C_M ) = 
\inv(gxC_M, b\sigma (gx)C_M )$.  In the case that $g \in GL_2(F)$
is diagonal, $\inv(xC_M , b\sigma(x)C_M )$ is not necessarily equal to
$\inv(gxC_M, gb\sigma (x)C_M )$, but $\inv(gxC_M,gb\sigma(x)C_M )
= \inv(gxC_M , b\sigma(gx)C_M)$.  If $y$ is such that $yC_M = gxC_M$, and
$y \in SL_2 (L)$, then $\inv(gxC_M,b\sigma(gx)C_M) = 
\inv(yC_M , b \sigma(y)C_M )$ is a new $w \in \tilde{W}$
such that $X_w (b\sigma) \neq \emptyset$.  We used the fact that if $yC_M = 
gxC_M$ then $b\sigma(gx)C_M = b\sigma(gxC_M) = b\sigma(yC_M) = 
b\sigma(y) C_M$. 

\begin{prop} We can choose $g$ a diagonal element of
$GL_2(F)$ such that if $\inv(xC_M,b\sigma(x)C_M) = C_M^i$, then
$\inv(gxC_M,b\sigma(gx) C_M) = C_M^{-i}$.
\end{prop}

\begin{proof} Let 
$g = \left(
\begin{matrix}
	\pi & 0 \\
	0 & 1
\end{matrix} \right).$  We see that $\inv(gx C_M , b\sigma (gx) C_M) = 
\rho(s^{-1}b\sigma(gx)C_M)$, where $s \in SL_2(L)$ is chosen 
such that
$s^{-1}gxC_M = C_M$.  We may choose $s$ such that $s^{-1} = qx^{-1}g^{-1}$.
So $\rho(s^{-1}b\sigma (gx)C_M) = \rho
(qx^{-1}g^{-1}bg\sigma(x) C_M) = \rho(qx^{-1}b\sigma(x)C_M)$, and if 
$\inv(xC_M , b\sigma(x) C_M) = C_M^i$, then 
$\rho(qx^{-1}b\sigma(x)C_M) = C_M^{-i}$.
\end{proof}

We will later adapt this proof to give a similar result for some groups of 
higher rank.

\subsection{Relationship Between $X_w(1\sigma)$ and $X_w(b\sigma)$}
\label{Relationship}

In this section we explore the relationship between the set of $w$ such
that $X_w(1\sigma)$ is non-empty and the set of $w$ such that $X_w(b\sigma)$
is non-empty, for some $b \neq 1$.

\begin{prop} If $w$ is in the positive Weyl chamber (i.e., $w = C_M^i$
for $i \geq 0$), 
then $X_w(1\sigma)$ is non-empty if 
and only if $X_{bw}(b\sigma)$ is non-empty.  If $w$ is in the negative Weyl
chamber (i.e., $w = C_M^i$ for $i < 0$), then $X_w(1\sigma)$ is 
non-empty if and only
if $X_{b^{-1}w}(b\sigma)$ is non-empty.  
\end{prop}

\begin{proof}  We start by remarking that this is clear from the
results of Section~\ref{Identity} and Section~\ref{Non-identity}.  But we 
produce an {\em a priori} proof that may be generalizable in some respects to 
higher rank groups.

Assume that $w$ is in the positive Weyl chamber.  
If $X_w(1\sigma)$ is non-empty then let $x$ be such that
$\inv(x, \sigma(x)) = w$.  Let $A_1$ be an apartment containing $xC_M$ 
and $C_M$.  
Let $G$ be the minimal gallery (in $A_1$) between $C_M$ and $xC_M$, and
let $C_M = G_1 , G_2, \ldots, G_s = xC_M$ be the chambers of $G$.  Let $i$ be 
maximal such that $G_i$ is contained in $\mathcal{B}_1$.  We could have $G_i
= xC_M$ in the case that $w = C_M$.  Let $A_2$ be an apartment containing 
$C_M$ and $G_i$.  If $i \neq s$ then we also require that $A_2$ not contain
$G_{i+1}$.  Let $g \in SL_2(L)$ be such that $gC_M = C_M$ and $g$ sends 
$A_2$ to $A_M$.  To see that $\inv(gxC_M, b\sigma(gx)C_M) = bw$, 
consider Figure~\ref{Relationship1} in which $b = 
\left(
\begin{matrix}
	\pi^s & 0 \\
	0 & \pi^{-s}
\end{matrix}
\right).$
\begin{figure}
\input{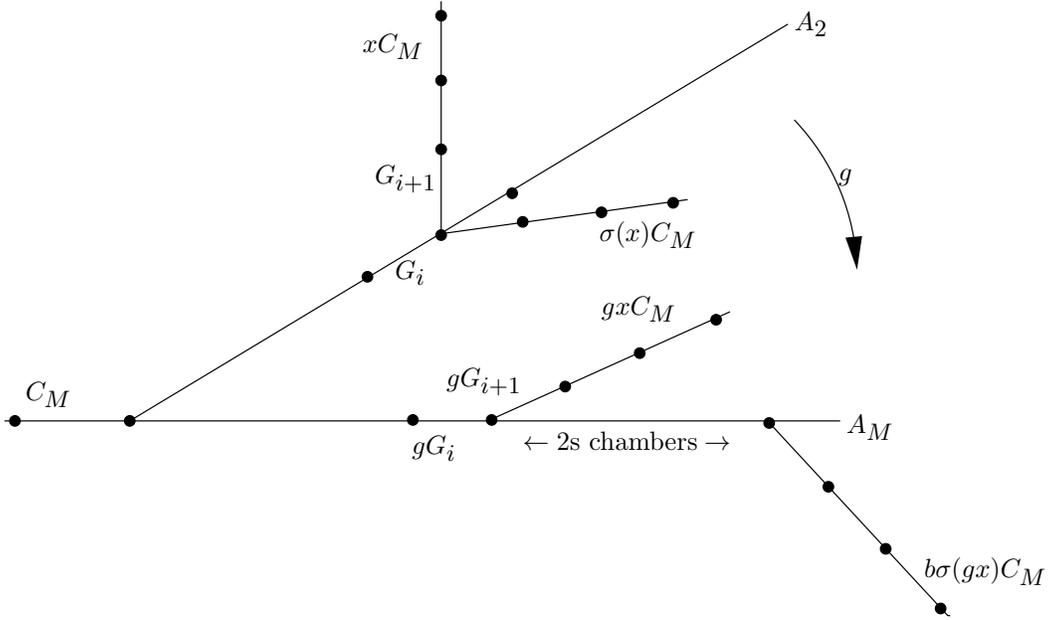}
\caption{For $w$ in the positive Weyl chamber}
\label{Relationship1}
\end{figure}
The minimal gallery $E_2$ between $gxC_M$ and $b\sigma(gx)C_M$
is $2s$ chambers longer than the minimal gallery $E_1$ between $xC_M$ 
and $\sigma(x)C_M$, and the
types of the vertices on the chamber $gxC_M$ are the same as those of the
corresponding vertices of $xC_M$.  So if $\rho(x^{-1}E_1)$ 
is the minimal
gallery in $A_M$ between $C_M$ and $w$ then 
$\rho((gx)^{-1}E_2)$ is the minimal gallery in $A_M$ between
$C_M$ and $bw$. To get the other implication, use $g^{-1}$.

Now assume that $w$ is in the negative Weyl chamber.  
Let $G$, $G_i$, $A_2$ and $g$ be as before.  In this situation,
Figure~\ref{Relationship1} becomes Figure~\ref{Relationship2}.
\begin{figure}
\centerline{\input{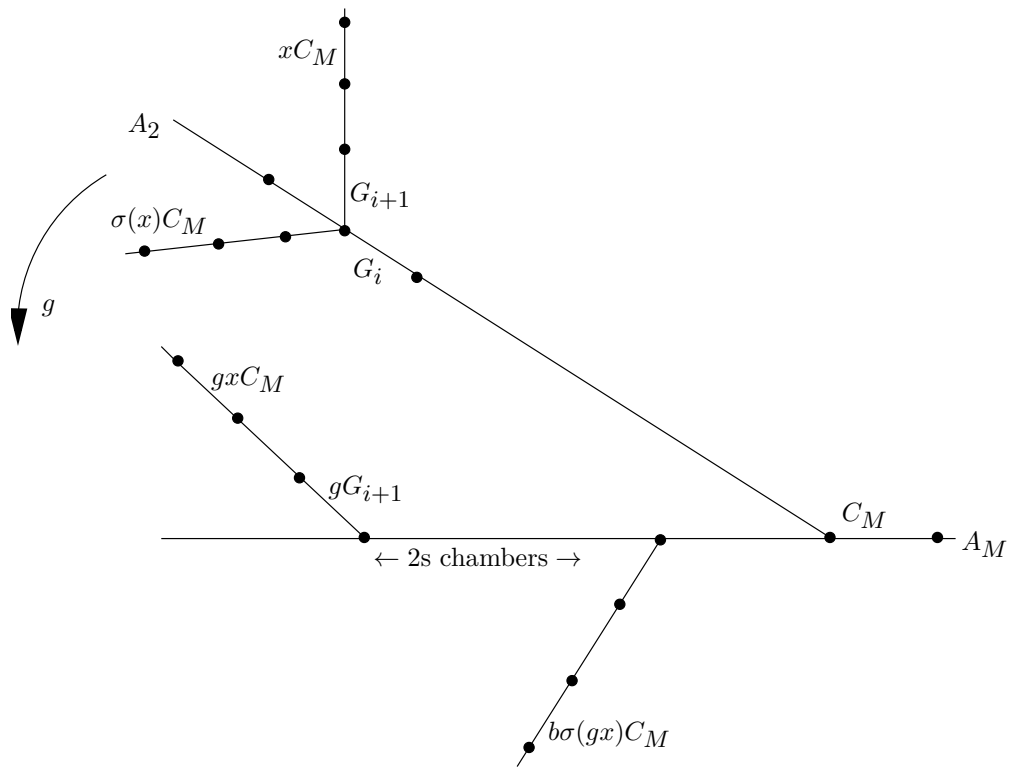}}
\caption{For $w$ in the negative Weyl chamber}
\label{Relationship2}
\end{figure}
As before, the gallery $E_2$ between $gxC_M$ and $b\sigma(gx)C_M$
is $2s$ chambers longer than the minimal gallery $E_1$ between $xC_M$ 
and $\sigma(x)C_M$, and the
types of the vertices on the chamber $gxC_M$ are the same as those of the
corresponding vertices of $xC_M$.  So if $\rho(x^{-1}E_1)$ 
is the minimal
gallery in $A_M$ between $C_M$ and $w$ then 
$\rho((gx)^{-1}E_2)$ is the minimal gallery in $A_M$ that begins 
at $C_M$ and goes in the same direction as $\rho(x^{-1}E_1)$,
but that is $2s$ chambers longer.  In the current case, the final chamber 
of this gallery is $\tilde{b}w$.

As before, to get the other implication, use $g^{-1}$.
\end{proof}

\section{$GL_2$ and $PGL_2$}\label{GL2}

The situation is slightly different for $GL_2$.  First, we must understand
differences in setup.  For $GL_2(L)$ it is possible to send $C_M$ to itself
without fixing it pointwise.  The matrix $m = \left(
\begin{matrix}
	0 & \pi \\
	1 & 0
\end{matrix} \right)$
does this, for example.  Another way of saying this is that 
contrary to the $SL_2$ situation, the Iwahori
associated to $C_M$ is no longer the set $\{ g : gC_M = C_M \}$.  

We still have $\inv: GL_2 (L)/I \times GL_2 (L)/I \rightarrow 
I \backslash GL_2(L) / I \simeq \tilde{W}$, but now $\tilde{W} \simeq W_a 
\rtimes \mathbb{Z}$, where $W_a$, the affine Weyl group, is in one-to-one
correspondence with the chambers of $A_M$ (recall that for $SL_2$, 
we had $W_a \simeq \tilde{W}$).  We therefore have 
$\inv: GL_2(L)/I \times GL_2(L)/I \rightarrow W_a \ltimes \mathbb{Z}$.  
This map
is given by $\inv(x,y) = (\rho(x^{-1}yC_M), v(\det(x^{-1}y)))$, 
where $v$ is the valuation.

The $\sigma$-conjugacy classes of $GL_2(L)$ are 
$$\left\{ \left(
\begin{matrix}
	\pi^{\alpha} & 0 \\
	0 & \pi^{\beta}
\end{matrix}
\right) : \alpha \geq \beta \right\} \cup \left\{ \left(
\begin{matrix}
	0 & \pi \\
	1 & 0
\end{matrix} \right)^{\alpha} : \alpha \in \mathbb{Z}, 
\alpha~\textrm{odd} \right\} = R.$$
This follows from a result of Kottwitz \cite{K1}.

For a given $\sigma$-conjugacy class $\{ x^{-1} b \sigma (x) : x \in GL_2 (L)
\}$, where $b \in R$, we are interested in which $(w, n) \in W_a \ltimes
\mathbb{Z}$ have non-empty $X_{(w,n)}(b\sigma)$.  To solve this, we will 
describe the set $\{ (\rho(x^{-1}b\sigma(x)C_M),
v(\det(x^{-1}b\sigma(x)))) \in W_a \ltimes \mathbb{Z} : x \in GL_2 (L) \}$
for any fixed $b$ in the above set $R$.  

First, note that $v(\det(x^{-1}b\sigma(x))) = v(\det(b))$, so the second
component is fixed for fixed $b$.  
The possible values of the first component can be 
determined using a process that is similar to that in Sections~\ref{Identity}
and~\ref{Non-identity}.  The difference is that $b$ can now take on more
values.  Also, it is {\em a priori} possible (although we will see that
it does not occur) for $\{ \rho(x^{-1}b\sigma(x)C_M) : x \in GL_2(L) \}$
to be bigger than $\{ \rho(x^{-1}b\sigma(x)C_M) : x \in SL_2(L) \}$,
since $GL_2(L)$ acts on $\mathcal{B}_{\infty}$ in ways that 
$SL_2(L)$ does not.

\begin{figure}
\centerline{\input{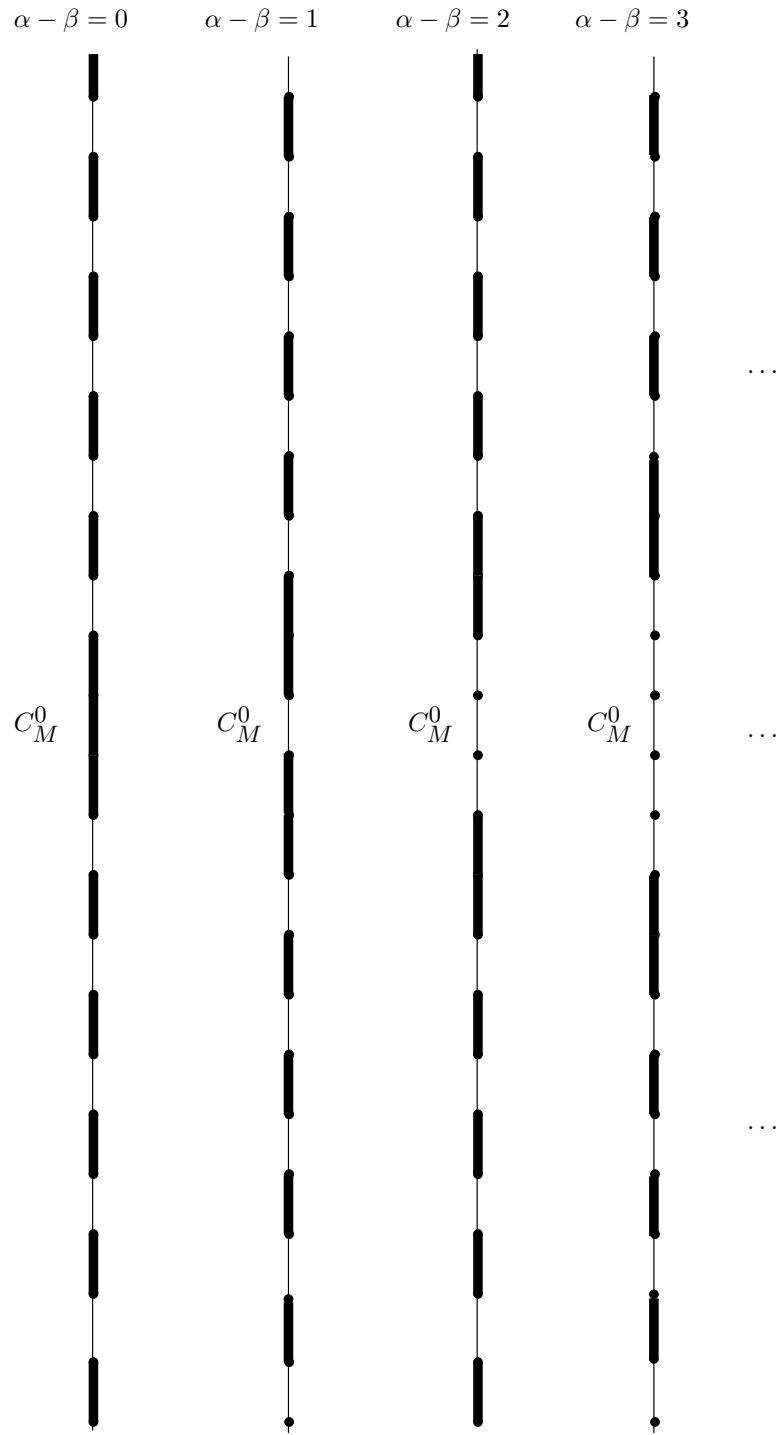}}
\caption{Summary of some results for $GL_2$}
\label{GL2Summary1}
\end{figure}

If $b = 
\left(
\begin{matrix}
	\pi^{\alpha} & 0 \\
	0 & \pi^{\beta}
\end{matrix}
\right)$
then $b$ shifts $A_M$ to the right by $\alpha - \beta$ units.  Using 
reasoning similar to that in Section~\ref{Non-identity}, for
these $b$ we see that $\rho(x^{-1}b\sigma (x)C_M)$ could be
any of the chambers $C_M^{\pm(\alpha - \beta)}$ or
$C_M^{\pm(\alpha - \beta + i)}$, for all $i$ odd, $i \geq 1$.  
Figure~\ref{GL2Summary1} summarizes these results.

Now we consider $b = \left(
\begin{matrix}
	0 & \pi \\
	1 & 0
\end{matrix} \right)^{\alpha}, \alpha \in \mathbb{Z}, \alpha~\textrm{odd}$.
The matrix $m = \left(
\begin{matrix}
	0 & \pi \\
	1 & 0
\end{matrix}
\right)$ flips $A_M$ about the center of $C_M$.  Therefore so does the matrix
$b$.  Let $G$ be the minimal gallery between 
$C_M$ and $xC_M$.  Then $b\sigma(G)$
is the minimal gallery between $C_M$ and $b\sigma(x)C_M$.  Let $C_M = G_1,
G_2, \ldots, G_s = xC_M$ be the chambers of $G$.  Then $G_s, G_{s-1}, \ldots,
G_1, b\sigma(G_2), b\sigma(G_3), \ldots, b\sigma(G_s)$ is the unique
minimal gallery between $xC_M$ and $b\sigma(x)C_M$.  Call this gallery
$\Gamma_x$.  It has length $2s-1$.  

\begin{figure}
\centerline{\input{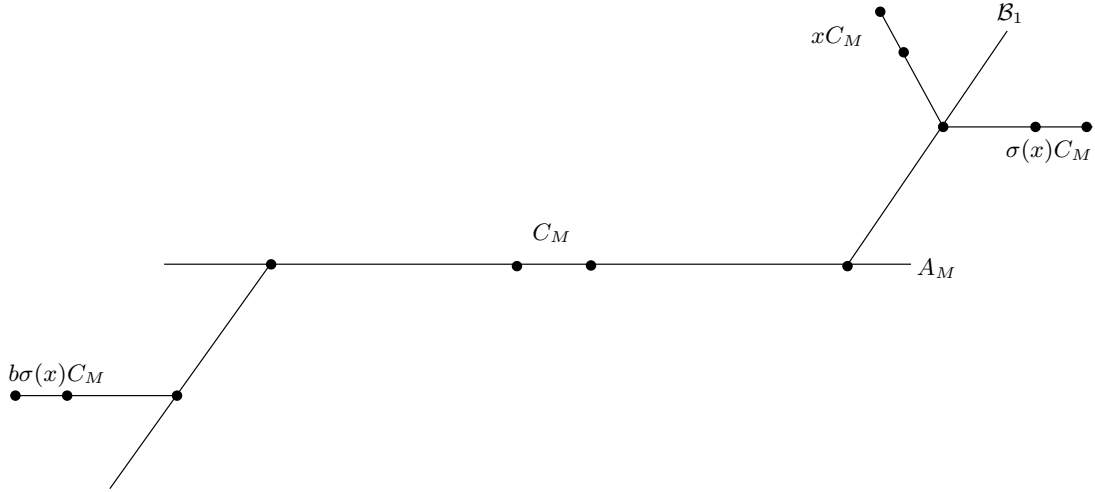}}
\caption{Computations for some $\sigma$-conjugacy classes for $GL_2$}
\label{AlphaProofHelp}
\end{figure}

\begin{figure}
\centerline{\input{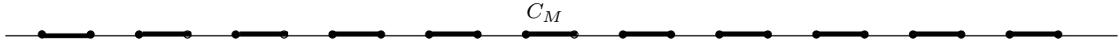}}
\caption{Summary of some more results for $GL_2$}
\label{Summ}
\end{figure}

By choosing $x$ appropriately, we 
can arrange for $s$ to be any integer $\geq 1$.  Further, if 
$\rho(x^{-1}\Gamma_x)$ has final chamber $C_M^i$, 
then $\rho(y^{-1} \Gamma_y)$ has final chamber $C_M^{-i}$, 
where $y = xm$.  Therefore, the possible values for
$\rho(x^{-1}b\sigma(x)C_M)$ are $C_M^i$ for any $i$ even.
This is independent of which odd $\alpha \in \mathbb{Z}$ we have chosen.
Figure~\ref{AlphaProofHelp} makes the calculations more clear, 
and Figure~\ref{Summ} summarizes
the possible values of the first component of $\inv$.  
The second component is
always $v(\det(x^{-1}b\sigma(x))) = v(\det(b))$.

For $PGL_2$, $\tilde{W} \cong W_a \ltimes \mathbb{Z}/2$, and 
$\inv(x,y) = (\rho(x^{-1}yC_M), v(\det(x^{-1}y)))$,
where in $PGL_2$, $v(\det(x^{-1}y))$ is only determined mod $2$.
The results for $PGL_2$ are the same as those for $GL_2$, but second-component
values are computed mod $2$.

\section{Reorganization of Results}

So far, we have stated all results by specifying, for each $b$, which $w$ 
have non-empty $X_w(b\sigma)$, or which $(w,n)$ have non-empty 
$X_{(w,n)}(b\sigma)$.  We reword these results by saying, for each $w$,
which $b$ have non-empty $X_w(b\sigma)$, and for each $(w,n)$, which
$b$ have non-empty $X_{(w,n)}(b\sigma)$.

We first reword results for $SL_2$.  If $b = \left(
\begin{matrix}
	\pi^s & 0 \\
	0 & \pi^{-s}
\end{matrix} \right)$, then let $\Sigma_s$ be 
the $\sigma$-conjugacy class of $b$.
Let $D_M^i$ denote the double-$I$-coset corresponding to the chamber
$C_M^i$.  Then one can see that $D_M^{\pm 2n} \subseteq \Sigma_n$ for
$n \in \mathbb{Z}$, $n \geq 0$. One can also see that for $n \geq 0$, 
$n \in \mathbb{Z}$, $D_M^{\pm(2n+1)}$ intersects the $\sigma$-conjugacy
classes $\Sigma_0, \Sigma_1, \ldots, \Sigma_n$, and no others.  
These results are summarized in Tables~\ref{ReverseResults1_SL2} 
and~\ref{ReverseResults2_SL2}.  
\begin{table}
\centerline{
\begin{tabular}{|c|c|} 
\hline
i & j \\
\hline
\hline
0 & 0 \\
\hline
2 & 1 \\
\hline
4 & 2 \\
\hline
6 & 3 \\
\hline
$\vdots$ & $\vdots$ \\
\hline
\end{tabular}}
\caption{$D_M^{\pm i}$ intersects exactly these $\Sigma_j$ non-trivially 
for $i$ even}
\label{ReverseResults1_SL2}
\end{table}
\begin{table}
\centerline{
\begin{tabular}{|c|c|}
\hline 
$i$ & $j$ \\
\hline 
\hline 
$1$ & $0$ \\
\hline 
$3$ & $0,1$ \\
\hline 
$5$ & $0,1,2$ \\
\hline 
$7$ & $0,1,2,3$ \\
\hline \
$\vdots$ & $\vdots$ \\
\hline
\end{tabular}}
\caption{$D_M^{\pm i}$ intersects exactly these $\Sigma_j$ non-trivially 
for $i$ odd}
\label{ReverseResults2_SL2}
\end{table}

In particular, $D_M^{2n}$ for 
$n \in \mathbb{Z}$ and $D_M^{\pm 1}$ are each
completely contained in some $\sigma$-conjugacy class.  The cosets
$D_M^{\pm (2n + 1)}$ for $n \in \mathbb{Z}$, $n \geq 0$ are spread 
over increasingly many $\sigma$-conjugacy classes as $n$ gets larger.

We now reword results for $GL_2$.  Let $\Sigma_{\alpha , \beta}$ be the
$\sigma$-conjugacy class of 
$\left(
\begin{matrix}
	\pi^{\alpha} & 0 \\
	0 & \pi^{\beta}
\end{matrix}
\right)$, and let $\Sigma_{\alpha}$ be the $\sigma$-conjugacy class of
$\left(
\begin{matrix}
	0 & \pi \\
	1 & 0
\end{matrix}
\right)^{\alpha}$.  Let $D_M^i$ denote the element of $W_a$ corresponding
to the chamber $C_M^i$ in $A_M$, and let $(D_M^i , j)$ be an element of 
$\tilde{W} \cong W_a \ltimes \mathbb{Z}$.  We also denote the corresponding
double-$I$-coset by $(D_M^i , j)$.  Then one can see that 
$(D_M^i , j) \subseteq
\Sigma_{\frac{j+i}{2},\frac{j-i}{2}}$ if $i \equiv j \mod 2$. If $i$ is even
and $j$ is odd, then $(D_M^i , j)$ intersects $\Sigma_j$ and 
$\Sigma_{\frac{j+k}{2},\frac{j-k}{2}}$ 
for $1 \leq k \leq i-1$ an odd integer.  It does not intersect non-trivially
with any other $\sigma$-conjugacy classes.  Finally, if $i$ is 
odd and $j$ is even, then $(D_M^i , j)$ intersects 
$\Sigma_{\frac{j+k}{2},\frac{j-k}{2}}$ for 
$0 \leq k \leq i-1$ an even integer.  It does not intersect 
non-trivially with any other $\sigma$-conjugacy classes.  These results are
summarized in Table~\ref{ReverseResults_GL2}.
 
\begin{table}
\centerline{
\begin{tabular}{|c||c|c|} 
\hline
$i$ & $j$ odd & $j$ even \\
\hline
\hline
$0$ & $\Sigma_j$ & $\Sigma_{\frac{j}{2},\frac{j}{2}}$ \\
\hline
$\pm 1$ & $\Sigma_{\frac{j+1}{2},\frac{j-1}{2}}$ & 
	$\Sigma_{\frac{j}{2},\frac{j}{2}}$ \\
\hline
$\pm 2$ & $\Sigma_j$ ; $\Sigma_{\frac{j+1}{2},\frac{j-1}{2}}$ &
	$\Sigma_{\frac{j+2}{2},\frac{j-2}{2}}$ \\
\hline
$\pm 3$ & $\Sigma_{\frac{j+3}{2},\frac{j-3}{2}}$ &
	$\Sigma_{\frac{j}{2},\frac{j}{2}}$ ; 
		$\Sigma_{\frac{j+2}{2},\frac{j-2}{2}}$ \\
\hline
$\pm 4$ & $\Sigma_j$ ; $\Sigma_{\frac{j+1}{2},\frac{j-1}{2}}$ ; 
		$\Sigma_{\frac{j+3}{2},\frac{j-3}{2}}$ & 
	$\Sigma_{\frac{j+4}{2},\frac{j-4}{2}}$ \\
\hline
$\pm 5$ & $\Sigma_{\frac{j+5}{2},\frac{j-5}{2}}$ &
	$\Sigma_{\frac{j}{2},\frac{j}{2}}$ ; 
		$\Sigma_{\frac{j+2}{2},\frac{j-2}{2}}$ ; 
		$\Sigma_{\frac{j+4}{2},\frac{j-4}{2}}$ \\
\hline
$\pm 6$ & $\Sigma_j$ ; $\Sigma_{\frac{j+1}{2},\frac{j-1}{2}}$ ; 
		$\Sigma_{\frac{j+3}{2},\frac{j-3}{2}}$ ; 
		$\Sigma_{\frac{j+5}{2},\frac{j-5}{2}}$ &
	$\Sigma_{\frac{j+6}{2},\frac{j-6}{2}}$ \\
\hline
$\vdots$ & $\vdots$ & $\vdots$ \\
\hline
\end{tabular}}
\caption{$\sigma$-conjugacy classes that intersect $(D_M^i , j)$ 
non-trivially}
\label{ReverseResults_GL2}
\end{table}

\chapter{Groups of Semisimple Rank $2$}

In Section~\ref{SL3} we will show for $b$ an element of a certain 
collection of $\sigma$-conjugacy class representatives of $SL_3(L)$,
exactly which double-$I$-cosets intersect $\{ x^{-1} b \sigma(x) : 
x \in SL_3(L) \}$ non-trivially.  
In Section~\ref{GL3}, we will show that these results
can also be made to give complete information of the same kind for certain
$\sigma$-conjugacy classes of $GL_3(L)$ and $PGL_3(L)$.  
Section~\ref{SP4} attempts to apply the same methods to $Sp_4$.
Certain results are conjectural in this setting because some of the
lengthy computations done for $SL_3$ have not been done in entirety
for $Sp_4$.  Section~\ref{GSP4} discusses how the methods of
Section~\ref{SL3} could be applied to other rank $2$ groups.
Section~\ref{Invariance} gives some invariance properties of the set
$\{ \inv(x,b\sigma(x)) : x \in G(L) \}$ that hold for any
simply-connected group $G$, and any $\sigma$-conjugacy class $b$.
These could be applied to higher rank groups, for instance.
Some of the methods of Section~\ref{SL3} could also be applied
to $SL_n$ or $Sp_{2n}$.  A comment to this effect is made when this is the 
case.

\section{$SL_3$}\label{SL3}

In this section we enumerate the non-empty $X_w(b\sigma)$ for $SL_3$, where
the letter $b$ will always refer to an element of $SL_3(L)$ of the form 
$$\left( 
\begin{matrix}
	\pi^{\alpha} & 0 & 0 \\
	0 & \pi^{\beta} & 0 \\
	0 & 0 & \pi^{-\alpha - \beta}
\end{matrix} 
\right),
$$
with $\alpha \geq \beta \geq \alpha - \beta$.  These are
representatives of distinct $\sigma$-conjugacy classes of $SL_3(L)$.
However, not every $\sigma$-conjugacy class is represented 
by one of these elements \cite{K1}.

The collection of chambers in $A_M$ that correspond to $w$ with 
non-empty $X_w(b\sigma)$ is computed in the following way.  
In Section~\ref{SMGandComp}, we give some necessary definitions.  In 
Section~\ref{Superset}, we use geometric methods to produce a superset of the
desired collection of chambers.  In Section~\ref{Rapoport} we use
algebraic methods that produce a subset of the desired collection for 
$b = 1$.  We
enlarge this subset in Section~\ref{Subset}.  This enlargement will be
equal to the superset from Section~\ref{Superset}, and therefore is 
the solution set of chambers.  
In Section~\ref{SubsetGeometric} we give a geometric method of arriving 
at the subset of Section~\ref{Subset} that also applies for $b \neq 1$.  In 
Section~\ref{Symmetry3}, we prove a symmetry result analogous to that done for 
$SL_2$ in Section~\ref{Symmetry}.  In Section~\ref{Relationship3}, we
develop a relationship between the collection of $w$ with non-empty
$X_w (1\sigma)$ and the collection of $w$ with non-empty $X_w (b\sigma)$.
This is done by illuminating an efficient method of computing
$\{ \inv (x,b\sigma(x)) : x \in SL_3(L) \}$ that works for any of the 
$b$ listed above.

\subsection{Standard Minimal Galleries and 
Composite Galleries}\label{SMGandComp}

Our goal is to determine which $X_w(b\sigma)$ are non-empty. The
technique that led to success in the rank $1$ case involved considering,
for every chamber $E \subseteq \mathcal{B}_{\infty}$, a gallery $\Gamma_{E}$
connecting $E$ to $b\sigma(E)$.  We used the unique minimal gallery
between these two chambers.  We will use a similar process for $SL_3$.
However, there is no longer a unique minimal gallery
$\Gamma_{E}$ connecting $E$ and $b\sigma(E)$. In this section we will 
specify a choice of gallery, which we will call the {\em composite gallery}.
The composite gallery is usually not minimal.

\begin{figure}
\setlength{\unitlength}{1in}
\begin{picture}(6,4.75)(0,0)
\centerline{\includegraphics[height=4.75in]{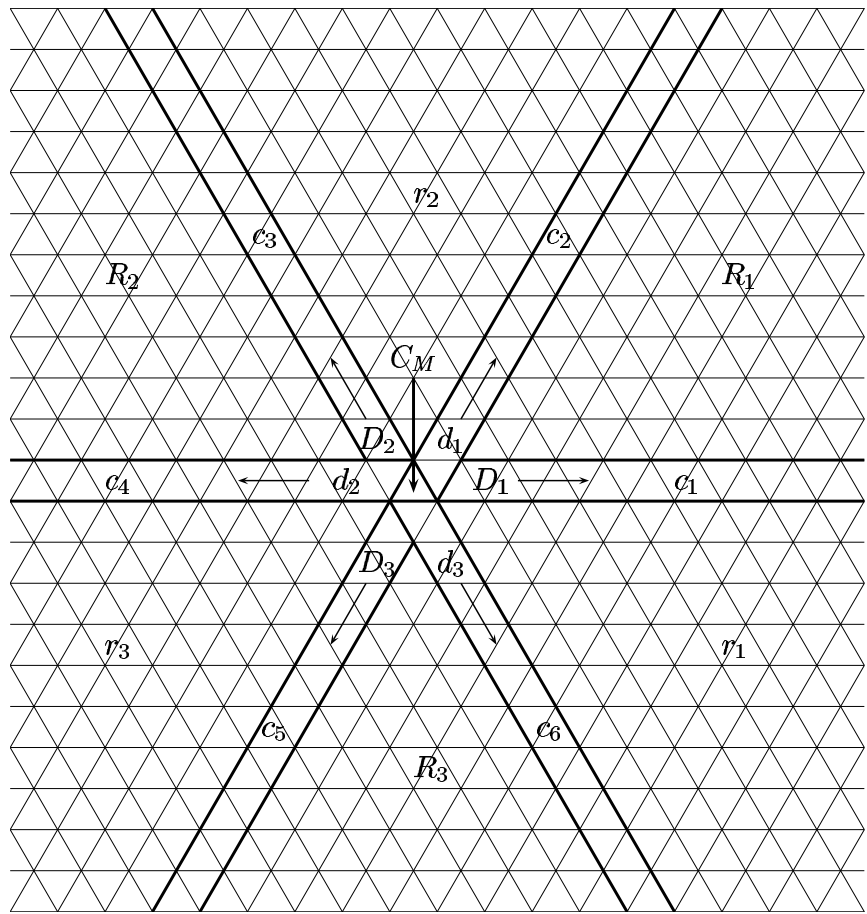}}
\end{picture}
\caption{Primary and secondary directions}
\label{directions}
\end{figure}

First define the three {\em primary directions} $D_1,D_2,D_3$ and the
three {\em secondary directions} $d_1,d_2,d_3$ in the main apartment
$A_M$ for $SL_3$ as marked by arrows in Figure~\ref{directions}.
Given a chamber
$E \subseteq A_M$, we define the {\em standard minimal gallery} (SMG) between
$C_M$ and $E$ as follows.  If $E$ is in one of the corridors marked
$c_1,\ldots,c_6$ in Figure~\ref{directions}, 
then the SMG is the unique minimal gallery
from $C_M$ to $E$.  If $E$ is in region $R_i$, proceed first in the
direction $D_i$, then in the direction $d_i$.  If $E$ is in region $r_i$,
proceed first in direction $D_i$, then in direction $d_{j}$, where $j
\equiv i-1 \mod 3$.

The following lemma allows us to define the SMG between $C_M$ and $E$ for
any chamber $E \subseteq \mathcal{B}_{\infty}$ in a natural way.  As in 
Chapter~\ref{RankOne}, $\rho$ is the 
retraction from $\mathcal{B}_{\infty}$ to $A_M$ centered at $C_M$.

\begin{lem} 
There is a unique gallery $G_E$ between $C_M$ and $E$
which is minimal, and such that $\rho (G_E)$ is the SMG
between $C_M$ and $\rho (E)$. 
\label{Moon}
\end{lem}

\begin{proof}
Let $A_1$ be an apartment containing $C_M$ and $E$.  Let $g_1 \in SL_3 (L)$
be such that $g_1 A_1 = A_M$ and $g_1 C_M = C_M$.  If $G_{\rho(E)}$ is the SMG
between $C_M$ and $\rho (E)$, then let $G_E = g_1^{-1}G_{\rho(E)}$.  This
proves existence.

To prove uniqueness, note that if $\tilde{G}_E$ is another minimal gallery
from $C_M$ to $E$, then $\tilde{G}_E \subseteq A_1$ (see 
Definition~\ref{parallel} and Theorem~\ref{paralleltheorem} 
in Section~\ref{SubsetGeometric}). If 
$\rho(\tilde{G}_E) = G_{\rho(E)}$ then $g_1^{-1}G_{\rho(E)} 
= \tilde{G}_E$.  But $g_1^{-1}G_{\rho(E)} = G_E$.
\end{proof}

To obtain a composite gallery, $\Gamma_E$, from the SMG 
between $C_M$ and $E$, proceed as follows.

\noindent {\em Case 1:} $b = 1$.  Let $\Gamma_{E}^1$ be the gallery
composed of the chambers of the SMG between $C_M$ and $E$ that are not in
$\mathcal{B}_1$.  Let $\Gamma_{E}^3 = \sigma(\Gamma_{E}^1)$.  Then
$\Gamma_E = \Gamma_{E}^1 \cup \Gamma_{E}^3$.  We will call the edge $e$
between $\Gamma_{E}^1$ and $\Gamma_{E}^3$ the {\em edge of departure of the
SMG between $C_M$ and $E$ from $\mathcal{B}_1$}.  

\noindent {\em Case 2:} $b \neq 1$.  Let $\Gamma_{E}^1$ be the gallery
composed of the chambers of the SMG between $C_M$ and $E$ that are not in
$A_M$. Let $\Gamma_{E}^3 = b\sigma(\Gamma_{E}^1)$.  Let $e$ be the unique
edge in $\Gamma_{E}^1$ that is contained in $A_M$.  Let $\Gamma_{E}^2$ be
any minimal gallery connecting $e$ and $b\sigma (e) = be$.  Then $\Gamma_E
= \Gamma_{E}^1 \cup \Gamma_{E}^2 \cup \Gamma_{E}^3$.  We call $e$ the 
{\em edge of departure of the SMG between $C_M$ and $D$ from $A_M$}.
The term {\em edge of departure} without further modification will be used
to refer either to an edge of departure from $A_M$ or to an edge of 
departure from $\mathcal{B}_1$.

We will also need the following definitions.

\begin{defn} If $E_1$ and $E_2$ are two chambers in $\mathcal{B}_{\infty}$
that share an edge $e$, then we say the {\em transition type} from $E_1$
to $E_2$ is the type of the two vertices not in $e$.
\end{defn}

\begin{defn}
If $G$ is a non-stuttering gallery in $\mathcal{B}_{\infty}$ 
consisting of $G_0, \ldots,
G_s$, and if the transition type from $G_{i-1}$ to $G_i$ is $t_i$ for
$1 \leq i \leq s$, then we say that {\em $G$ has type $t_1, \ldots ,
t_{s}$}.
\end{defn}

Note that if $G \subseteq A_M$, then $G_0$ and $t_1, \ldots, t_{s}$
determine $G$.  It is also possible to specify a gallery type by giving
a picture of a non-stuttering gallery in $A_M$.

\subsection{A Superset of the Solution Set}\label{Superset}

Just as in the $SL_2$ case, $SL_3$ is simply connected.  So given $b$, an
answer to the question of which $w \in \tilde{W}$ have non-empty
$X_w(b\sigma)$ can be given by specifying the corresponding chambers in
$A_M$.  As such, the answer we are looking for is the
set of chambers $S = \{ \rho (x^{-1}b\sigma(x)C_M) : x \in
SL_3(L) \}$.  In the $SL_2$ case, for every $x \in SL_2(L)$, we had a
unique minimal gallery $\Gamma_{x}$ connecting $xC_M$ and
$\sigma(x)C_M$.  Since this gallery was minimal, we were able to determine
$\rho (x^{-1}b\sigma(x)C_M)$ just by folding $x^{-1}
\Gamma_{x}$ from $C_M$ down into $A_M$ in the unique possible way.  In
the $SL_3$ case, we let $\Gamma_{x}$ instead be the composite gallery
$\Gamma_E$ (for $E = xC_M$) 
constructed in the previous section.  Since this composite gallery is not
necessarily minimal, one can not determine
$\rho (x^{-1}b\sigma(x)C_M)$ from its type and starting point alone.  One can,
however, get a set of possible chambers $S_x$ for
$\rho (x^{-1}b\sigma(x)C_M)$ by enumerating all the possible
foldings of a gallery of the same type as $x^{-1}\Gamma_{x}$, and
putting the final chamber of each into $S_x$. The set $S_1 = \cup_{x \in
SL_3(L)} S_x$ contains the set $S$ in which we are interested.

\begin{figure}
\setlength{\unitlength}{1in}
\begin{picture}(4,2.5)(0,0)
\centerline{\includegraphics[height=2.5in]{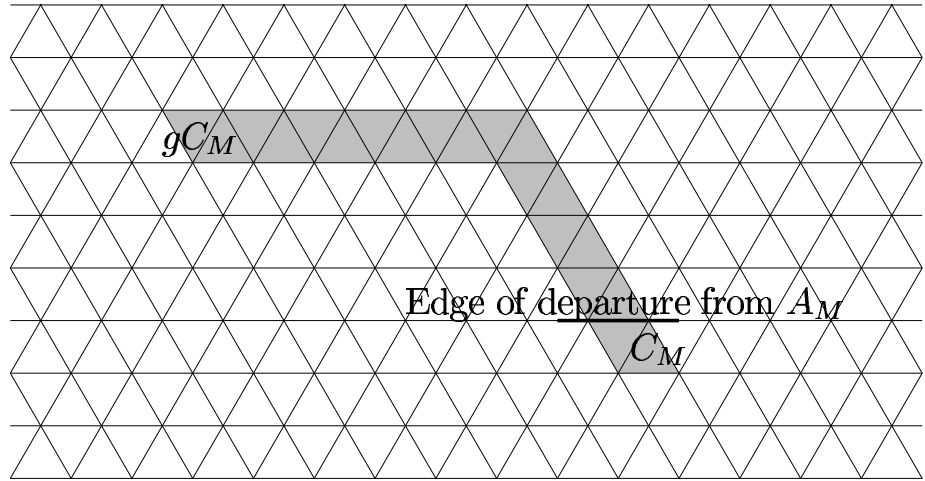}}
\end{picture}
\caption{Example SMG type and edge of departure}
\label{GammaD1}
\end{figure}

\begin{figure}
\setlength{\unitlength}{1in}
\begin{picture}(4.75,3.25)(0,0)
\centerline{\includegraphics[height=3.25in]{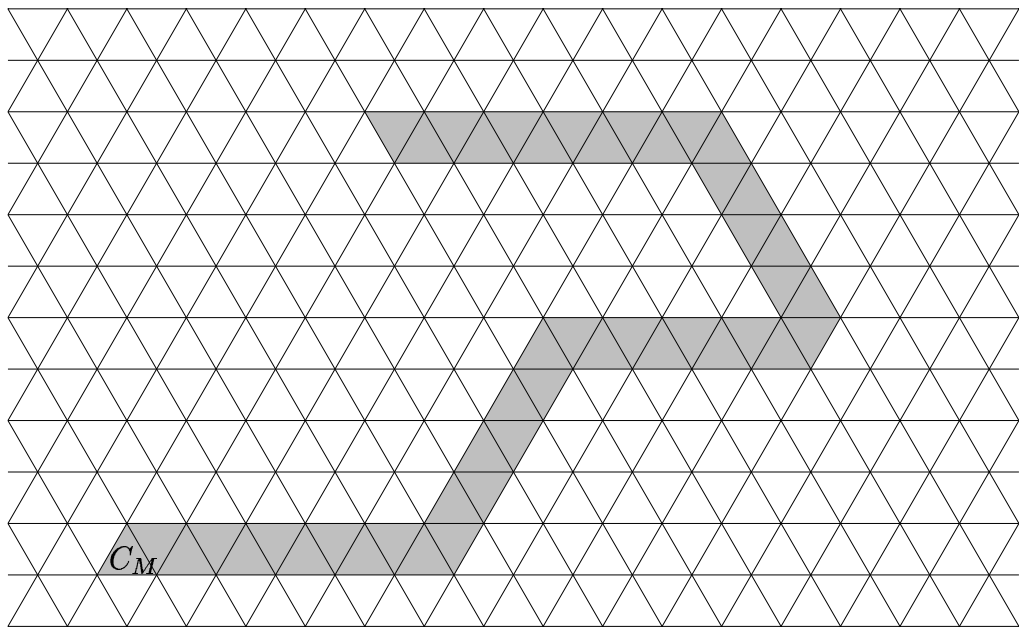}}
\end{picture}
\caption{Composite gallery example}
\label{CompositeGalleryExample}
\end{figure}

\begin{figure}
\setlength{\unitlength}{1in}
\begin{picture}(5,3.5)(0,0)
\centerline{\includegraphics[height=3.5in]{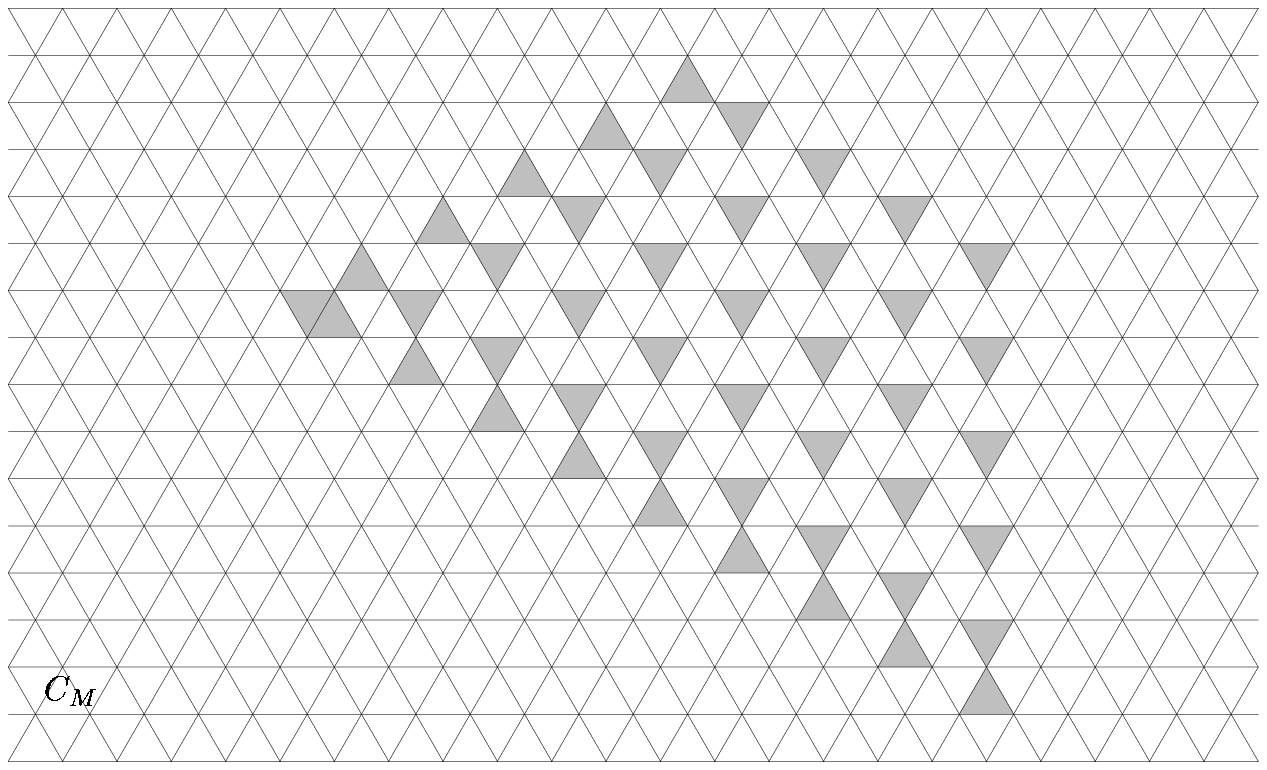}}
\end{picture}
\caption{Results from the SMG in Figure~\ref{GammaD1}}
\label{SingleExampleResults}
\end{figure}

We include an example of this computational process.  Let 
$$b = \left( 
\begin{matrix}
	\pi^3 & 0 & 0 \\
	0 & \pi^{-1} & 0 \\
	0 & 0 & \pi^{-2}
\end{matrix}
\right),$$
and let $x$ be such that the SMG connecting $C_M$ to $xC_M$ has type and
edge of departure from $A_M$ as indicated in Figure~\ref{GammaD1}.  
Then the resulting composite gallery $\Gamma_x$ 
has the same type as the gallery pictured in
Figure~\ref{CompositeGalleryExample}.  
Therefore, 
$\rho (x^{-1}b\sigma(x)C_M)$ must be one of the chambers marked
on Figure~\ref{SingleExampleResults}.
These results are achieved by computing all possible values of 
$\rho(E)$, where $E$ is the last chamber of some gallery that
begins at $C_M$ and that has the same type as $\Gamma_x$ (pictured in 
Figure~\ref{CompositeGalleryExample}).

The apparent problem with this process is that it seems to be an infinite
computation.  The set $S_1$ is a union of the $S_x$, where $x \in SL_3(L)$
is arbitrary.  Of course, we need only consider one representative $x$
from each coset $SL_3(L)/I$, but there are still infinitely many such
cosets.  We can further optimize using the fact that only the {\em type} of
$x^{-1}\Gamma_{x}$ is important for computing $S_x$.  This
type is determined by the type of the SMG between $C_M$ and $xC_M$, and
the departure edge of this SMG from $A_M$ (or from $\mathcal{B}_1$ if $b =
1$).  We call a pair consisting of an SMG type and a departure edge
a {\em type-edge pair}.  The set we are trying to compute is 
$S_1 = \cup S_{(t,e)}$, where the union is over all type-edge pairs $(t,e)$,
and where $S_{(t,e)} = S_x$ for some $x \in SL_3(L)$ such that the SMG
from $C_M$ to $xC_M$ has type $t$ and departure edge $e$.

The benefit of this point of view is that the $S_{(t,e)}$ can be separated 
into finitely many infinite classes, each of which can be computed all 
at once.  We give an example of two of these 
infinite classes, $I_1$ and $I_2$, for $b$ as above.  
The SMG types in $I_1$ are those represented by the SMGs 
of the chambers in the region
$R_2$ in Figure~\ref{directions}.  The edges of departure from $A_M$ that 
we will consider are
the horizontal ones.  The SMG types in $I_2$ are those represented
by the SMGs of the chambers in corridor $c_4$.  The pairs $(t,e) \in I_2$ 
will have arbitrary $e$.

\begin{figure}
\centerline{\input{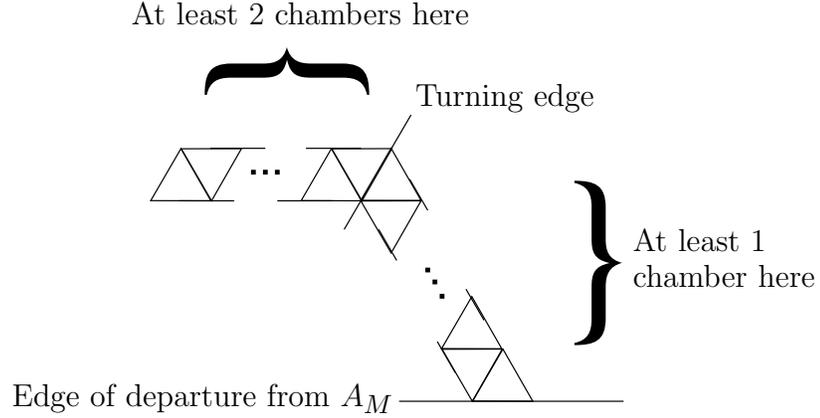}}
\caption{$\Gamma_{E}^1$ possibilities for $I_1$}
\label{GammaE1InfiniteClass}
\end{figure}

We consider $I_1$ first.  As in Case 2 of the definition of the 
composite gallery, let
$\Gamma_{E}^1$ be the gallery composed of the chambers of the SMG in
question that are not in $A_M$ (i.e., those after the edge of departure).  
So $I_1$ gives rise to the
$\Gamma_{E}^1$ in Figure~\ref{GammaE1InfiniteClass}.  

Let $W$ be the finite Weyl group, and let $w \in W$, $a \in T(F)$ be
such that $awC_M = \rho(E)$.  We break $I_1$ into six sub-classes
according to $w$.  Let 
$$f = \left( 
\begin{matrix}
	-1 & 0 & 0 \\
	0 & 0 & 1 \\
	0 & 1 & 0
\end{matrix} \right).$$

\begin{figure}
\centerline{\input{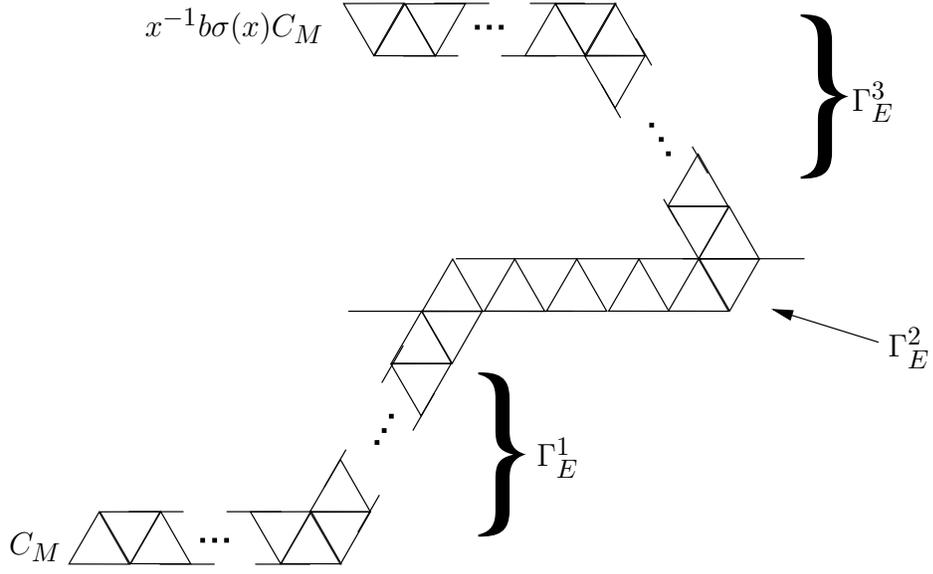}}
\caption{Some composite gallery possibilities for $I_1$}
\label{CompositeGalleryI1_1}
\end{figure}

\begin{figure}
\setlength{\unitlength}{1in}
\begin{picture}(6,7.5)(0,0)
\centerline{\includegraphics[height=7.5in]{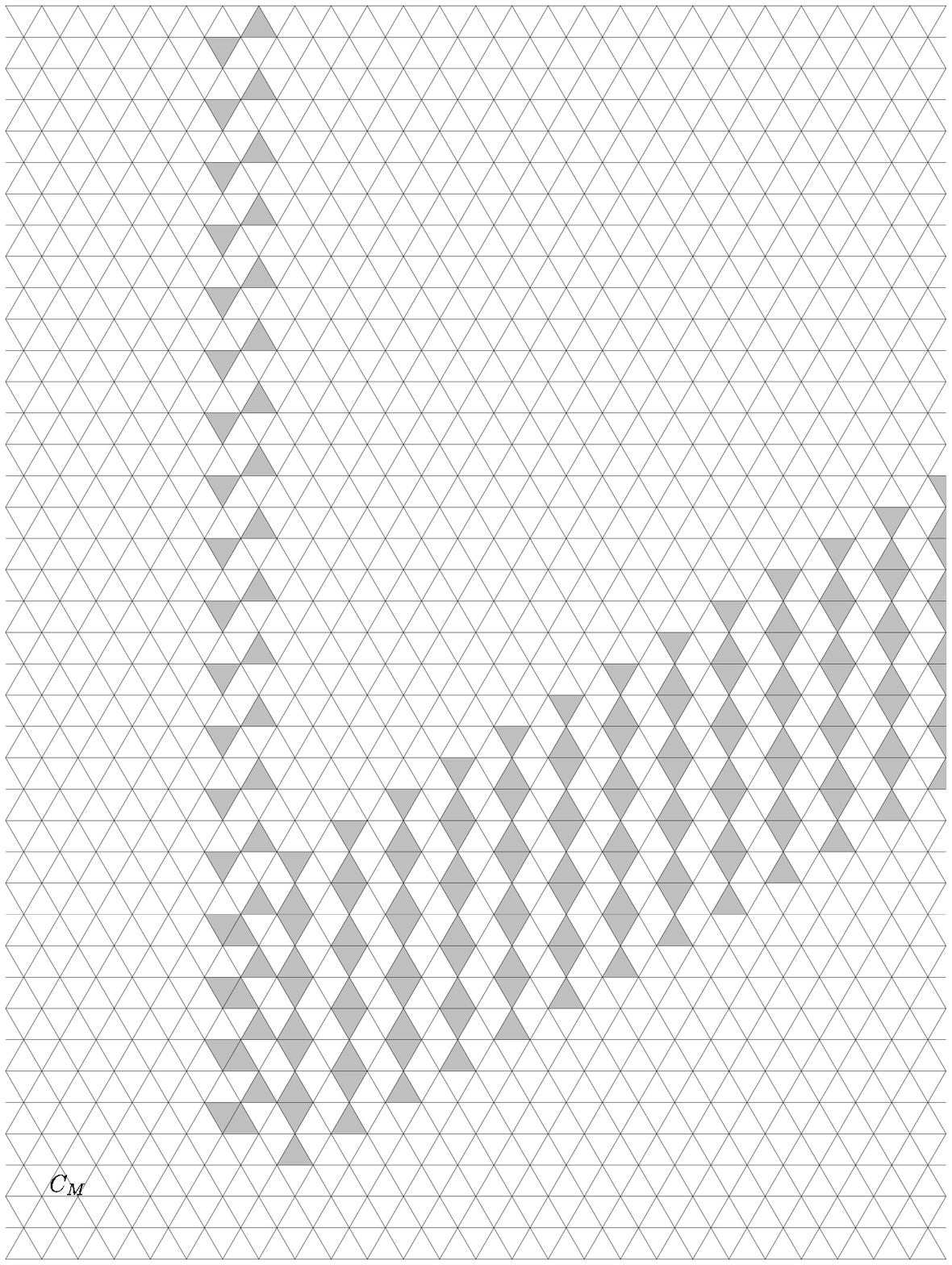}}
\end{picture}
\caption{Class $I_1$, where $\Gamma_{E}^1$ has $3$ chambers after the 
turning edge}
\label{ThreeAfterTE}
\end{figure}

\begin{figure}
\setlength{\unitlength}{1in}
\begin{picture}(6,7.5)(0,0)
\centerline{\includegraphics[height=7.5in]{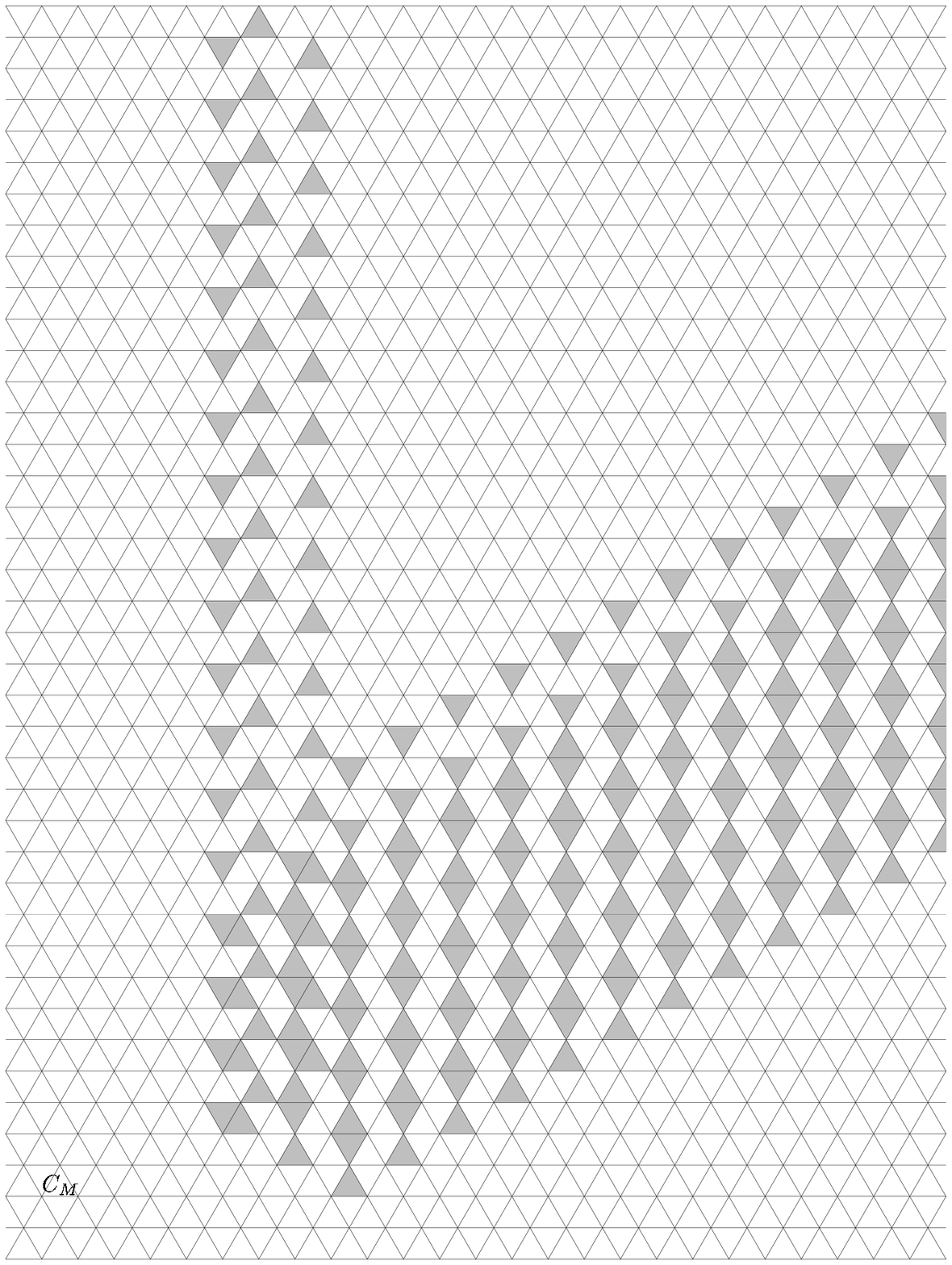}}
\end{picture}
\caption{Class $I_1$, where $\Gamma_{E}^1$ has $5$ chambers after the 
turning edge}
\label{FiveAfterTE}
\end{figure}

\begin{figure}
\setlength{\unitlength}{1in}
\begin{picture}(6,7.5)(0,0)
\centerline{\includegraphics[height=7.5in]{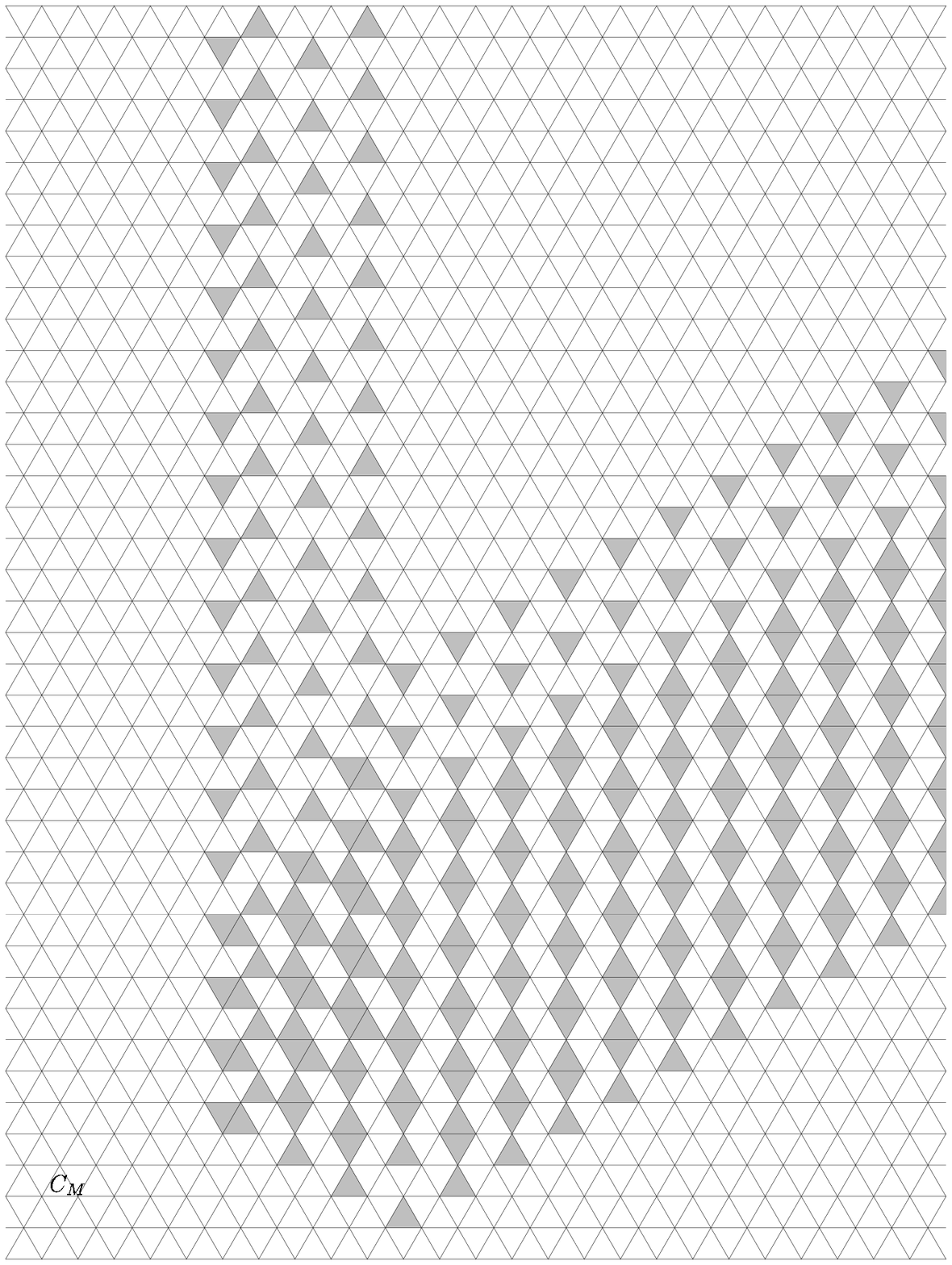}}
\end{picture}
\caption{Class $I_1$, where $\Gamma_{E}^1$ has $7$ chambers after the 
turning edge}
\label{SevenAfterTE}
\end{figure}

\begin{figure}
\setlength{\unitlength}{1in}
\begin{picture}(6,7.5)(0,0)
\centerline{\includegraphics[height=7.5in]{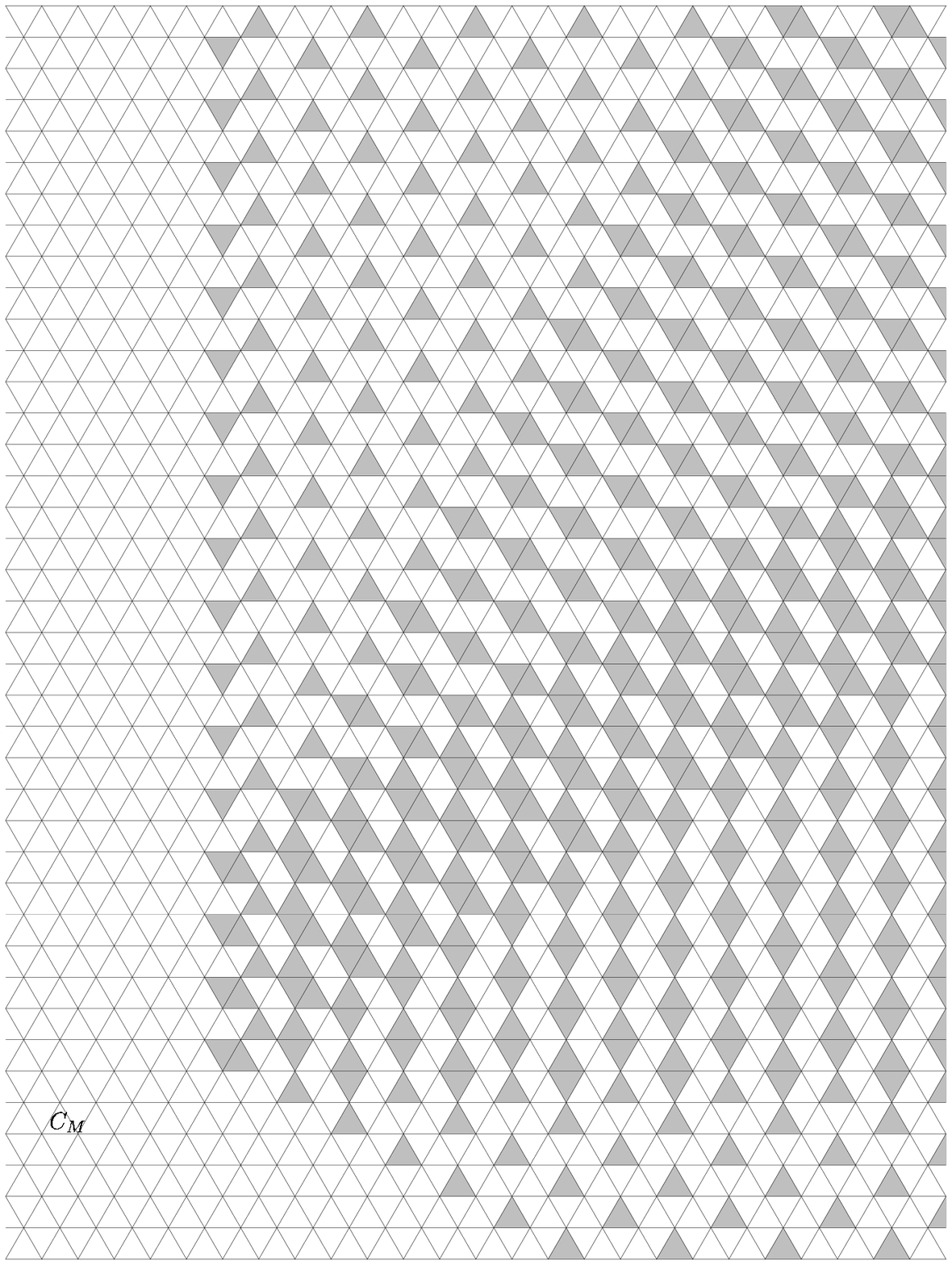}}
\end{picture}
\caption{Summary of results for $I_1$, an odd number of 
chambers after the turning edge}
\label{OddAfterTE}
\end{figure}

\noindent {\em Case 1:} $w = f$. We get a composite 
gallery of the type shown in
Figure~\ref{CompositeGalleryI1_1}.
If the {\em turning edge} is
labeled in Figure~\ref{GammaE1InfiniteClass}, 
then we consider in Figure~\ref{ThreeAfterTE} the instances where
$\Gamma_{E}^1$ has $3$ chambers after the turning edge.  
We consider in
Figure~\ref{FiveAfterTE} the instances where $\Gamma_{E}^1$ 
has $5$ chambers after the
turning edge.  
We consider on Figure~\ref{SevenAfterTE} the instances where
$\Gamma_{E}^1$ has $7$ chambers after the turning edge.  
After doing
these computations, it is easy to see what the situation is for the
instances where $\Gamma_{E}^1$ has $2n+1$ chambers after the turning
edge, for $n \geq 1$.  The results for these $n$ are put together in 
Figure~\ref{OddAfterTE}.

\begin{figure}
\setlength{\unitlength}{1in}
\begin{picture}(6,7.5)(0,0)
\centerline{\includegraphics[height=7.5in]{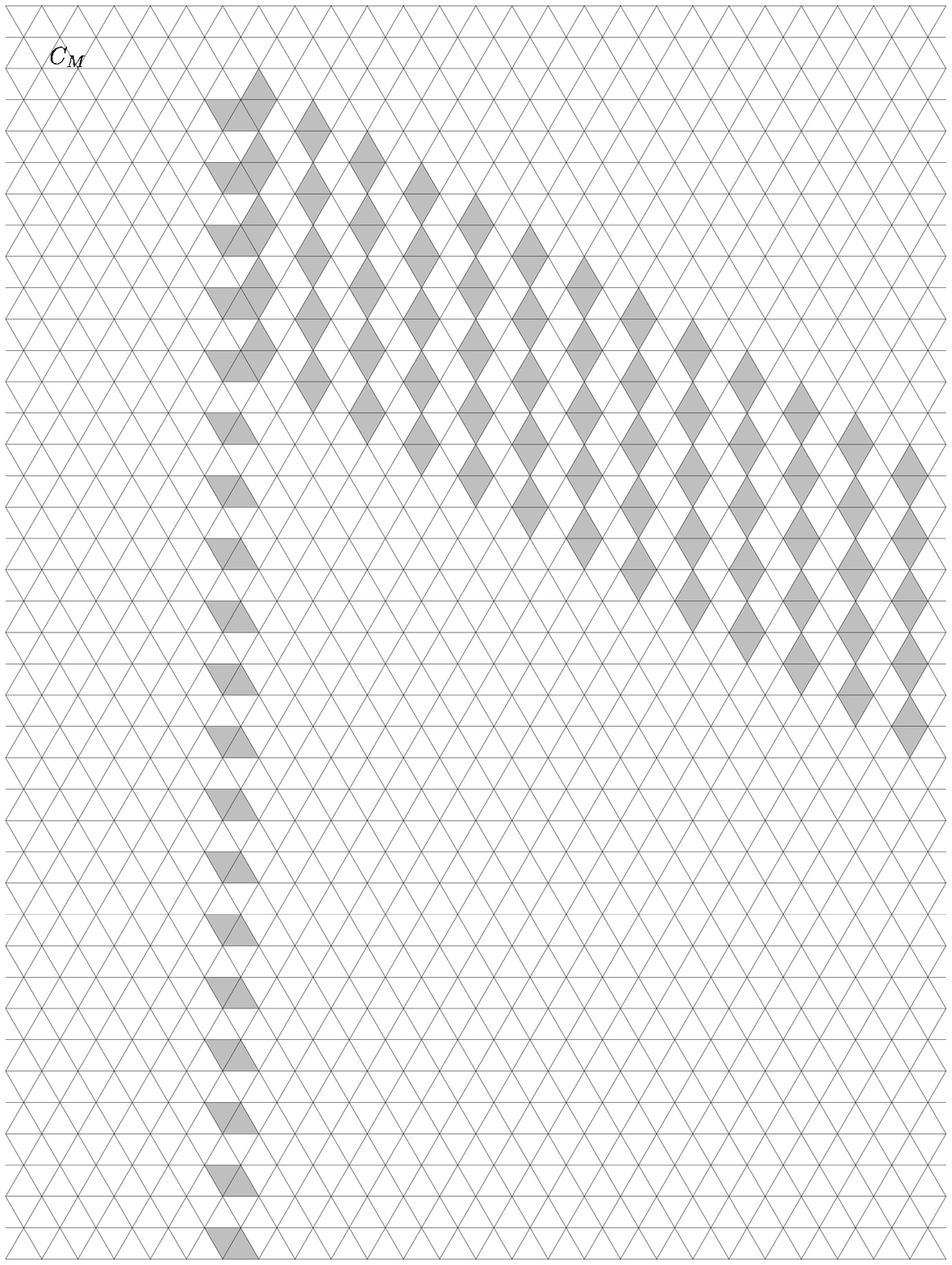}}
\end{picture}
\caption{Class $I_1$, where $\Gamma_{E_1}$ has $2$ chambers after the 
turning edge}
\label{TwoAfterTE}
\end{figure}

\begin{figure}
\setlength{\unitlength}{1in}
\begin{picture}(6,7.5)(0,0)
\centerline{\includegraphics[height=7.5in]{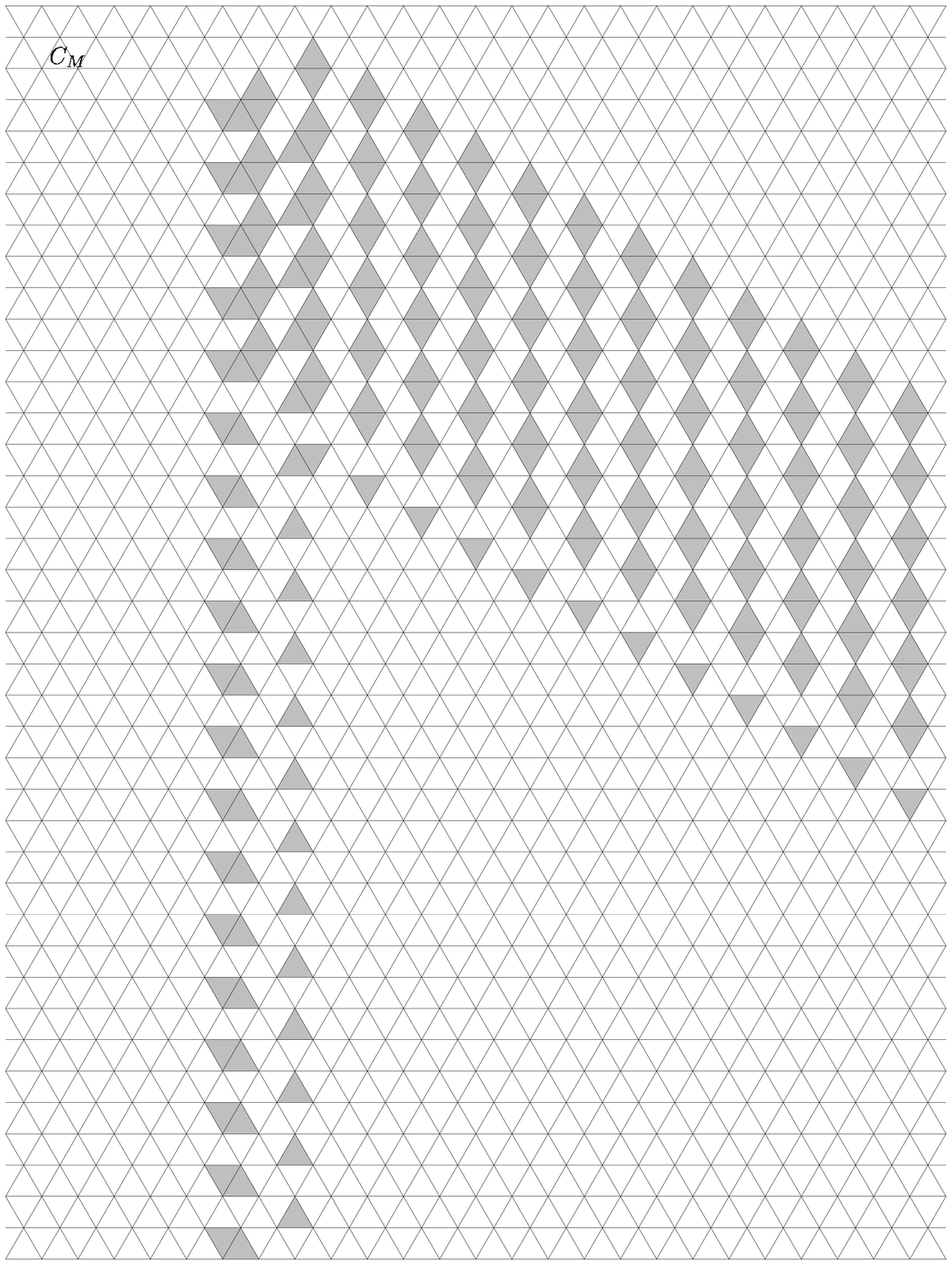}}
\end{picture}
\caption{Class $I_1$, where $\Gamma_{E_1}$ has 
$4$ chambers after the turning edge}
\label{FourAfterTE}
\end{figure}

\begin{figure}
\setlength{\unitlength}{1in}
\begin{picture}(6,7.5)(0,0)
\centerline{\includegraphics[height=7.5in]{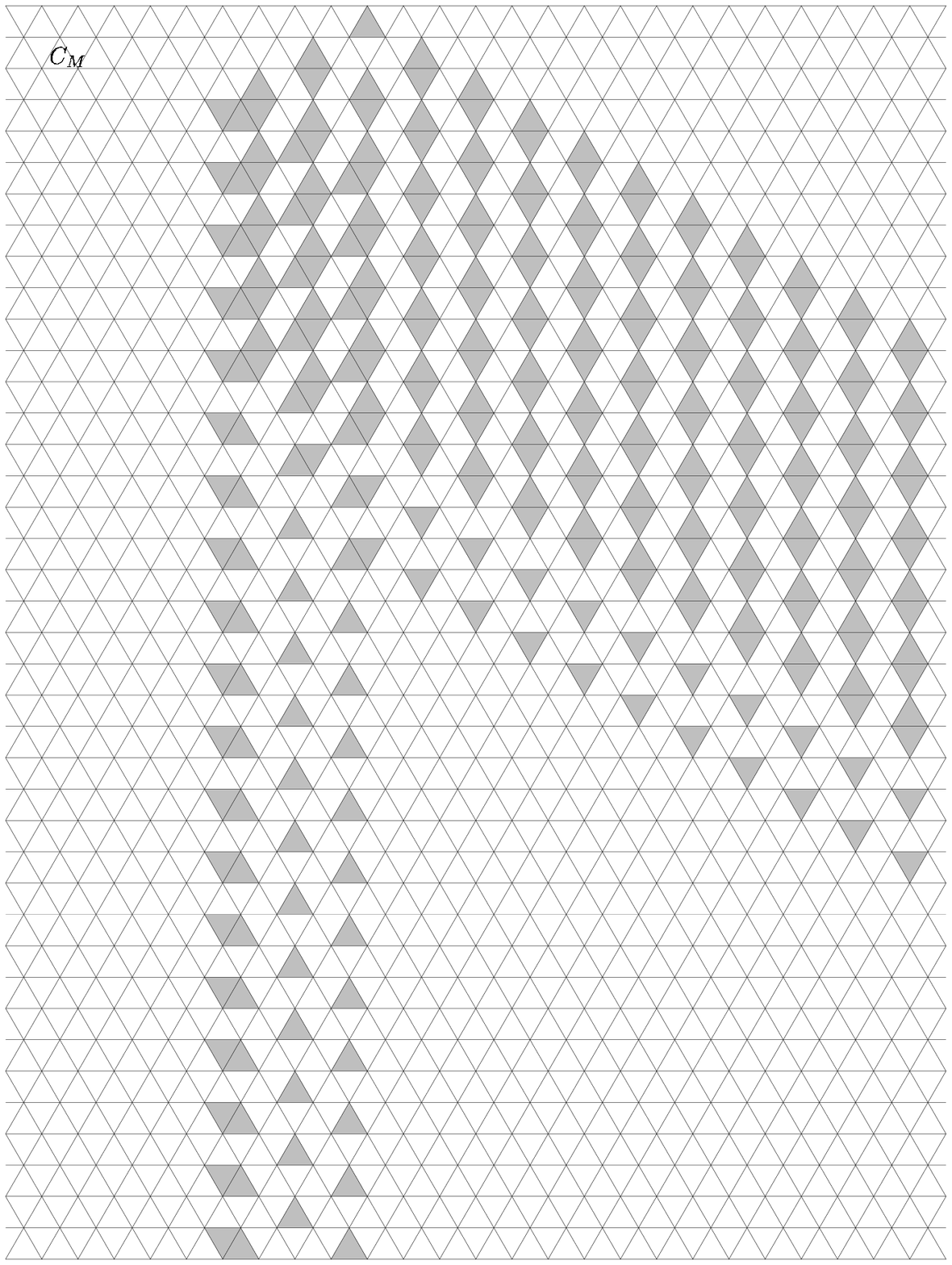}}
\end{picture}
\caption{Class $I_1$, where $\Gamma_{E_1}$ has $6$ chambers after the 
turning edge}
\label{SixAfterTE}
\end{figure}

\begin{figure}
\setlength{\unitlength}{1in}
\begin{picture}(6,7.5)(0,0)
\centerline{\includegraphics[height=7.5in]{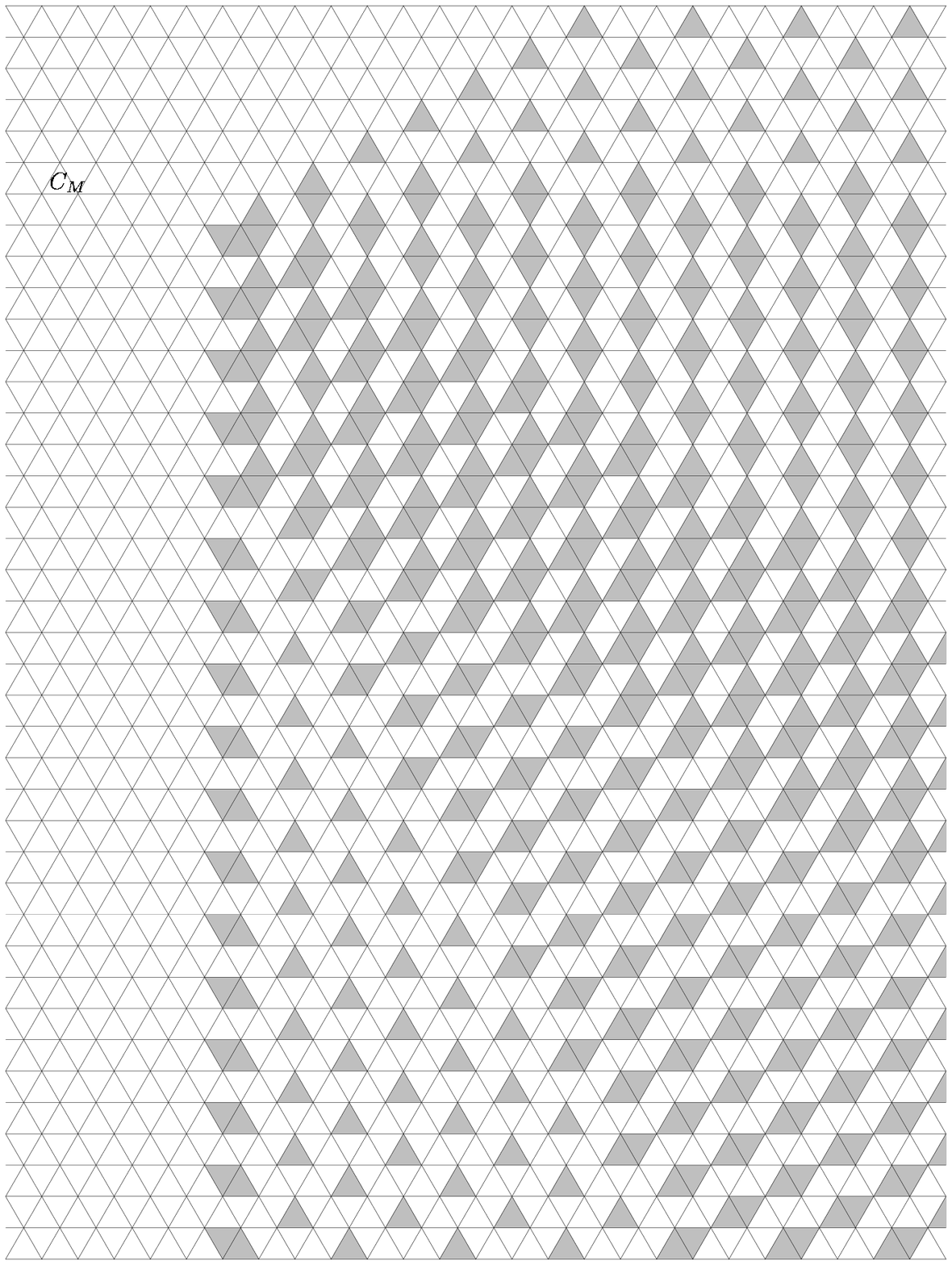}}
\end{picture}
\caption{Summary of results for $I_1$, an even number of chambers 
after the turning edge}
\label{EvenAfterTE}
\end{figure}

\noindent {\em Case 2:} $w = 1$.  We get a composite gallery of the type
shown in Figure~\ref{CompositeGalleryI1_2}.
\begin{figure}
\centerline{\input{fig21.pstex_t}}
\caption{More composite gallery possibilities for $I_1$}
\label{CompositeGalleryI1_2}
\end{figure}
Figures~\ref{TwoAfterTE}
,~\ref{FourAfterTE}, 
and~\ref{SixAfterTE} have instances where
$\Gamma_{E}^1$ has $2$, $4$, and $6$ chambers after the turning
edge, respectively.  Figure~\ref{EvenAfterTE} has the 
amalgamated results for $2n$ chambers after
the turning edge for all $n \geq 1$. 

It is easy to see that the results for the cases in which $w$ is some
other order $2$ element of $W$ are just rotations of the results for case
$1$ by $120^{\circ}$ and $240^{\circ}$ about the center point of the
chamber $C_M$.  The results for the cases in
which $w$ is some other order $3$ element of $W$ are rotations of the case
$2$ results by $120^{\circ}$ and $240^{\circ}$ about the center point of
$C_M$.  

We break $I_2$ into the same six sub-classes.  

\begin{figure}
\setlength{\unitlength}{1in}
\begin{picture}(3.5,2)(0,0)
\centerline{\includegraphics[height=2in]{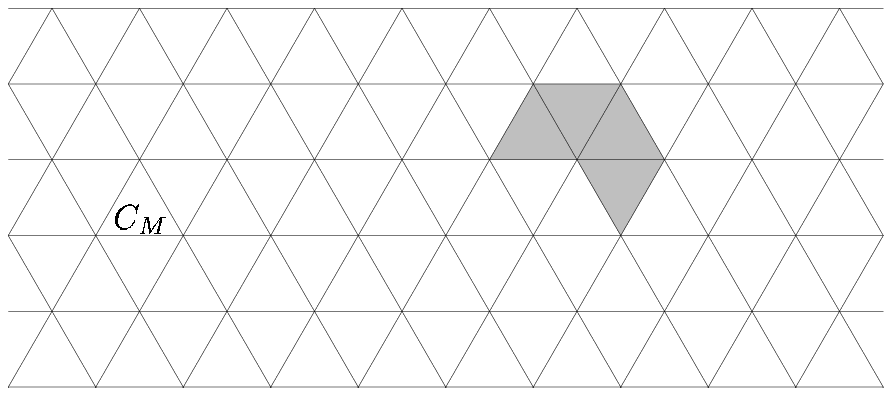}}
\end{picture}
\caption{Class $I_2$, where $\Gamma_{E_1}$ has 
$0,1$ or $2$ chambers after the departure edge, $w = 1$}
\label{TwoAfterDE}
\end{figure}

\begin{figure}
\setlength{\unitlength}{1in}
\begin{picture}(3.5,2)(0,0)
\centerline{\includegraphics[height=2in]{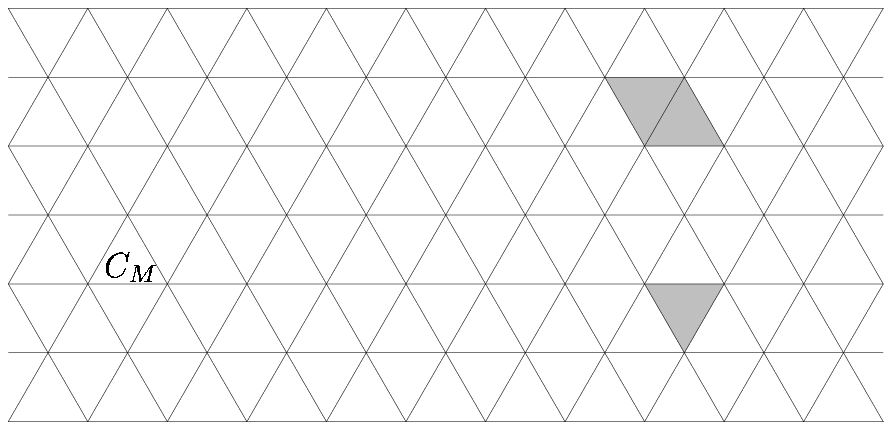}}
\end{picture}
\caption{Class $I_2$, where $\Gamma_{E_1}$ has 
$3$ or $4$ chambers after the departure edge, $w = 1$}
\label{FourAfterDE}
\end{figure}

\begin{figure}
\setlength{\unitlength}{1in}
\begin{picture}(4,2.75)(0,0)
\centerline{\includegraphics[height=2.75in]{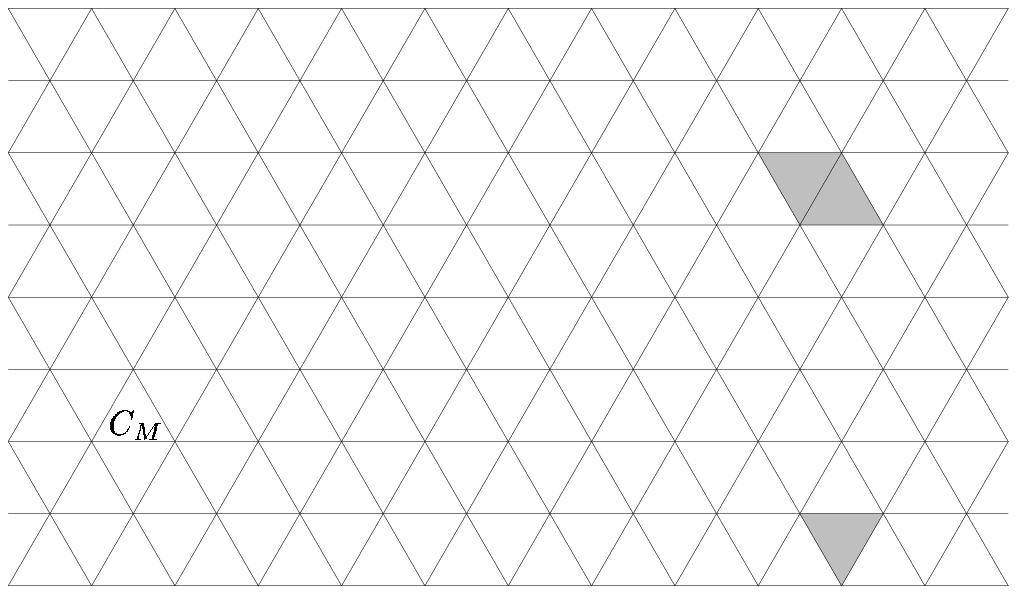}}
\end{picture}
\caption{Class $I_2$, where $\Gamma_{E_1}$ has 
$5$ or $6$ chambers after the departure edge, $w = 1$}
\label{SixAfterDE}
\end{figure}

\begin{figure}
\setlength{\unitlength}{1in}
\begin{picture}(5.5,6)(0,0)
\centerline{\includegraphics[height=6in]{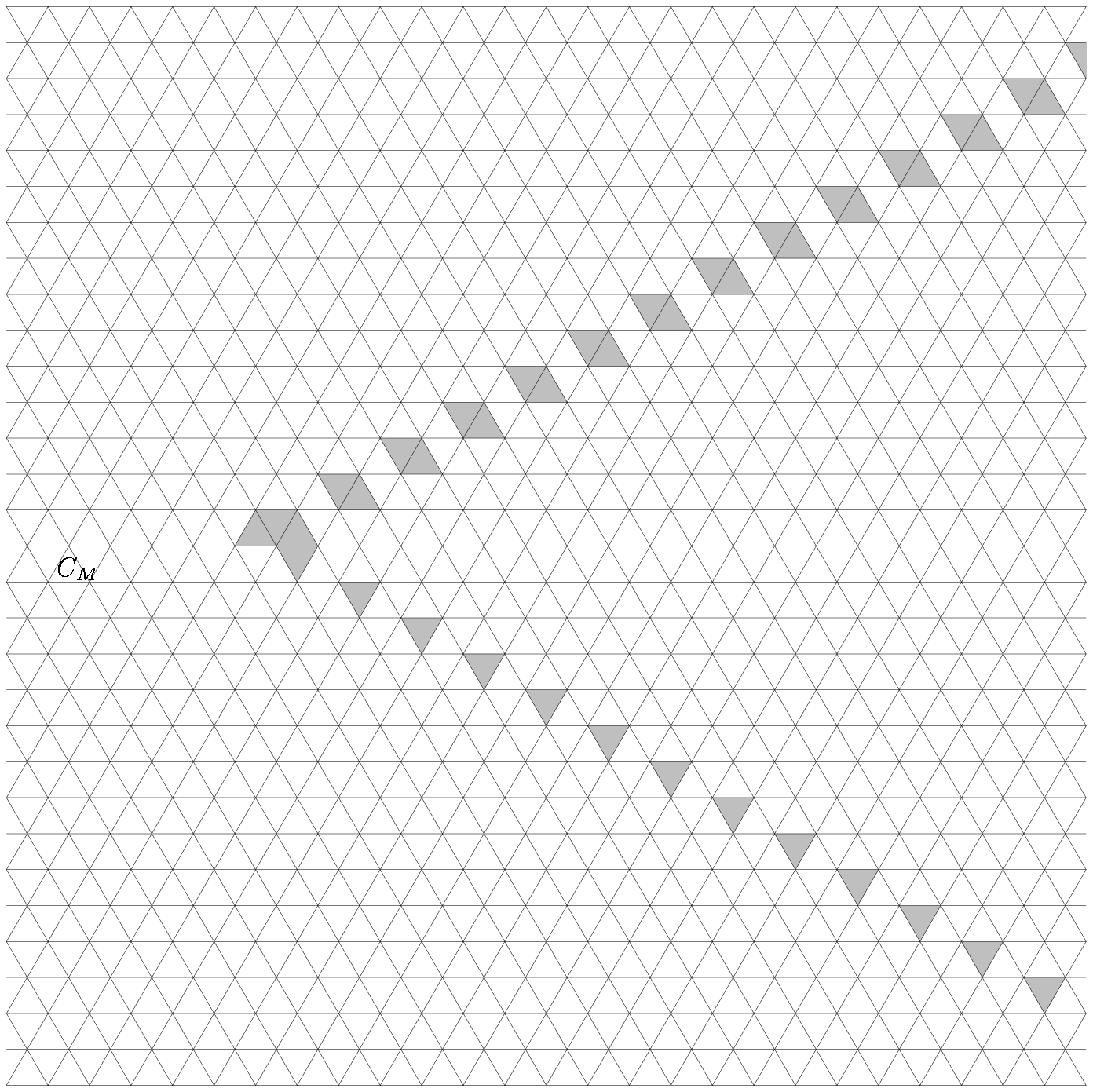}}
\end{picture}
\caption{Summary of results for $I_2$, $w = 1$}
\label{EvenAfterDE}
\end{figure}

\noindent {\em Case 1:} $w = 1$.  We get a composite gallery of one
of the two types
shown in Figures~\ref{CompositeGalleryI2_1} 
and~\ref{CompositeGalleryI2_2}, where $D$ is the last chamber in the
SMG from $C_M$ to $xC_M$ that is contained in $A_M$.
\begin{figure}
\centerline{\input{fig26.pstex_t}}
\caption{Some composite gallery possibilities for $I_2$}
\label{CompositeGalleryI2_1}
\end{figure}
\begin{figure}
\centerline{\input{fig101.pstex_t}}
\caption{More composite gallery possibilities for $I_2$}
\label{CompositeGalleryI2_2}
\end{figure}
We consider in Figure~\ref{TwoAfterDE} the instances where 
$\Gamma_E^1$ has $0, 1$ or $2$ chambers
after the departure edge, in Figure~\ref{FourAfterDE} 
the instances where it has $3$ or $4$,
and Figure~\ref{SixAfterDE} the instances where it has $5$ or $6$.  
Figure~\ref{EvenAfterDE} is the combined results 
for $n$ chambers after the departure edge, for $n \geq 0$.

\begin{figure}
\setlength{\unitlength}{1in}
\begin{picture}(5.5,5.25)(0,0)
\centerline{\includegraphics[height=5.25in]{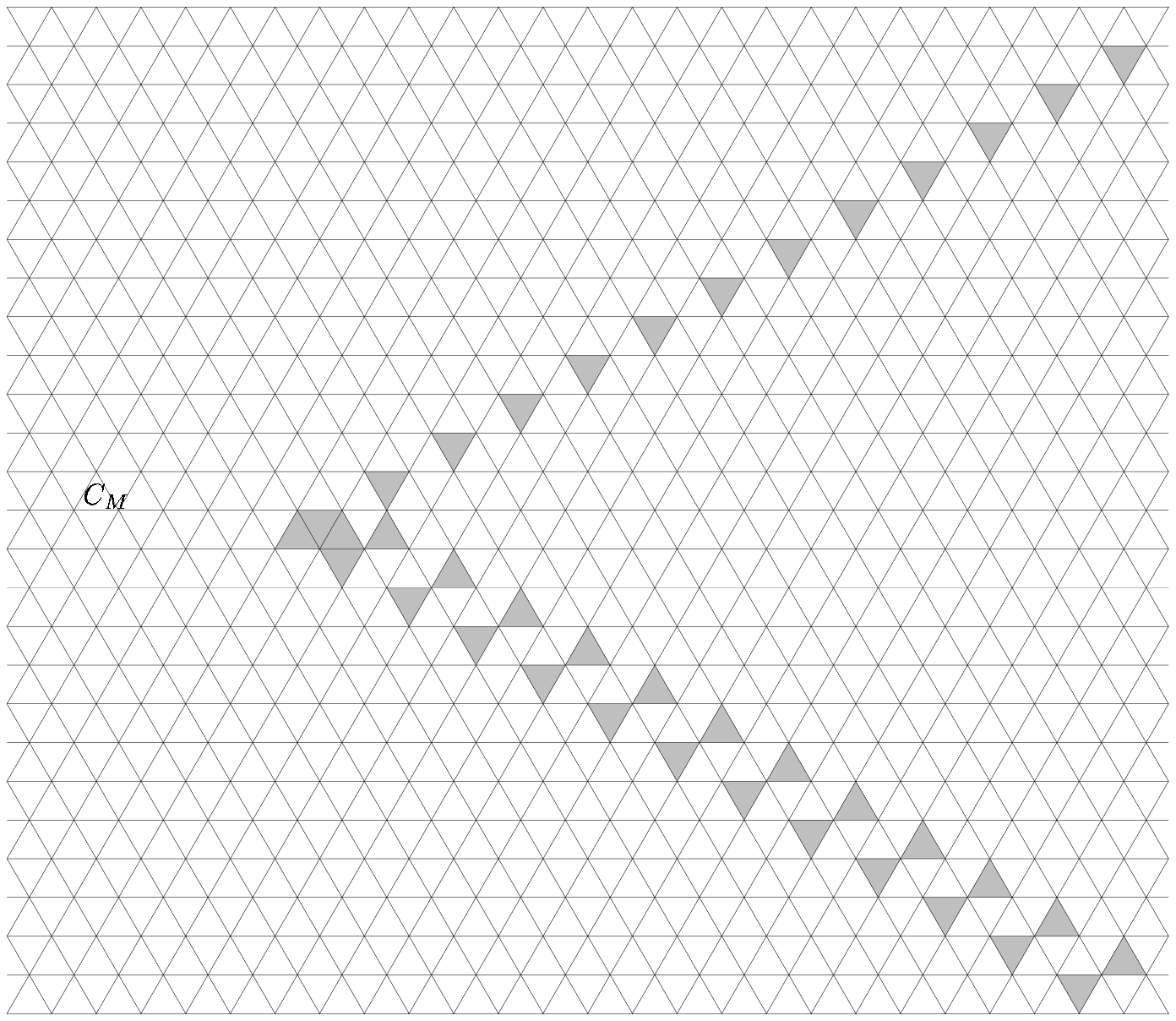}}
\end{picture}
\caption{Summary of results for $I_2$, $w = f$}
\label{OddAfterDE}
\end{figure}

\noindent {\em Case 2:} $w = f$.
The combined results for this case are in Figure~\ref{OddAfterDE}.

\begin{figure}
\setlength{\unitlength}{1in}
\begin{picture}(6,7.5)(0,0)
\centerline{\includegraphics[height=7.5in]{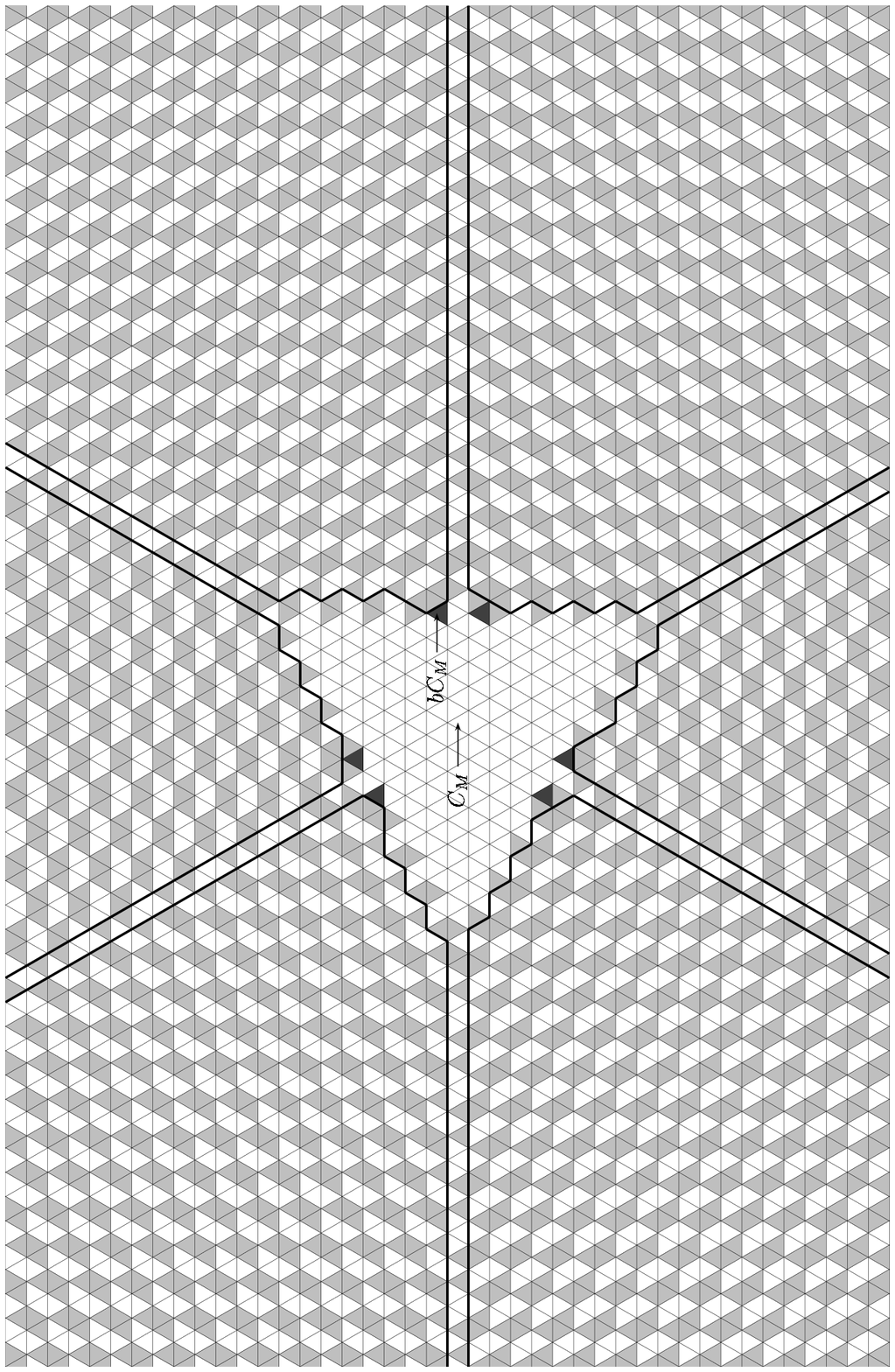}}
\end{picture}
\caption{Summary of results for $I_1$ and $I_2$}
\label{I1I2Whole}
\end{figure}

Again, if $w$ is some other element of $W$, then the results are rotations
of one of the above cases by $120^{\circ}$ and $240^{\circ}$ about the center
of $C_M$.
Conglomerating all results from $I_1$ and $I_2$ gives the results pictured
in Figure~\ref{I1I2Whole}. The chambers
which are shaded more darkly in this figure are the chambers
$w^{-1}bwC_M$ for $w \in W$.

We will not compute other infinite classes of type-edge pairs here, but we
will describe how the collection of type-edge pairs can be divided into
infinite classes.  For each corridor $c_i$, we get two infinite classes,
$I_{c_i}^1$ and $I_{c_i}^2$.  The SMG types in each of these 
classes are those corresponding to SMGs whose terminus is in $c_i$.
If $\eta_1$ and $\eta_2$ are the two possible angles of edges
of departure of these SMGs, then $I_{c_i}^j$ has type-edge pairs
whose edge of departure 
component has angle $\eta_j$.  Note that $I_2$ is the union of two infinite
classes of this kind, namely $I_{c_4}^1$ and $I_{c_4}^2$.

If $G$ is an SMG type corresponding to an SMG whose final chamber is in
one of the regions $R_i$ or $r_i$, then $G$ has a turning point.  The edge
of departure for $G$ could be before, after, or at the turning point.
However, we need not consider type-edge pairs for which the edge of departure
is after or at the turning point, because $\Gamma_E^1$ for these type-edge
pairs is the same as that arising from some type-edge pair in 
some $I_{c_i}^j$.

So for each $i$ we have $I_{R_i}^1$ and $I_{R_i}^2$ with SMG types
corresponding to an SMG whose final chamber is in $R_i$.  The allowed edges
of departure are those before the turning edge, with angle specified 
by the superscript.  Note that $I_1$ is of this type. 
We define $I_{r_i}^1$ and $I_{r_i}^2$ analogously for $i = 1,2,3$.

When one conglomerates the results from all infinite classes of composite
galleries for
$$b = \left( 
\begin{matrix}
	\pi^3 & 0 & 0 \\
	0 & \pi^{-1} & 0 \\
	0 & 0 & \pi^{-2}
\end{matrix}
\right),$$
one gets an upper bound set $S_1$.  It turns out that infinite classes
other than those subsumed by $I_1$ and $I_2$ do not
contribute any chambers beyond those contributed by $I_1$ and $I_2$.  So
Figure~\ref{I1I2Whole} is the complete superset for $b$ with $\alpha = 3$ and
$\beta = -1$.

\begin{figure}
\setlength{\unitlength}{1in}
\begin{picture}(6,7.5)(0,0)
\centerline{\includegraphics[height=7.5in]{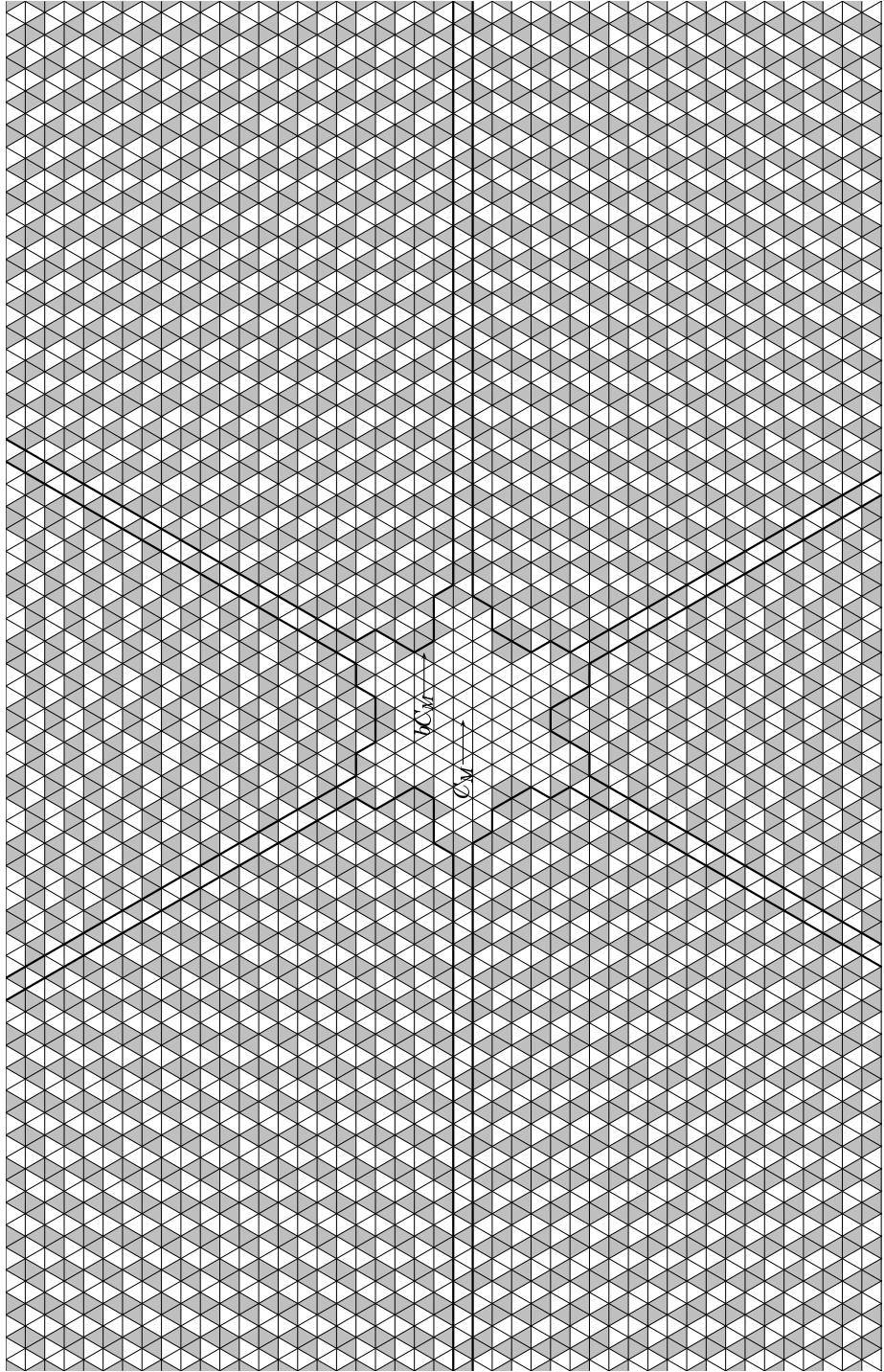}}
\end{picture}
\caption{Result for $\alpha = 2$, $\beta = 0$}
\label{SomeResults2_0}
\end{figure}

When one conglomerates the results of the classes $I_1$ and $I_2$ for 
$$b= \left(
\begin{matrix}
	\pi^2 & 0 & 0 \\
	0 & \pi^0 & 0 \\
	0 & 0 & \pi^{-2}
\end{matrix}
\right),$$ 
one gets the chambers pictured in Figure~\ref{SomeResults2_0}.  The 
results of other infinite classes have not been computed in this 
case, but is reasonable to expect that they would not contribute any 
chambers additional to those contributed by $I_1$ and $I_2$.  This will
be discussed further later.

\begin{figure}
\setlength{\unitlength}{1in}
\begin{picture}(6,7.5)(0,0)
\centerline{\includegraphics[height=7.5in]{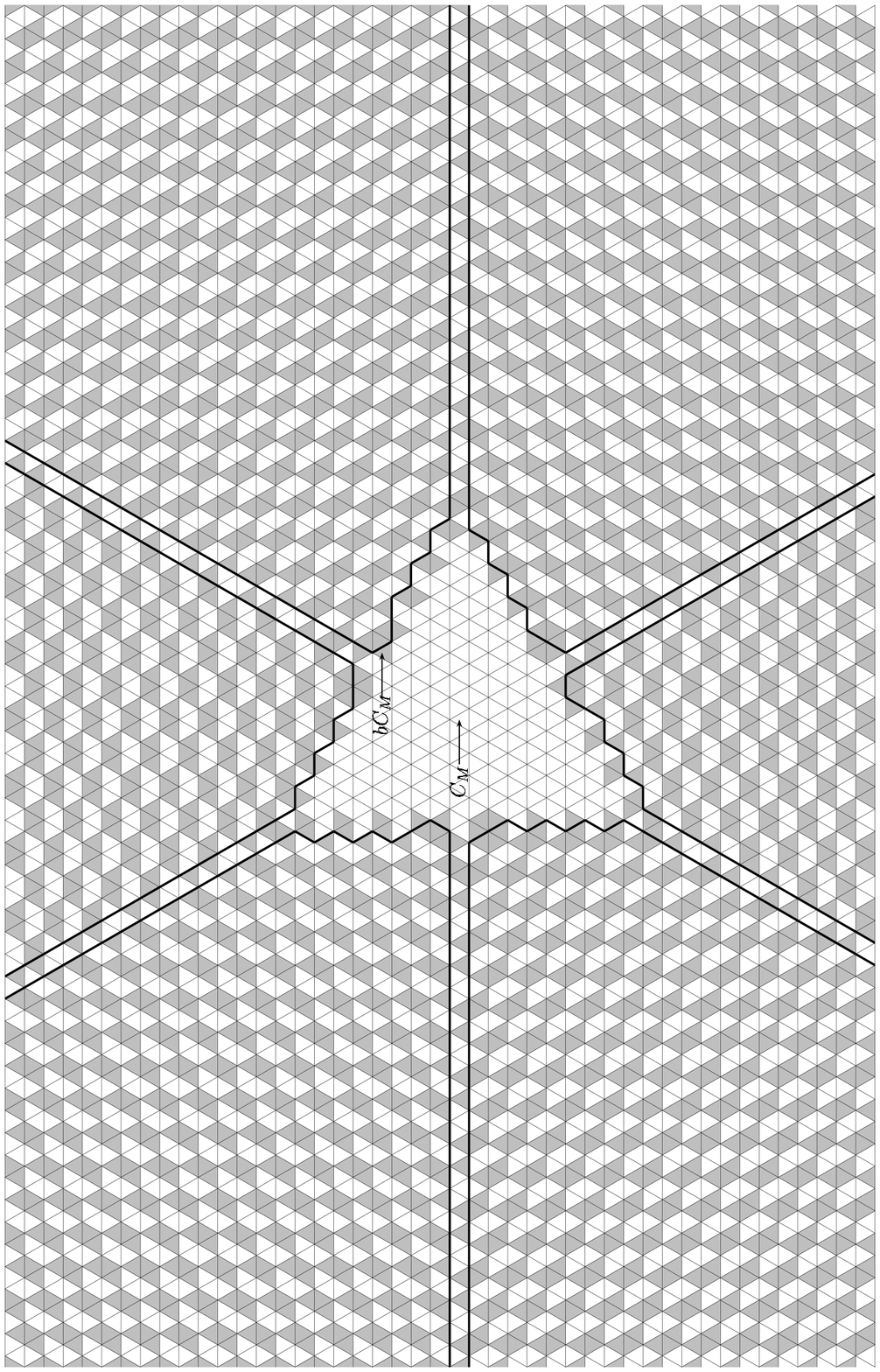}}
\end{picture}
\caption{Result for $\alpha = 2$, $\beta = 1$}
\label{SomeResults2_1}
\end{figure}

When one conglomerates the results of the classes $I_1$ and $I_2$ for 
$$b= \left(
\begin{matrix}
	\pi^2 & 0 & 0 \\
	0 & \pi^1 & 0 \\
	0 & 0 & \pi^{-3}
\end{matrix}
\right),$$ 
one gets the chambers pictured in Figure~\ref{SomeResults2_1}.  Again,
other infinite classes have not been computed, but it is reasonable not to,
as will be discussed later.

\begin{figure}
\setlength{\unitlength}{1in}
\begin{picture}(6,7.5)(0,0)
\centerline{\includegraphics[height=7.5in]{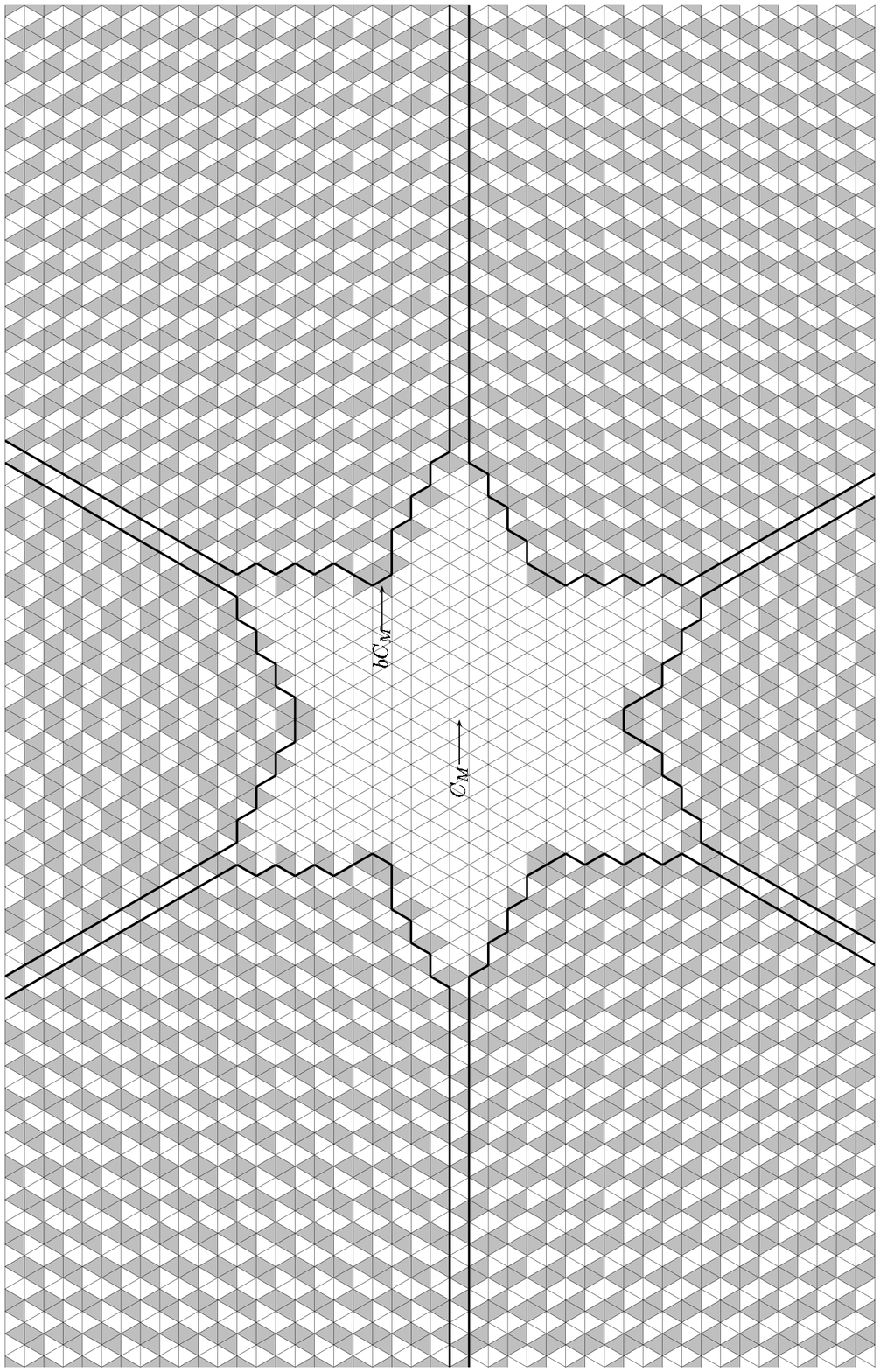}}
\end{picture}
\caption{Result for $\alpha = 4$, $\beta = 0$}
\label{SomeResults4_0}
\end{figure}

When one conglomerates the results of the classes $I_1$ and $I_2$ for 
$$b= \left(
\begin{matrix}
	\pi^4 & 0 & 0 \\
	0 & \pi^0 & 0 \\
	0 & 0 & \pi^{-4}
\end{matrix}
\right),$$ 
one gets the chambers pictured in Figure~\ref{SomeResults4_0}.  Again,
other infinite classes have not been computed, but it is reasonable not to.

\begin{figure}
\setlength{\unitlength}{1in}
\begin{picture}(6,7.5)(0,0)
\centerline{\includegraphics[height=7.5in]{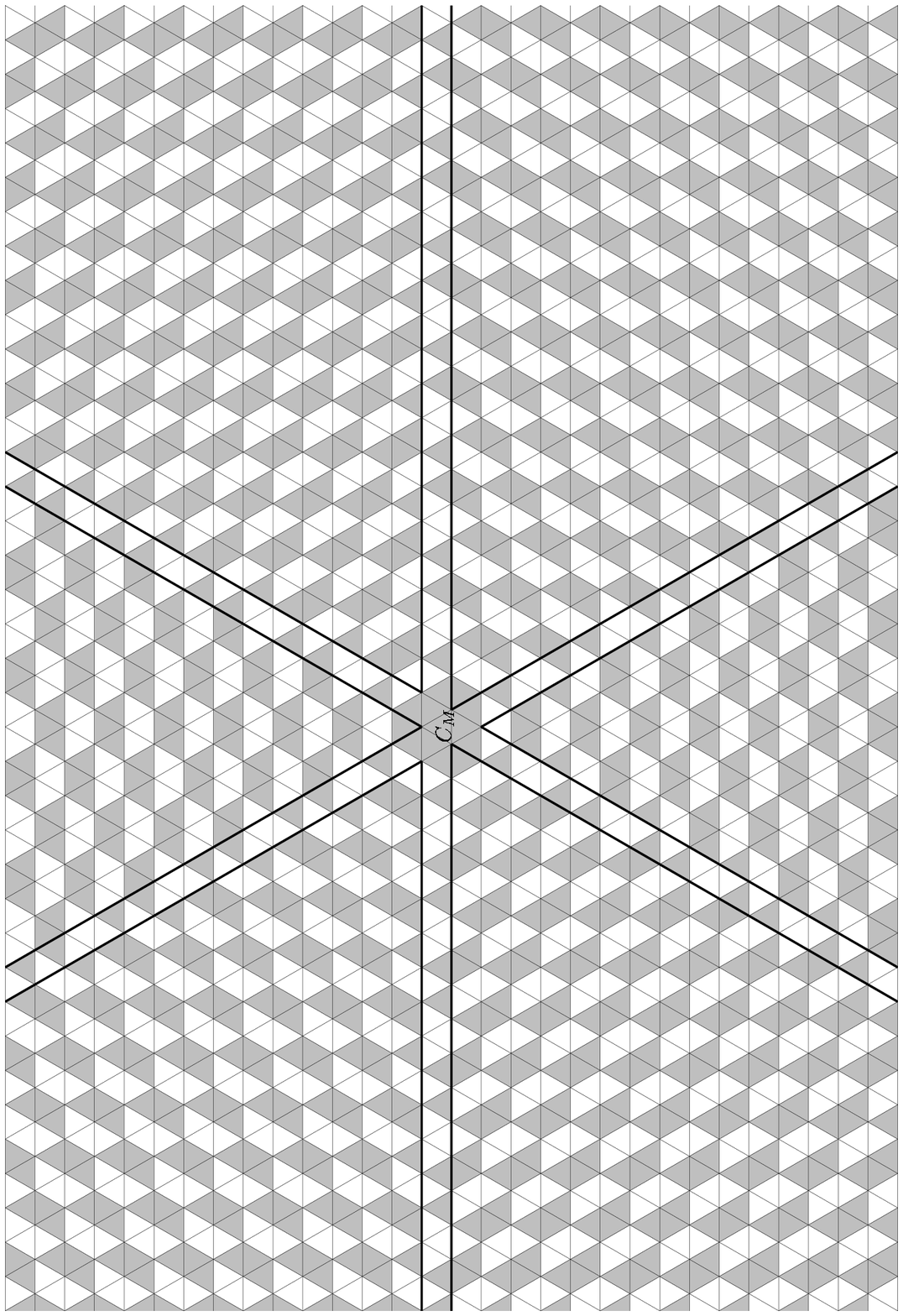}}
\end{picture}
\caption{Result for $b=1$}
\label{TotalResults0_0}
\end{figure}

When one conglomerates the results of the classes $I_1$ and $I_2$ for $b=1$, 
one gets the chambers pictured in 
Figure~\ref{TotalResults0_0}.  We have also computed
the results of all other infinite classes of composite galleries in the 
$b = 1$ case, as we did for the $\alpha = 3$, $\beta = -1$ case.  Again, 
one can see (after the fact) that the infinite classes beyond 
those subsumed by $I_1$ and $I_2$ do not contribute any additional chambers.

The computation of the results of $I_1$ and $I_2$ for $b = 1$ is a little
different from the cases done thus far, all of which are similar to
the example $\alpha = 3$, $\beta = -1$ that was done in detail.  
The general shapes for 
composite galleries for $b = 1$ with $(t,e) \in I_1$ 
are shown  in Figure~\ref{CompositeGalleryShapeb1wf} for $w = f$ and 
Figure~\ref{CompositeGalleryShapeb1w1} for $w = 1$.
\begin{figure}
\centerline{\input{fig102.pstex_t}}
\caption{General shape of composite galleries for $b = 1$, $w = f$, $I_1$}
\label{CompositeGalleryShapeb1wf}
\end{figure}
\begin{figure}
\centerline{\input{fig103.pstex_t}}
\caption{General shape of composite galleries for $b = 1$, $w = 1$, $I_1$}
\label{CompositeGalleryShapeb1w1}
\end{figure}
The general shapes for $b = 1$ with $(t,e) \in I_2$
are shown in Figure~\ref{CompositeGalleryShapeb1w1I2} for $w = 1$.  The 
$w = f$ for $I_2$ galleries would be similar.  All these galleries 
are different in general shape from their $\alpha = 3$, $\beta = -1$
counterparts.  Note that the $I_2$ galleries are minimal, 
and therefore would be easy to fold.
\begin{figure}
\centerline{\input{fig104.pstex_t}}
\caption{General shapes of composite galleries for $b = 1$, $w = 1$, $I_2$}
\label{CompositeGalleryShapeb1w1I2}
\end{figure}

\begin{figure}
\setlength{\unitlength}{1in}
\begin{picture}(6,7.5)(0,0)
\centerline{\includegraphics[height=7.5in]{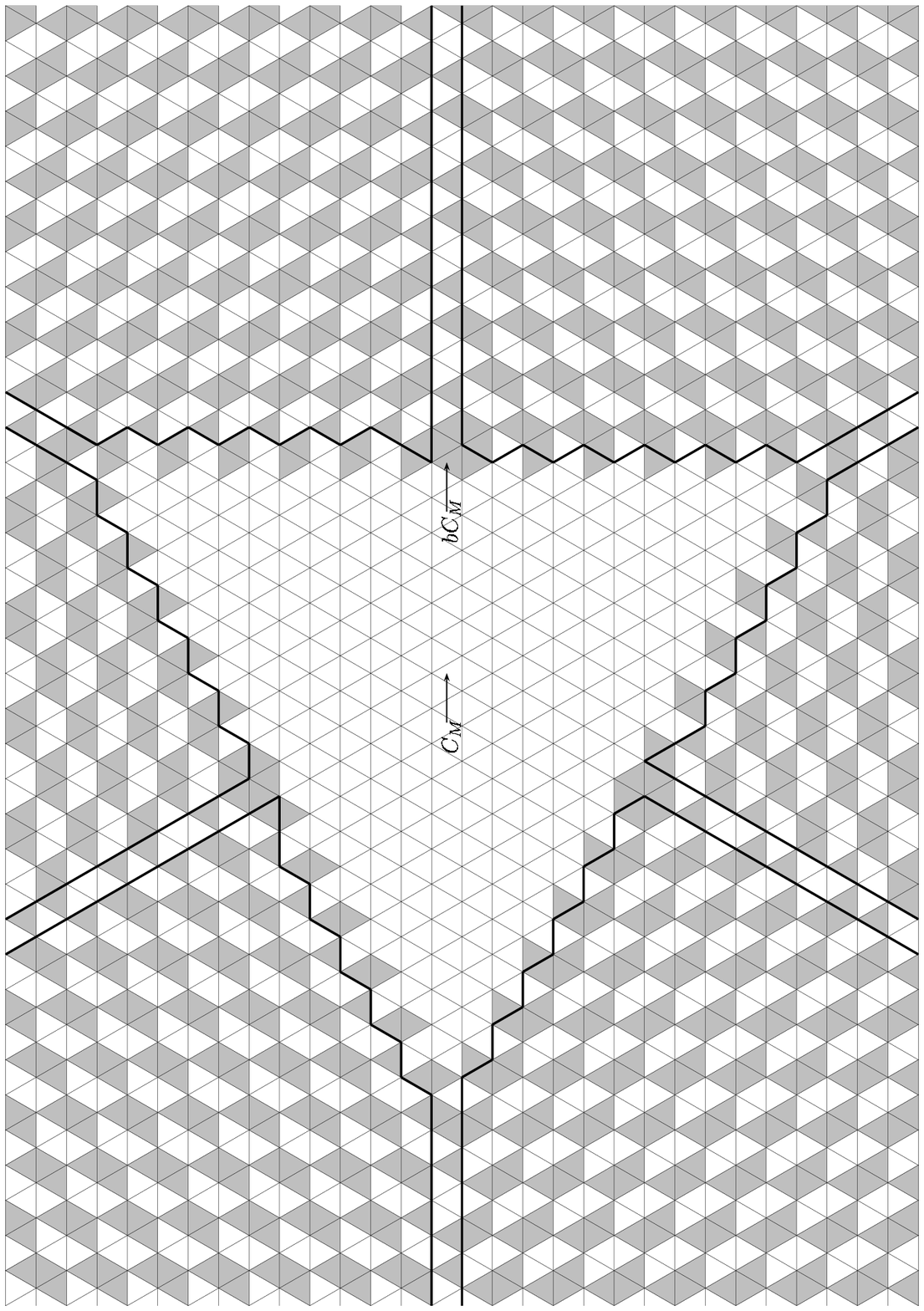}}
\end{picture}
\caption{Result for $\alpha = 4$, $\beta = -2$}
\label{SomeResults4_minus2}
\end{figure}

When one conglomerates the results of the classes $I_1$ and $I_2$ for 
$$b= \left(
\begin{matrix}
	\pi^4 & 0 & 0 \\
	0 & \pi^{-2} & 0 \\
	0 & 0 & \pi^{-2}
\end{matrix}
\right),$$
one gets the chambers pictured in Figure~\ref{SomeResults4_minus2}.

\begin{figure}
\centerline{\input{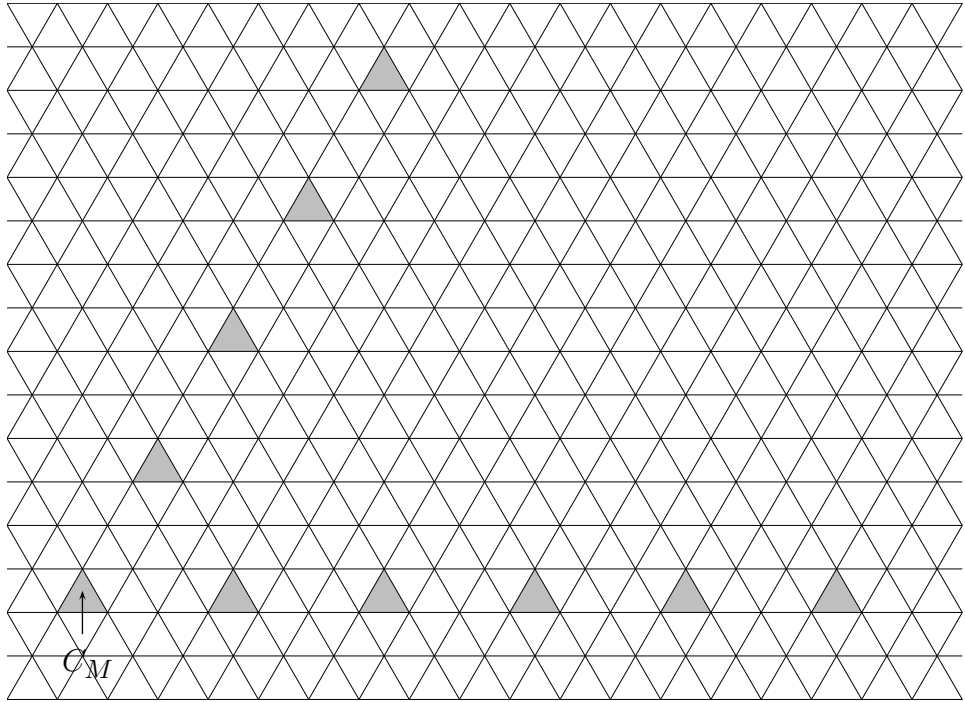}}
\caption{$bC_M$ for degenerate $b$}
\label{degen}
\end{figure}

The computation of the results of $I_1$ and $I_2$ for this last value of $b$
is also a little different from the cases done thus far because 
$b \neq 1$, but $\alpha + 2\beta = 0$.  This is
to say that $bC_M$ sits along the bottom edge of the positive Weyl
chamber.  We will call any $b$ with $\alpha + 2\beta = 0$ or
$\alpha + 2\beta = 2\alpha + \beta$ {\em degenerate}.  The values of
$bC_M$ for degenerate $b$ are pictured in Figure~\ref{degen}.  Note that
the condition $\alpha + 2\beta = 2\alpha + \beta$ corresponds to the $bC_M$
along the top-left boundary of the positive Weyl chamber.

\begin{figure}
\centerline{\input{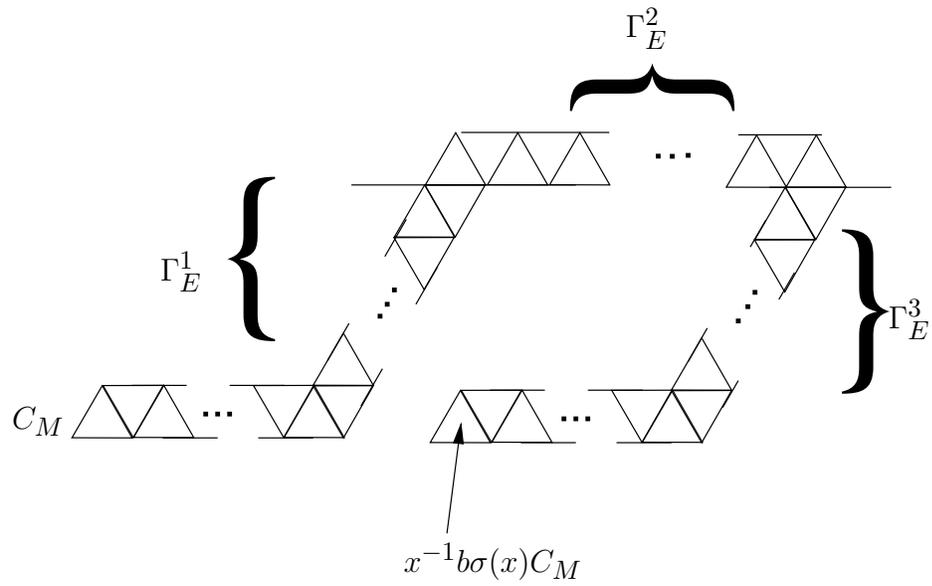}}
\caption{General shape of composite galleries for $b \neq 1$ with 
$\alpha + 2\beta = 0$, $w = f$, $I_1$}
\label{bottomdegenI1wf}
\end{figure}

\begin{figure}
\centerline{\input{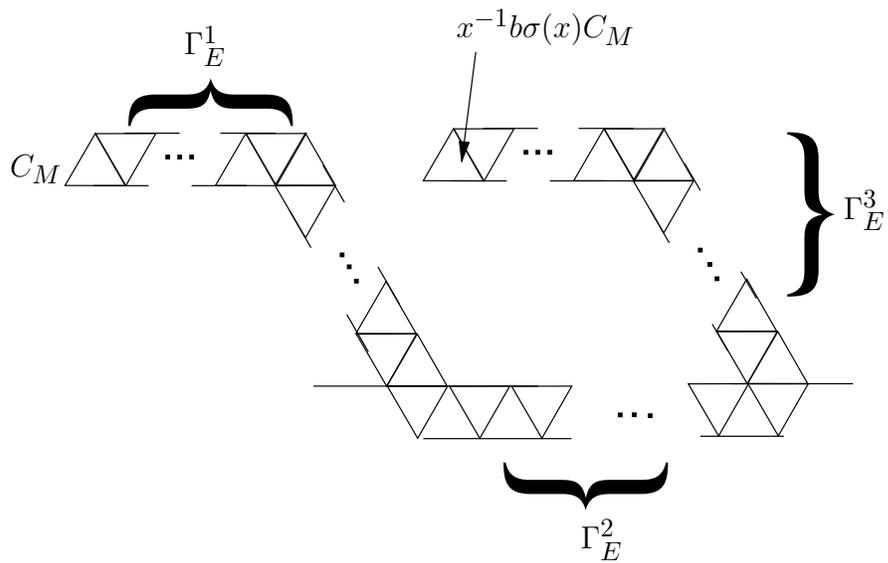}}
\caption{General shape of composite galleries for $b \neq 1$ with 
$\alpha + 2\beta = 0$, $w = 1$, $I_1$}
\label{bottomdegenI1w1}
\end{figure}

The general shapes of composite galleries for $b \neq 1$, 
$\alpha + 2\beta = 0$, $(t,e) \in I_1$ are pictured in 
Figure~\ref{bottomdegenI1wf} for
$w = f$ and Figure~\ref{bottomdegenI1w1} 
for $w = 1$.  The composite galleries for
$I_2$ in these cases are similar enough in general shape to those of 
non-degenerate $b$ that folding considerations are essentially the same.

In the case that $\alpha + 2\beta = 2\alpha + \beta$, 
$b \neq 1$, the composite galleries
for $I_2$ become different in general shape from those of non-degenerate
$b$, but those for $I_1$ are similar.

Even though the general
shape of composite galleries is different in each of the four categories:
$b$ non-degenerate, $b = 1$, $\alpha + 2\beta = 0$ but $b \neq 1$,
$\alpha + 2\beta = 2\alpha + \beta$ but $b \neq 1$, there is an efficient way
of computing the folding results of these galleries that 
turns out to be very similar
in all four cases.  This efficient method is somewhat different from that
used previously to compute the $\alpha = 3$, $\beta = -1$ example (although
the efficient method will also apply to that example), and 
will be described in 
Section~\ref{Relationship3}.

The lines superimposed on the patterns in 
Figures~\ref{I1I2Whole},~\ref{SomeResults2_0},~\ref{SomeResults2_1},~\ref{SomeResults4_0},~\ref{TotalResults0_0}, and~\ref{SomeResults4_minus2} 
are only there to
make the patterns easier to view and to compare.  The chambers
which are shaded more darkly in Figure~\ref{I1I2Whole} are the chambers
$w^{-1}bwC_M$ for $w \in W$.  These are not shaded on the other figures
because one can tell where they are by analogy (or by an easy computation).  
Note that for 
Figure~\ref{TotalResults0_0} and Figure~\ref{SomeResults4_minus2}, 
some of the $w^{-1}bwC_M$ are equal
to each other.  These are degenerate cases, and will be discussed further
in Section~\ref{Relationship3}.

It is reasonable to expect that other infinite
classes would not contribute further to the results 
in Figures~\ref{SomeResults2_0},~\ref{SomeResults2_1}, 
or~\ref{SomeResults4_0}, and~\ref{SomeResults4_minus2} 
as they do not contribute further to the results
in Figures~\ref{I1I2Whole} and~\ref{TotalResults0_0}.

It is easy to generalize what the results of classes $I_1$ and
$I_2$ would be for any of the $b$ we are working with (listed in the beginning
of Section~\ref{SL3}).  This will also be done in detail in 
Section~\ref{Relationship3}.  Since the computations
involved in discovering
the total result of all infinite classes of 
type-edge pairs are very lengthy for any given $b$, 
no computation general to all $b$ has been done.
It seems very reasonable that in any case the total results would be
the same as the $I_1$ and $I_2$ results together, as mentioned previously.

\subsection{Relationship Between $X_w (1\sigma)$ and 
$X_w (b\sigma)$, and an Efficient Way of Computing 
Supersets}\label{Relationship3}

In this section we try to illuminate the cause of some of the 
similarities among the supersets $S_1$ for
different values of $b$.  We do this by discussing a more efficient method
of computing the infinite classes $I_1$ and $I_2$ that were computed 
in Section~\ref{Superset}.  The basic idea is as follows.  Let $\Gamma_{E_1}$
and $\Gamma_{E_2}$ be two composite galleries (constructed in 
Section~\ref{SMGandComp}) that need to be folded in order to compute one 
of the supersets of Section~\ref{Superset}.  Let $p_i$ for $i = 1,2$ 
be the number of chambers appearing in $\Gamma_{E_i}^1$ between the edge
of departure and the turning edge, let $q_i$ for $i = 1,2$ be the
number of chambers appearing in $\Gamma_{E_i}^1$ after the turning edge,
and let $w_i$ for $i = 1,2$ be elements of $W$ such that $a_i w_i C_M 
= \rho(E_i)$ for some $a_i \in T(F)$.
If $p_1 = p_2$ and $w_1 = w_2$, then the
folding results of $\Gamma_{E_1}$ and $\Gamma_{E_2}$ can be computed
in very similar ways.  In fact, for most fixed values of $p_1 = p_2$,
for any fixed values of $w_1 = w_2$, and for most $b$, 
one can compute folding results for arbitrary $q_i$
all at once.  This is done by replacing all the composite
galleries with the given $p$ and $w$ with a single ``half-infinite gallery''
which one can show in advance will have the same folding results as all
these composite galleries together.  The folding results of the ``half-infinite
gallery'' are easy to compute.

At the end of Section~\ref{SMGandComp} we defined 
the {\em transition type} between 
two chambers that share an edge, and the {\em type} of a non-stuttering 
gallery $G$ in $\mathcal{B}_{\infty}$.  We also need the following definitions.
\begin{defn}
A {\em finite gallery type} is, equivalently, 

$1)$ An initial chamber $C_0 \subseteq A_M$ and a finite sequence
of transition types, which we will call $\{t_1,t_2,\ldots, n \}$

$2)$ A non-stuttering finite gallery in $A_M$.
\end{defn}
\begin{defn}
A {\em left-infinite gallery type} is, equivalently,

$1)$ A terminal chamber $C_0 \subseteq A_M$ and a left-infinite sequence
$\{ \ldots , t_{-3} , t_{-2} , t_{-1} , t_0 \}$ of transition types.

$2)$ A left-infinite non-stuttering gallery in $A_M$.
\end{defn}

Let $\{ C_M ; t_1 , t_2 , \ldots , t_n \}$ be a finite gallery type 
whose initial gallery is $C_M$.  Let $C_M = C_0 , C_1 , \ldots , C_n$
be the corresponding non-stuttering gallery in $A_M$.  We call this the 
{\em standard folding} of the gallery type 
$\{ C_M ; t_1 , t_2 , \ldots , t_n \}$.  There are {\em non-standard
foldings} of $\{ C_M ; t_1, t_2 , \ldots , t_n \}$ determined as follows.
We first assume inductively that $C_i$ is determined.  Note that the base case,
$C_0 = C_M$ actually is determined.  Then let $L_{i+1}$
be the wall containing the edge of $C_i$ that does not contain the 
vertex of type $t_{i+1}$.  If $C_M$ and $C_i$ are on the same side of
$L_{i+1}$, then let $C_{i+1}$ be the reflection of $C_i$ about $L_{i+1}$.
If $C_M$ and $C_i$ are on opposite sides of $L_{i+1}$, then let
$C_{i+1}$ be either $C_i$ itself, or the reflection of $C_i$ about 
$L_{i+1}$.  The choices that arise in this process may lead to many
different non-standard foldings of $\{ C_M ; t_1 , t_2 , \ldots , t_n \}$.
\begin{defn}
The {\em folding results} of $\{ C_M ; t_1 , t_2 , \ldots , t_n \}$ are
the collection of possible final chambers $C_n$ that can arise from all
possible (standard and non-standard) foldings of $\{ C_M ; t_1 , \ldots ,
t_n \}$.
\end{defn}
We now let $\{ \ldots , t_{-3} , t_{-2} , t_{-1} , t_0 ; C_0 \}$
be a left-infinite gallery type, and we let $\{ \ldots , C_{-3} , C_{-2} ,
C_{-1} , C_0 \}$ be the corresponding non-stuttering gallery in $A_M$.
Then it is easy to see that 
for any $n$, $\{ C_{-n} ; t_{-n+1} , t_{-n+2} , \ldots , t_0 \}$
is a finite gallery type.  Let $R_n$ be the folding results of this
finite gallery type.  Let $R = \cup_{n=1}^{\infty} R_n$.  We call $R$
the {\em folding results} of 
$\{ \ldots , t_{-3} , t_{-2} , t_{-1} , t_0 ; C_0 \}$.

Each type-edge pair $(t,e) \in I_1 \cup I_2$ gives $\Gamma_x$ a composite
gallery, and the folding results of that type-edge pair, $S_{(t,e)}$
are the folding results of the finite gallery type that $\Gamma_x$ 
represents.  In fact, we frequently have $S_{(t_1,e_1)} = S_{(t_2,e_2)}$
for $(t_1,e_1) \neq (t_2,e_2)$.  If $(t,e) \in I_1$, then there are 
three characteristics of $(t,e)$ that are clearly relevant to the type
of $x^{-1}\Gamma_x$, and therefore relevant to $S_{(t,e)}$.
First, the number $p$ of chambers between $e$ and the turning edge of
$t$ will have significance.  Second, the number $q$ of chambers after
the turning edge of $t$ will have significance.  And third, the 
Weyl group element $w \in W$ such that $awC_M = C_n$ will
have significance, where $a$ is a diagonal matrix, and $C_n$ is the last 
chamber of $t$. It is easy to see that if $p_{(t_1,e_1)} = p_{(t_2,e_2)}$,
$q_{(t_1,e_1)} = q_{(t_2,e_2)}$, and $w_{(t_1,e_1)} = w_{(t_2,e_2)}$
then $S_{(t_1,e_1)} = S_{(t_2,e_2)}$, for $(t_i,e_i) \in I_1$.
Let $R^1_{p,q,w}$ be this set.

For $(t,e) \in I_2$, we only have $q$, the number of chambers after $e$
in $t$; and $w$.  If $q_{(t_1,e_1)} = q_{(t_2,e_2)}$, and 
$w_{(t_1,e_1)} = w_{(t_2,e_2)}$ for $(t_i,e_i) \in I_2$, then
$S_{(t_1,e_1)} = S_{(t_2,e_2)}$.  Let $R^2_{q,w}$ be this set.

For fixed $p$, $w$, let $R^1_{p,w} = \cup_q R^1_{p,q,w}$, and for
fixed $w$ let $R_w^2 = \cup_q R^2_{q,w}$.  We will show that
$R^1_{p,w}$ and $R^2_w$ can usually be computed as the folding results 
of some half-infinite gallery type.  Let $r \in W$ be rotation by
$120^{\circ}$ counterclockwise, and let $f \in W$ be the reflection
about the horizontal line through $v_M$.  For $b$ with $\alpha + 2\beta > 0$,
and given $p \geq 1$ odd, we define $\Omega_{1}^2$, $\Omega_{f}^2$, 
$\Omega_{p,1}^1$ and $\Omega_{p,f}^1$ to be the half-infinite galleries
shown in Figure~\ref{Omegas} (This figure shows the specified
half-infinite galleries for a particular choice of $b$, namely
$\alpha = 4$ and $\beta = -1$.  But it is clear how the definition
would work for other $b$ with $\alpha + 2\beta > 0$.)
The half-infinite galleries $\Omega_{r}^2$, $\Omega_{rf}^2$, 
$\Omega_{p,r}^1$ and $\Omega_{p,rf}^1$ are rotations of those
in Figure~\ref{Omegas} by $240^{\circ}$ counterclockwise about the
center point of $C_M$, and 
$\Omega_{r^2}^2$, $\Omega_{r^2f}^2$, $\Omega_{p,r^2}^1$ and 
$\Omega_{p,r^2f}^1$ are rotations of those in Figure~\ref{Omegas}
by $120^{\circ}$ counterclockwise about the center point of $C_M$. 
If $\alpha + 2\beta = 0$ then $\Omega_{1}^2 = \Omega_{f}^2$,
$\Omega_{r}^2 = \Omega_{rf}^2$, and $\Omega_{r^2}^2 = \Omega_{r^2f}^2$ 
(the picture would be a degenerate version of Figure~\ref{Omegas}).
\begin{prop}
If $\alpha + 2\beta > 0$ or $b = 1$, then $R_{p,w}^1$ is the folding results
of $\Omega_{p,w}^1$, and $R_w^2$ is the folding results of
$\Omega_w^2$.  If $\alpha + 2\beta = 0$ and $b \neq 1$, then 
$(\cup_p R_{p,w}^1) \cup R_w^2$ is the union of the folding results of
the $\Omega_{p,w}^1$ for all $p$, and the folding results of $\Omega_w^2$.
It is still true that $R_w^2$ is the folding results of $\Omega_w^2$.
\end{prop}
Before we prove this proposition, we establish some more terminology.
We first let $\{ C_M ; t_1, t_2, \ldots , t_n \}$ be a finite gallery type.
Let $C_0 = C_M, C_1, \ldots , C_m$ be a (standard or non-standard) folding
of the sub-gallery type $\{ C_M ; t_1, t_2, \ldots , t_m \}$, where $m<n$.
Consider the finite gallery type $\{ C_m ; t_{m+1}, \ldots , t_n \}$,
and let $C_m = D_m$ , $D_{m+1}$ , $D_{m+2} , \ldots , D_n$ be the standard
folding.  If $D_{m+j}$, for $1 \leq j \leq n-m$ is on the same side of
the edge between $D_{m+j}$ and $D_{m+j-1}$ as $C_M$, then we say 
$m+j$ is a {\em choice point given the history $C_M , C_1 , \ldots, C_m$.}
Note that this all applies even if $m=0$.  Also note that if 
$1 \leq k \leq n$, then $k$ could be a choice point given some histories,
but not given others.  We now prove the proposition.
\begin{proof}
It suffices to consider $w = 1$ and $w = f$.  In fact, we will only 
consider $w = 1$, since the $w = f$ case is similar.  We first assume
that $\alpha + 2\beta > 0$ or $b = 1$, and that $(t,e) \in I_1$.  
Suppose $x \in SL_3(L)$ has SMG with type $t$ and
edge of departure $e$.  Also suppose $(t,e)$ has characteristics
$p$, $q$, and $w = 1$.  Let $\Gamma_x$ be the composite gallery.  So
in the notation of Section~\ref{SMGandComp}, 
$\Gamma_x = \Gamma_E^1 \cup \Gamma_E^2 \cup \Gamma_E^3$, where $E = xC_M$.
Consider the choice points given history $C_M = C_0$.  These all
lie in $\Gamma_E^3$, because of the shape of $\Gamma_x$ (illustrated
for $\alpha = 3$, $\beta = -1$, $q = 11$ and $p = 7$ in 
Figure~\ref{CompositeGalleryExample}).  The choice points given history
$C_M = C_0$ for the same $b$ and $x$ that were used in 
Figure~\ref{CompositeGalleryExample} are illustrated in 
Figure~\ref{GammaxChoicePoints}.  Note that $\Omega^1_{p,1}$ ends at the
same chamber that the standard folding of $\Gamma_x$ ends at, and note that
the choice points in $\Gamma_x$ that are after or at the turning
edge of $\Gamma_E^3$ are also choice points in $\Omega^1_{p,1}$.
This is illustrated for the example $\alpha = 3$, $\beta = -1$,
$w = 1$, $p = 7$ and $q = 11$ in Figure~\ref{LotsofChoicePoints}.
Note also that if $i$ is a choice point in $\Gamma_E^3$ that occurs
before the turning edge of $\Gamma_E^3$, and if $L_i$ is the wall 
containing the edge corresponding to $i$, then $L_i$ also passes
through $\Omega^1_{p,1}$.  Reflecting the section of $\Gamma_E^3$ that
occurs after $i$ about $L_i$ gives rise to a finite gallery type with
no choice edges additional to those obtained by reflecting the portion 
of $\Omega^1_{p,1}$ that occurs after $L_i$ about $L_i$.  This is illustrated 
in Figure~\ref{ChoiceBeforeTurning}.  
Therefore, the folding results of $\Omega^1_{p,1}$
contain $R^1_{p,1}$.
\begin{figure}
\centerline{\input{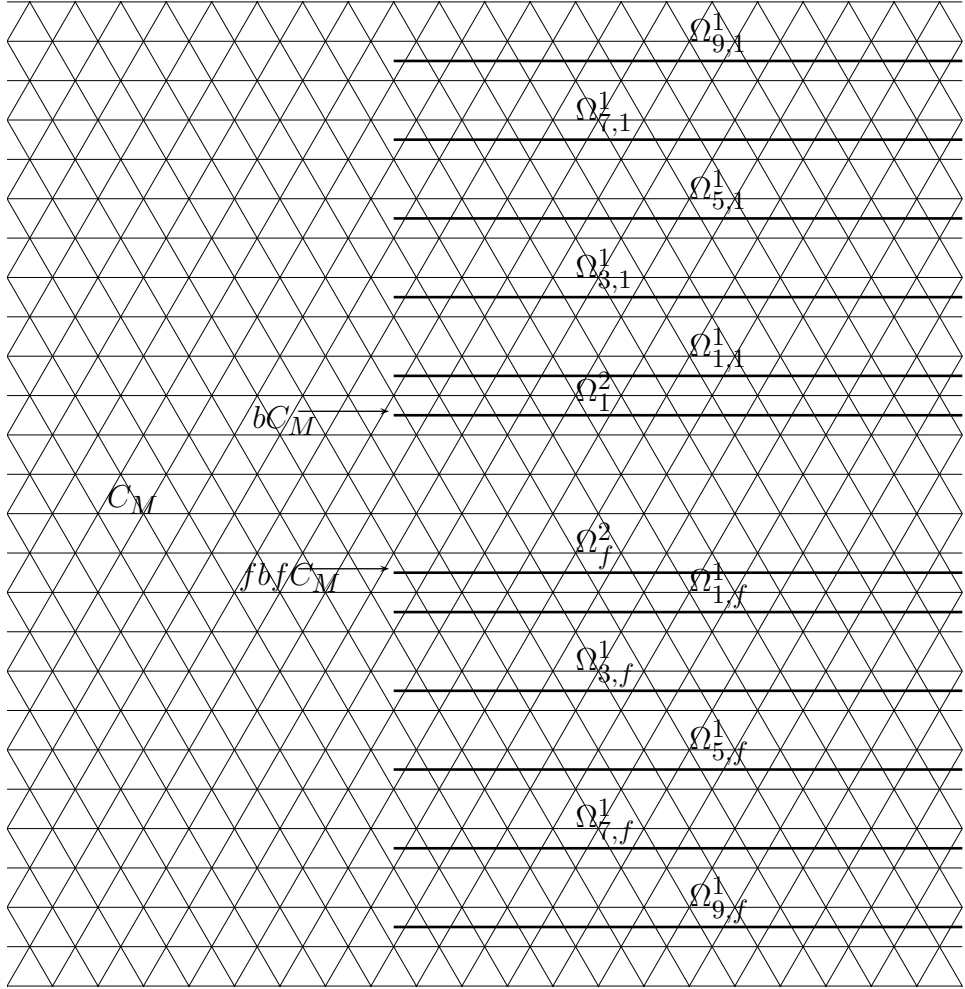}}
\caption{$\Omega_1^2$, $\Omega_f^2$, $\Omega_{p,1}^1$ and $\Omega_{p,f}^1$}
\label{Omegas}
\end{figure}
\begin{figure}
\centerline{\input{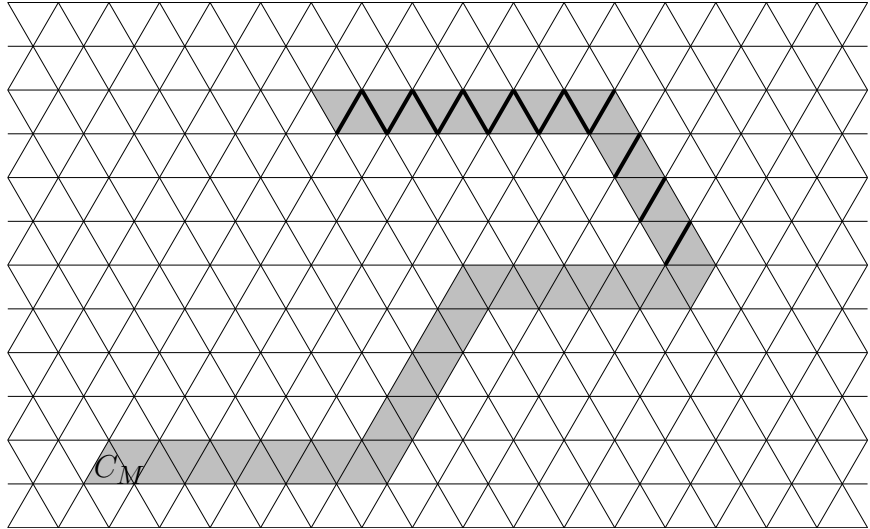}}
\caption{Choice points for a sample composite gallery}
\label{GammaxChoicePoints}
\end{figure}
\begin{figure}
\centerline{\input{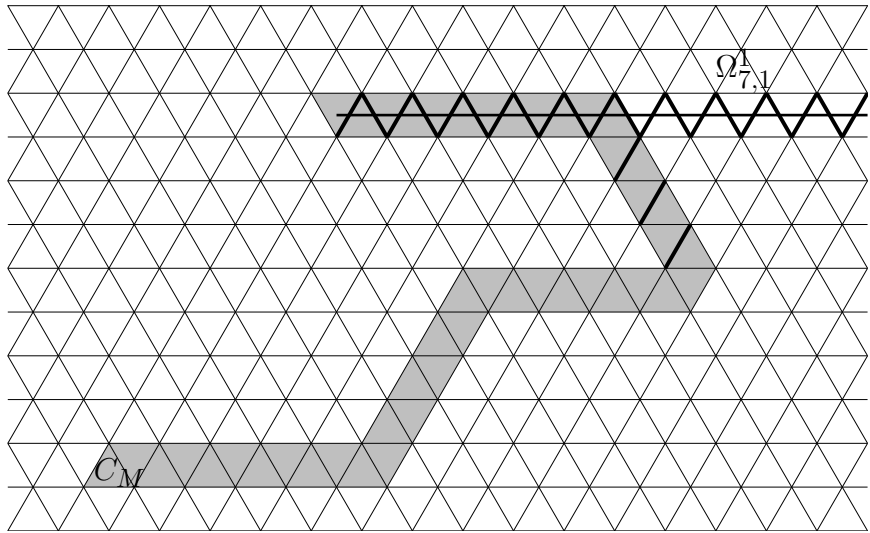}}
\caption{Choice points for $\Omega^1_{p,1}$ and $\Gamma_x$}
\label{LotsofChoicePoints}
\end{figure}
\begin{figure}
\centerline{\input{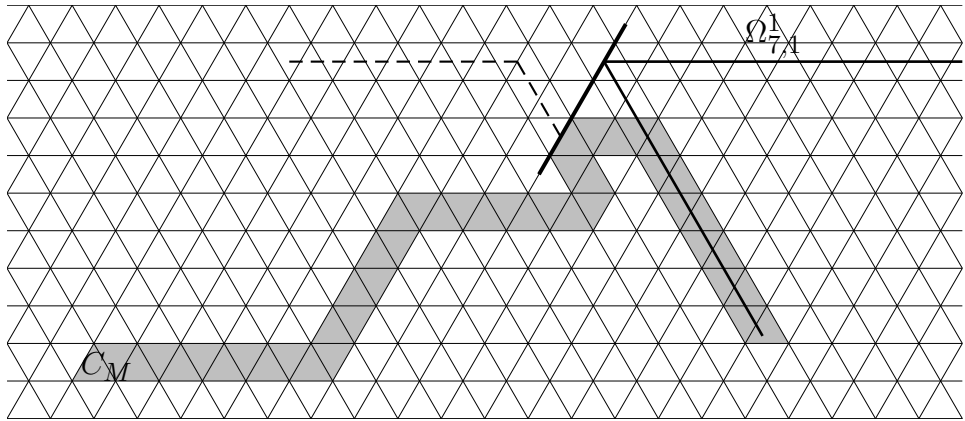}}
\caption{How to address choice points of $\Gamma_x$ that occur
before the turning point of $\Gamma_E^3$}
\label{ChoiceBeforeTurning}
\end{figure}

We must now show that $R^1_{p,1}$ contains the folding results of
$\Omega^1_{p,1}$.  It suffices to show that every choice point 
in $\Omega^1_{p,1}$ is also a choice point in some $\Gamma_x$
for some $xC_M$ with type-edge pair having characteristics 
$p$, $w = 1$, and some $q$.  But this can be accomplished by letting
$q$ get arbitrarily large.

To see that $R_1^2$ is the folding result of $\Omega_1^2$ for 
$\alpha + 2\beta \neq 2\alpha + \beta$ or for $b = 1$, proceed as
follows.  First note that $\Gamma_x$ for $x$ with $p = 0$ is minimal,
and the final gallery of $\rho(x^{-1}\Gamma_x)$ is $bC_M$.
This is also the final gallery of the standard folding of
$\Omega_1^2$.  Then let the chambers in the standard folding of 
$\Omega_1^2$ be called $\{ \ldots ,D_3,D_2,D_1\}$, where
$D_1 = bC_M$.  Let the wall between $D_i$ and $D_{i+1}$
be called $L_i$.  Let $G_i$ be the left-infinite gallery
in $A_M$ that one gets by reflecting the chambers $D_j$ for $j \leq i$
across $L_i$. Note that if $p \geq 1$, the only choice points in
$\Gamma_x$ are duplicated in $G_p$.  This shows that $R_1^2$ is equal
to the set of folding results of $\Omega_1^2$.

The $\alpha + 2\beta = 0$ for $b \neq 1$, $(t,e) \in I_1$
case is more complicated, but makes use of the
same ideas, and is therefore omitted.  The $2\alpha + \beta = 
\alpha + 2\beta$ case for $b \neq 1$, $(t,e) \in I_2$ is omitted 
for the same reason.
\end{proof}

This gives a more efficient way to compute the folding results of
all the type-edge pairs in $I_1 \cup I_2$ for any fixed $b$, and
also illustrates, to some extent, the source of the similarity
between the sets $S_1$ for different values of $b$.

\subsection{A Method Suggested by Rapoport and Kottwitz}\label{Rapoport}

The set $S_1$ of chambers in $A_M$ is, {\em a priori}, only a superset
of the solution set $S$ that we are seeking.  In this section we give the
results of a methodology suggested by Rapoport and Kottwitz for 
producing a subset $S_2 \subseteq S$.  In other words, this method will show
that certain chambers in $A_M$ must be contained in $S$.

The method in this section is only effective in the case $b = 1$.  In this 
case, techniques that will be described in the next section will serve
to enlarge $S_2$ to the point that it becomes equal to $S_1$, 
and therefore to $S$.  
In the $b \neq 1$ cases, we will use other
methods to arrive at a subset $S_2$ (see Section~\ref{SubsetGeometric}).

Let $a \in SL_3(F)$ be of the form
$$a = \left( 
\begin{matrix}
	\pi^m & 0 & 0 \\
	0 & \pi^n & 0 \\
	0 & 0 & \pi^{-m-n}
\end{matrix} \right),$$
where there are no conditions on $m$ and $n$.  Let $w$ be one of the following
matrices:
$$r = \left(
\begin{matrix}
	0 & 0 & 1 \\
	1 & 0 & 0 \\
	0 & 1 & 0 
\end{matrix} \right) \hspace{.25in} ; \hspace{.25in}
r^2 = \left(
\begin{matrix}
	0 & 1 & 0 \\
	0 & 0 & 1 \\
	1 & 0 & 0
\end{matrix} \right).$$
Therefore, $w$ represents a $3$-cycle in $W$, the finite Weyl group.
According to a result of Kottwitz, the matrix $aw$ belongs to a {\em basic}
$\sigma$-conjugacy class if there is some $l$ such that 
$aw \sigma(aw) \sigma^2 (aw) \ldots \sigma^{l-1}(aw)$ is central
in $SL_3(L)$ \cite{K1}, \cite{K2}.  See \cite{K1} and \cite{K2}
for a definition of {\em basic}, but for the present considerations,
one only needs to know that the $\sigma$-conjugacy class containing 
$b = 1$ is the
only basic $\sigma$-conjugacy class of $SL_3(L)$.

\begin{lem}\label{CentralityLemma}
There exists an $l$ such that 
$aw \sigma(aw) \sigma^2 (aw) \ldots \sigma^{l-1}(aw)$ is central for any
choice of $m$, $n$.
\end{lem}

\begin{proof}
First we note that $\sigma(aw) = aw$, since $aw \in SL_3(F)$.  So we want to 
show that $(aw)^l$ is central for some $l$.  But 
\begin{eqnarray*}
	(aw)^l & = & awaw(aw)^{l-2} \\
	& = & awaw^{-1}w^2(aw)^{l-2} \\
	& = & aa^ww^2(aw)^{l-2},
\end{eqnarray*}
where the superscript $w$ denotes $w$ acting on $a$ by conjugation.
Proceeding, we have
\begin{eqnarray*}
	aa^ww^2(aw)^{l-2} & = & aa^ww^2aw(aw)^{l-3} \\
	& = & aa^ww^2aw^{-2}w^3(aw)^{l-3} \\
	& = & aa^wa^{w^2}w^3(aw)^{l-3} \\
	& = & \ldots \\
	& = & aa^wa^{w^2}a^{w^3} \cdots a^{w^{l-1}}w^l.
\end{eqnarray*}
Choosing $l=3$, we get $aa^wa^{w^2}w^3 = 1$, which is central.

Since $w$ is a $3$-cycle, conjugating by $w$ serves to permute the 
diagonal entries of $a$ in a cyclic way.  This is why 
$aa^wa^{w^2}w^3$ is central.
\end{proof}

\begin{figure}
\setlength{\unitlength}{1in}
\begin{picture}(5,3.75)(0,0)
\centerline{\includegraphics[height=3.75in]{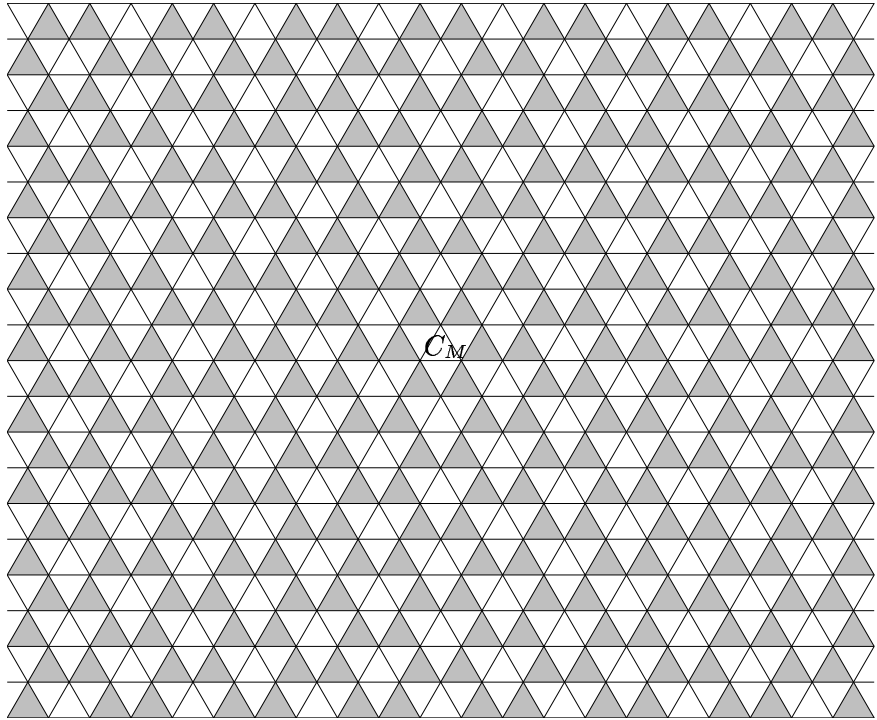}}
\end{picture}
\caption{Main results of the method suggested by Rapoport and Kottwitz}
\label{RapoportMethodResults}
\end{figure}

So the matrices $aw$ are in the $\sigma$-conjugacy class of $1$ for 
any choice of $m$, $n$.  Therefore, the double-$I$-cosets
$\{ IawI : m,n \in \mathbb{Z} \}$ all meet the identity $\sigma$-conjugacy
class non-trivially.  These double-$I$-cosets are pictured in 
Figure~\ref{RapoportMethodResults}.

With certain conditions on $m$ and $n$, it is also possible for $(aw)^l$
to be central for some $l$ if $w$ is not a $3$-cycle in $W$:

\begin{lem}\label{balllemma}
If $a = 1$ then there exists $l$ such that $(aw)^l$ is central for any 
$w \in W$.
\end{lem}

\begin{proof}
Take $l = 6$.
\end{proof}

\begin{figure}
\setlength{\unitlength}{1in}
\begin{picture}(3,2)(0,0)
\centerline{\includegraphics[height=2in]{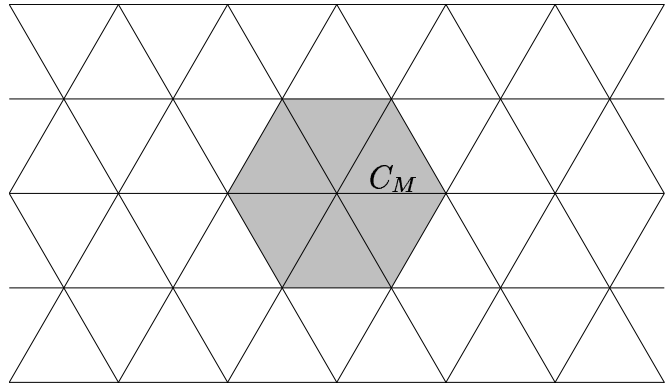}}
\end{picture}
\caption{The $a=1$ results of the method suggested by Rapoport and Kottwitz}
\label{RapoportMethodResultsII}
\end{figure}

The implication of this is that the double-$I$-cosets corresponding to
the chambers in Figure~\ref{RapoportMethodResultsII} all intersect the
$\sigma$-conjugacy class of $1$.  Let $f$ be as in Section~\ref{Superset}.

\begin{lem}\label{RapRes}
If $m = 0$ then $aa^w = 1$ for $w = f$. If $m = -n$ then $aa^w = 1$ 
for $w = rf$. If $n = 0$ then $aa^w = 1$ for $w = r^2 f$. 
\end{lem}

\begin{proof}
Just compute the relevant matrix products.
\end{proof}

\begin{figure}
\setlength{\unitlength}{1in}
\begin{picture}(4,4)(0,0)
\centerline{\includegraphics[height=4in]{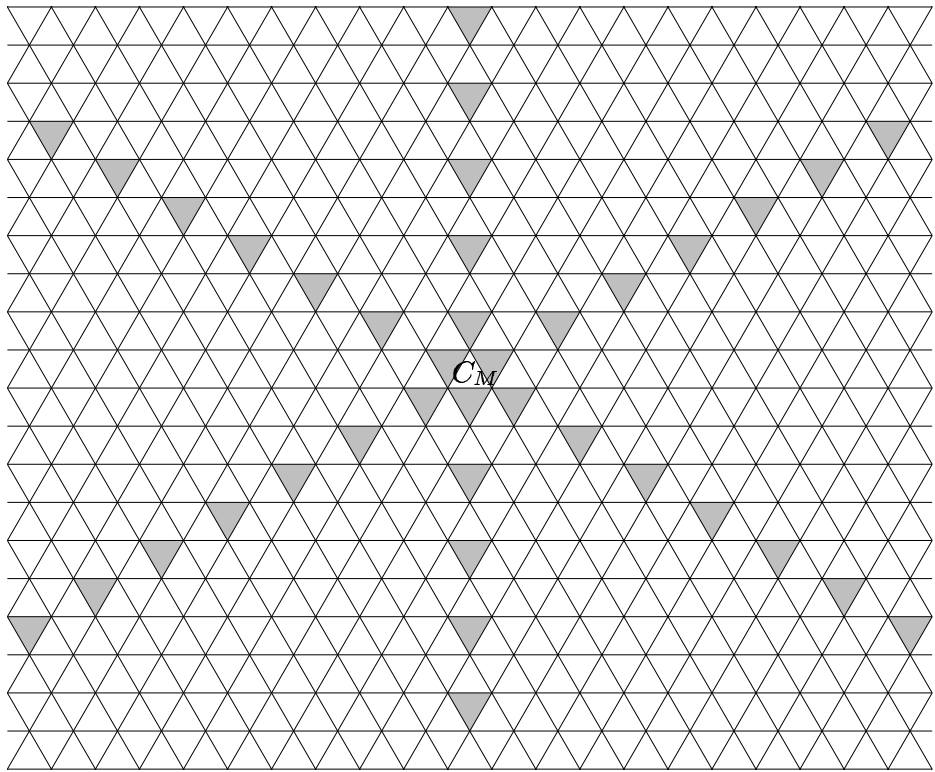}}
\end{picture}
\caption{Results from the method suggested by Rapoport and 
Kottwitz under the circumstances of Lemma~\ref{RapRes}}
\label{RapoportMethodResultsIII}
\end{figure}

The implication of this lemma is that the double-$I$-cosets
corresponding to the chambers in Figure~\ref{RapoportMethodResultsIII}
all intersect the $\sigma$-conjugacy class of $1$.  Note, in addition,
that we have a concrete element of the intersection of each of these
double-$I$-cosets with $\{ x^{-1} \sigma(x) : x \in SL_3(L) \}$.

It is easy to see that $(aw)^l$ is not central for any $l$, 
unless we are in one of the cases
enumerated above. 

It is clear how one would go about formulating and proving
the lemmas of this section for $SL_n$.  One would get results for
every possible cycle decomposition in $W = S_n$, and the ``most'' 
result chambers ($\mathbb{Z}^{n-1}$ worth of them) 
would arise from the $n$-cycles, since then $a$ is allowed
to be arbitrary.  At the other end of the spectrum, the finitely 
many Weyl group images of the main chamber arise from $a = 1$.
In subsequent sections, the results from this section will also be generalized
to $Sp_{2n}$ and $G_2$.

\subsection{A Subset of the Solution Set for $b=1$}\label{Subset}

Let $S_2$ be the collection of all double-$I$-cosets shown in the 
previous section to meet $\{ x^{-1} \sigma(x) : x \in SL_3 (L) \}$.
So $S_2$ is the union of the chambers pictured in 
Figures~\ref{RapoportMethodResults},~\ref{RapoportMethodResultsII}, 
and~\ref{RapoportMethodResultsIII}.  We know that 
$S_2 \subseteq S \subseteq S_1$, and in this section we enlarge
$S_2$ in such a way that it remains a subset of $S$.  The 
enlarged $S_2$ will turn out to be equal to $S_1$.
Note that this section concerns itself only with $b = 1$.

Let $w \in SL_3(L)$ be one of the matrices $r$, $r^2$ from the previous
section.  Let $a$ be as in the previous section.  Then we know 
$\tilde{b} = aw$ is $\sigma$-conjugate to $b = 1$, so 
$\{x^{-1} \sigma(x) : x \in SL_3(L) \} = 
\{x^{-1} \tilde{b} \sigma(x) : x \in SL_3(L) \}$, and also
$\{ I x^{-1} \sigma(x) I : x \in SL_3 (L) \} = 
\{I x^{-1} \tilde{b} \sigma(x) I : x \in SL_3(L) \}$.  But we know that the 
chamber in $A_M$ corresponding to $Ix^{-1} \tilde{b} \sigma(x)I$
is $\rho (x^{-1} \tilde{b} \sigma(x) C_M)$.  Just as in 
Section~\ref{SMGandComp}, let $\Gamma_x^1$ be the part of the SMG
from $C_M$ to $xC_M$ that is not in $A_M$, let $\Gamma_x^3 = 
\tilde{b}\sigma(\Gamma_x^1)$, and, if $e$ is the only edge of 
$\Gamma_x^1$ in $A_M$, let $\Gamma_x^2$ be a minimal gallery from $e$
to $\tilde{b}\sigma(e) = \tilde{b}e$.  If $\Gamma_x = 
\Gamma_x^1 \cup \Gamma_x^2 \cup \Gamma_x^3$, then the possible foldings 
of galleries of the same type as 
$x^{-1} \Gamma_x$ give possible candidates for additions to the set
$S_2$.  For some choices of $xC_M$ there is only one possible
folding of galleries of the same type as 
$x^{-1}\Gamma_x$, and the final chamber of this folding
is not already in $S_2$.  Such a situation would give an additional chamber
to add to $S_2$.

\begin{figure}
\centerline{\input{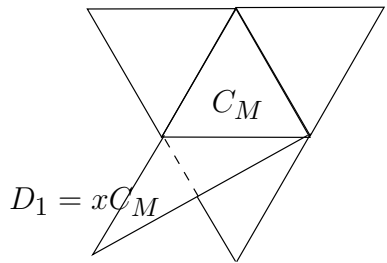}}
\caption{The choice $D_1$ of $xC_M$}
\label{BottomAppendage}
\end{figure}

\begin{figure}
\centerline{\input{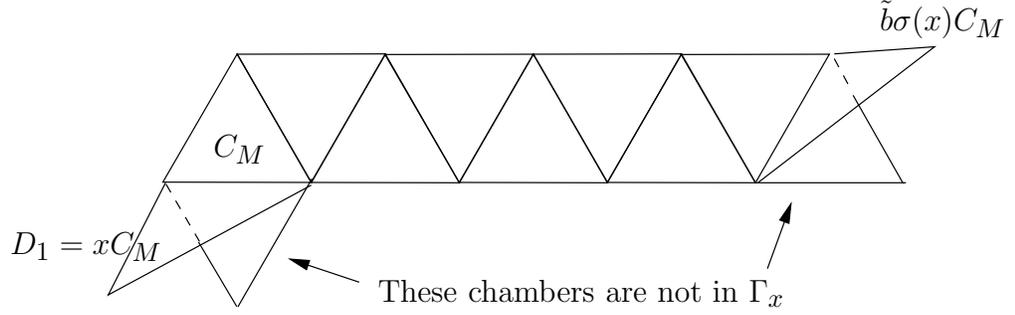}}
\caption{The resulting $\Gamma_x$ for $D_1$ and a particular choice of $aw$}
\label{ExampleAppendage}
\end{figure}

\begin{figure}
\centerline{\input{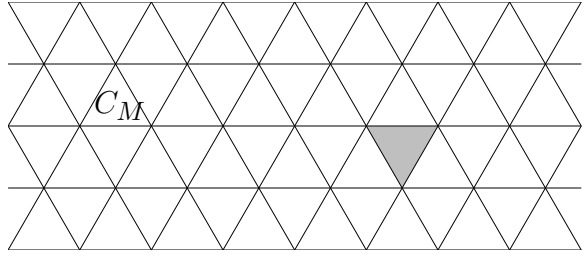}}
\caption{The result of folding galleries of the same type as 
$x^{-1}\Gamma_x$ for the $\Gamma_x$ pictured in Figure~\ref{ExampleAppendage}}
\label{ResultsofExample}
\end{figure}

For instance, we can choose $xC_M$ to be the chamber in 
Figure~\ref{BottomAppendage} that is adjacent to $C_M$, but 
not in $A_M$.
If $w = r^2$ and $$a = \left(
\begin{matrix}
	\pi^3 & 0 & 0 \\
	0 & \pi^{-1} & 0 \\
	0 & 0 & \pi^{-2}
\end{matrix} \right),$$ then $\Gamma_x$ is shown in 
Figure~\ref{ExampleAppendage}.
Any gallery of the same type as $x^{-1}\Gamma_x$ folds down into $A_M$ in the
same way, so we know that $Ix^{-1}\tilde{b}\sigma(x)I$ corresponds to
the chamber shown in Figure~\ref{ResultsofExample}, which was not
obtained by the methods of the previous section.

\begin{figure}
\setlength{\unitlength}{1in}
\begin{picture}(4,4.25)(0,0)
\centerline{\includegraphics[height=4.25in]{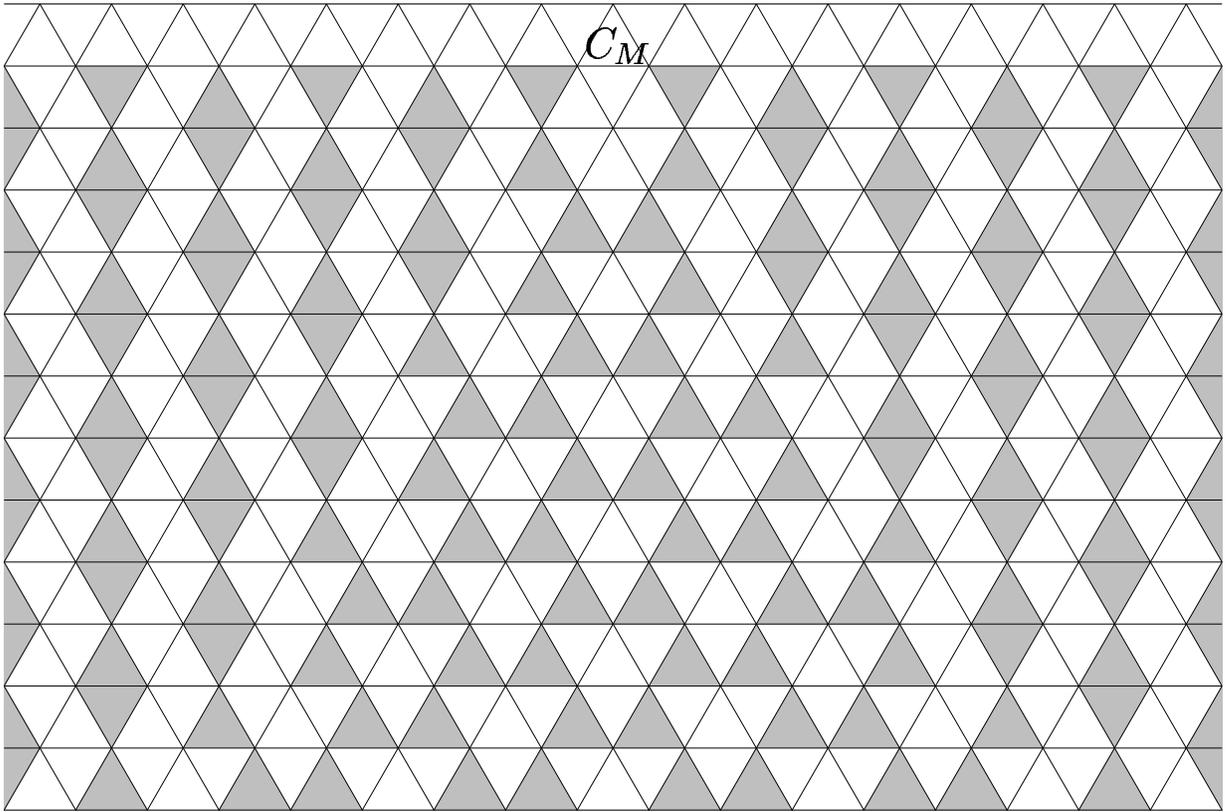}}
\end{picture}
\caption{Results for $xC_M = D_1$ and $w = r^2$ or for $xC_M = D_1$ and
$w = r$}
\label{D1W1}
\end{figure}

We now keep $w$ as above, but let $a$ take the general form $$a = \left(
\begin{matrix}
	\pi^m & 0 & 0 \\
	0 & \pi^n & 0 \\
	0 & 0 & \pi^{-m-n}
\end{matrix} \right).$$  
For these $\tilde{b} = aw$ we can compute the possible
ways that galleries of the same type as $x^{-1}\Gamma_x$ 
can fold into $A_M$, where, as before, $x$ is 
chosen so that $xC_M$ is the chamber $D_1$ 
pictured in Figure~\ref{BottomAppendage}.
One can easily see that for each $m$, $n$, there is only one way that
such galleries can fold down.  The chambers obtained through this
process are marked in Figure~\ref{D1W1}.  Many of these did not arise from the
methods of the previous section.

\begin{figure}
\centerline{\input{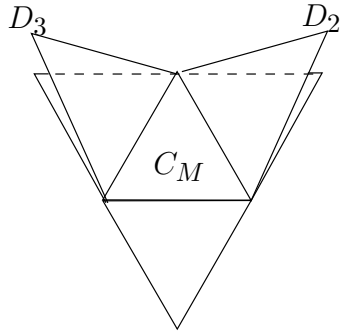}}
\caption{The choices $D_2$ and $D_3$ of $xC_M$}
\label{RightAppendageLeftAppendage}
\end{figure}

We could also choose $x$ such that $xC_M$ is one of the chambers
$D_2$, $D_3$ in Figure~\ref{RightAppendageLeftAppendage}, and we could choose
$w = r$.  For each of these possible choices, we could make considerations
similar to the above.  The results of this process are pictured
in subsequent figures according to Table~\ref{ResultReferenceTable}.
\begin{table}
\centerline{
\begin{tabular}{|c|c|c|} 
\hline
$w$ & $xC_M$ & Figure containing results \\
\hline
\hline
$r^2$ & $D_1$ & Figure~\ref{D1W1} \\
\hline
$r^2$ & $D_2$ & Figure~\ref{D2W1} \\
\hline
$r^2$ & $D_3$ & Figure~\ref{D3W1} \\
\hline
$r$ & $D_1$ & Figure~\ref{D1W1} \\
\hline
$r$ & $D_2$ & Figure~\ref{D2W1} \\
\hline
$r$ & $D_3$ & Figure~\ref{D3W1} \\
\hline
\end{tabular}}
\caption{Which figure contains which results}
\label{ResultReferenceTable}
\end{table}
\begin{figure}
\setlength{\unitlength}{1in}
\begin{picture}(4,4.25)(0,0)
\centerline{\includegraphics[height=4.25in]{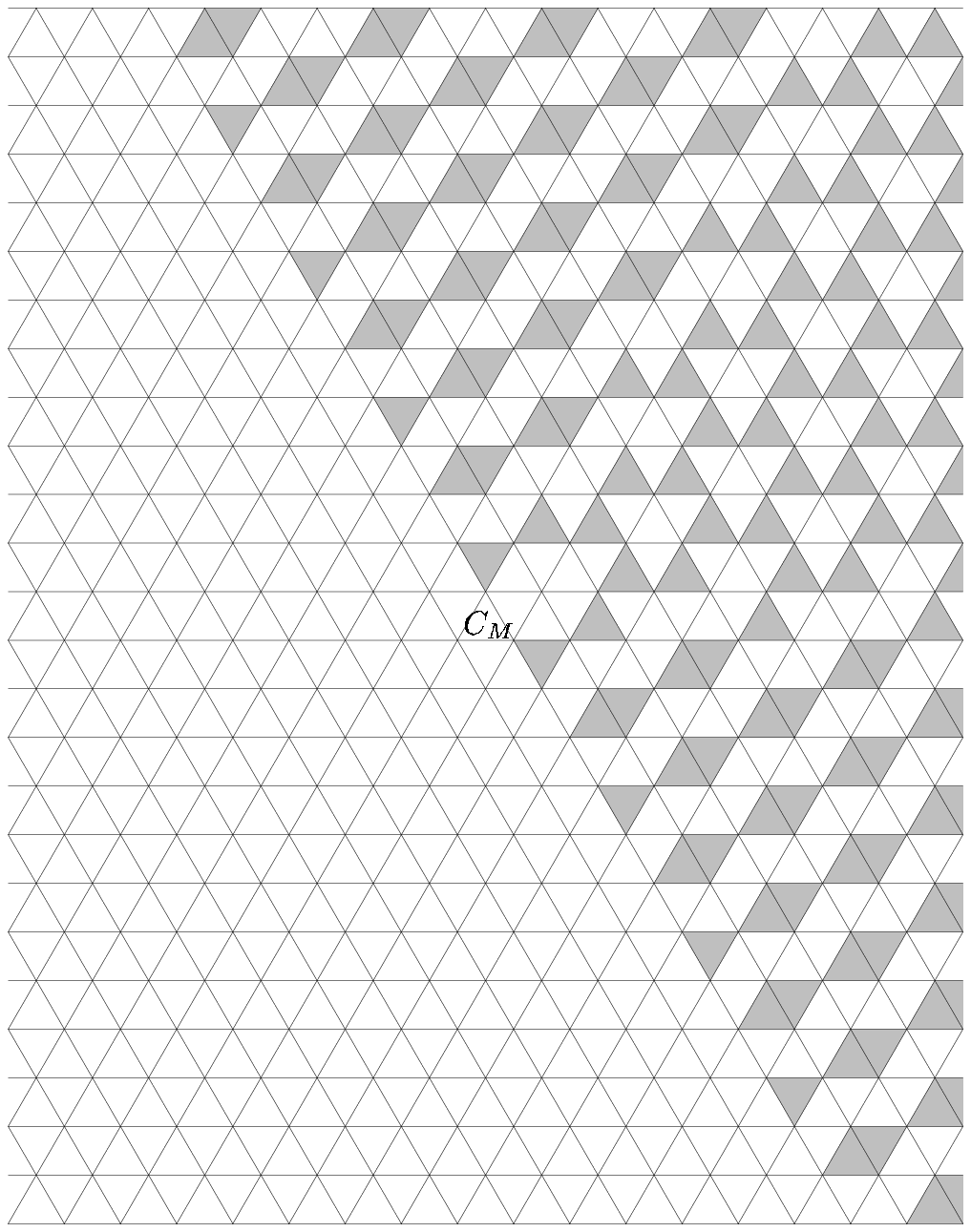}}
\end{picture}
\caption{Results for $xC_M = D_2$ and $w = r^2$ or for $xC_M = D_2$ and
$w = r$}
\label{D2W1}
\end{figure}
\begin{figure}
\centerline{\input{fig46.tex}}
\caption{Results for $xC_M = D_3$ and $w = r^2$ or for $xC_M = D_3$ and
$w = r$}
\label{D3W1}
\end{figure}

\begin{figure}
\setlength{\unitlength}{1in}
\begin{picture}(6,4)(0,0)
\centerline{\includegraphics[height=4in]{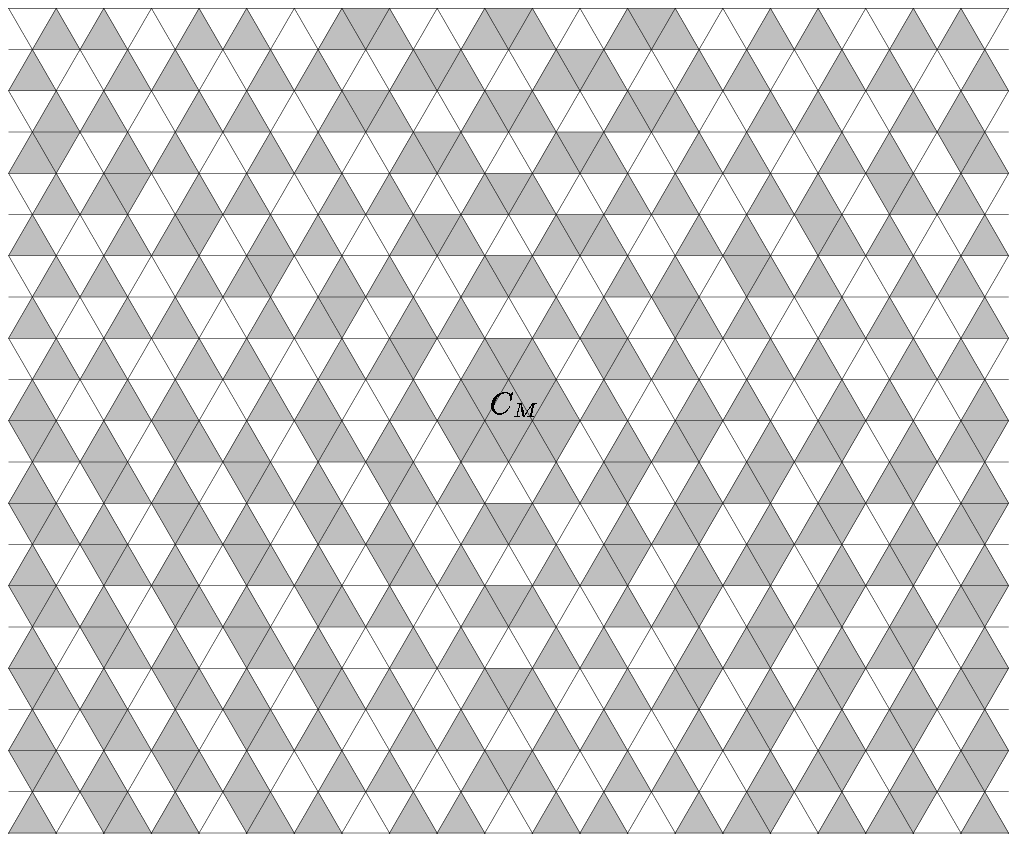}}
\end{picture}
\caption{Combination of all results referred to in 
Table~\ref{ResultReferenceTable}}
\label{FirstS2Bar}
\end{figure}

If we combine all these results with those of the previous section, we get 
the chambers pictured in Figure~\ref{FirstS2Bar}.  
Note that this combined set of
chambers is still not equal to $S_1$.

\begin{figure}
\centerline{\input{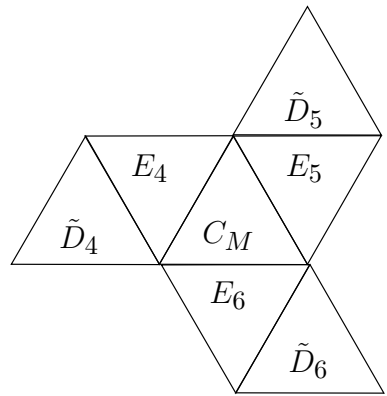}}
\caption{Some choices for $xC_M$ that are two chambers away from $C_M$}
\label{LengthTwoAppendages}
\end{figure}

If $\tilde{D}_4$, $\tilde{D}_5$, and $\tilde{D}_6$ are the chambers in $A_M$ 
labeled 
in Figure~\ref{LengthTwoAppendages}, and if $g_4 , g_5, g_6 \in I$ are chosen
such that $g_i E_i \neq E_i$, then we can choose $x$ such that $xC_M$ is any
of the $D_i = g_i \tilde{D}_i$.  For each of these $xC_M$, we could make 
considerations similar to the above for $w = r$, $r^2$.  The results of this
process include all the chambers in Figure~\ref{TotalResults0_0} 
that are not present in Figure~\ref{FirstS2Bar}.
This proves that the chambers in Figure~\ref{TotalResults0_0}
represent the exact collection of double-$I$-cosets that intersect
non-trivially with the $\sigma$-conjugacy class of $1$. 

Rather than 
using $xC_M$ a distance of two chambers from $C_M$, it is also
possible to obtain the chambers in Figure~\ref{TotalResults0_0}
that are not in Figure~\ref{FirstS2Bar} in the following way.
Let $\tilde{b}$ now be such that $\tilde{b}C_M$ is one of the chambers
in Figure~\ref{FirstS2Bar} that was not obtained using the methods 
of the previous section.  Let $xC_M \nsubseteq A_M$ be $D_1$, $D_2$ or 
$D_3$, and compute as previously.  This gives rise to
the chambers in Figure~\ref{TotalResults0_0} that are not in 
Figure~\ref{FirstS2Bar}.  Note that this constitutes an iteration
of our ``length one appendage'' methods, and gives rise to the 
same additional chambers that our ``length two appendage'' methods did.

The methods in this section assume the existence of a large and 
well-distributed set $S_2$.
If one has a smaller $S_2$, the enlargement effort will be less
effective.  One could apply the methods of this section to 
$\sigma$-conjugacy classes other than $b = 1$, but one would have to somehow
come up with a reasonably large $S_2$ first.  

\subsection{A Geometric Construction of a Subset of the Solution 
Set}\label{SubsetGeometric}

In this section, we prove using geometric methods that the chambers in 
the superset $S_1$ produced in Section~\ref{Superset} that arise
from the infinite classes $I_1$ and $I_2$ are all actually contained 
in the solution set $S$.  Since the classes $I_1$ and $I_2$ seemed to give
rise to all chambers in $S_1$, this proves that $S_1 = S$.

\begin{figure}
\centerline{\input{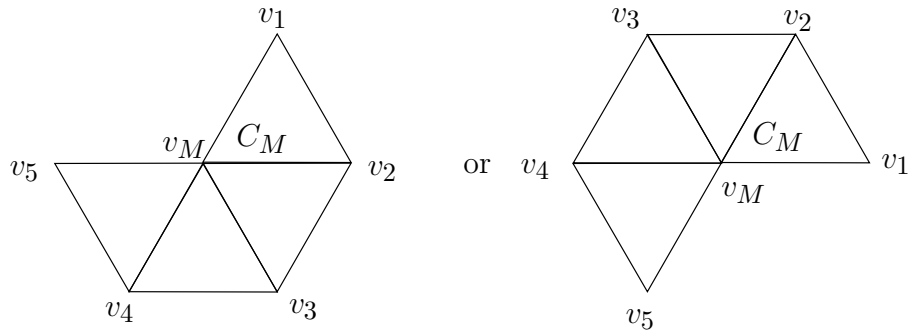}}
\caption{A Useful Lemma}
\label{UsefulLemma}
\end{figure}

We begin with the following very useful lemma.
\begin{lem}\label{UsefulLemmaStatement} 
Let $G$ be one of the galleries pictured in Figure~\ref{UsefulLemma}.  
Let $\tilde{G}
= gG$, where $g \in GL_3(L)$ is arbitrary.  There exists a unique vertex 
$\tilde{v}$ such that $\tilde{v}$ is adjacent to each of $gv_5$,
$gv_1$, and $gv_M$.  Here, $v_M$ is the main vertex in the building.
\end{lem}
\begin{proof}
It suffices to work with $g = 1$.  The vertices adjacent to $v_M$
correspond to non-trivial proper subspaces of $(\overline{\mathbb{F}}_q)^3$,
where $\mathbb{F}_q$ is the residue field of $F$ 
(so $\overline{\mathbb{F}}_q$ is the residue field of $L$).  For such a 
vertex $v$, let $v$ also denote the corresponding subspace of 
$(\overline{\mathbb{F}}_q)^3$.  If $v$ and $w$ are two vertices adjacent
to $v_M$, then $v$ and $w$ are adjacent to each other if and only if
$v \subseteq w$ or $w \subseteq v$ as subspaces.  We therefore have the 
following cases:

\noindent {\em Case 1:} $\dim(v_1) = 1$.  Then $\dim(v_2) = 2$, 
$\dim(v_3) = 1$, $\dim(v_4) = 2$, $\dim(v_5) = 1$, and we can and must
choose $\tilde{v} = \langle v_1 , v_5 \rangle$.

\noindent {\em Case 2:} $\dim(v_1) = 2$.  Then $\dim(v_2) = 1$, 
$\dim(v_3) = 2$, $\dim(v_4) = 1$, $\dim(v_5) = 2$, and we can and must
choose $\tilde{v} = v_1 \cap v_5$.
\end{proof}
\noindent The idea of this 
lemma is that whenever we see an arrangement of 
chambers shaped like one of those in Figure~\ref{UsefulLemma}, we can fill it
in uniquely with two more chambers to create a hexagon.

We also need the following definition:
\begin{defn}
Given a minimal gallery $G$, the {\em parallelogram}, $P(G)$ associated
with $G$ is $P(G) = \cap A$, where the intersection is over all apartments 
$A$ containing $G$.
\end{defn}\label{parallel}
\begin{thm}
$P(G) = \cup \tilde{G}$, where the union is over all minimal galleries 
$\tilde{G}$ that stretch from the first chamber of $G$ to the last
chamber of $G$.
\end{thm}\label{paralleltheorem}
\begin{proof}
See \cite{K3}.
\end{proof}

\noindent Note that $P(G)$ is not necessarily actually a 
parallelogram.  It is a parallelogram with possibly 
a single chamber removed from either
or both of the acute angle corners.

We now discuss the the program that we will follow in the rest of this
section.  We will begin by restricting our attention to the collection
of $\sigma$-conjugacy classes consisting of those $b$ with 
$\alpha + 2\beta \neq 0$, and $b = 1$.  We also restrict our attention 
to type-edge pairs in $I_1$.  If $x \in SL_3(L)$ gives rise to a 
type-edge pair $(t,e) \in I_1$, and if $b$ is restricted as specified,
then we will define two invariants $\gamma_1$ and $\gamma_2$ of $x$.
It will turn out that $\gamma_1$ and $\gamma_2$ determine 
$\rho(x^{-1}\Gamma_x)$.  
We will also show that given any $\gamma_1$ and $\gamma_2$,
we can choose $x$ such that $x$ and $b$ give rise to $(t,e)$, $\gamma_1$,
and $\gamma_2$.  This proves that for $b = 1$ or $b$ with 
$\alpha + 2\beta \neq 0$, all the results in $S_1$ which come from $I_1$
are in fact in $S$.  We then use similar methods, with a few added
complications, to address $b \neq 1$ with $\alpha + 2\beta =0$,
still focusing our attention on $(t,e) \in I_1$.  We then turn
to $(t,e) \in I_2$, and we consider separately the cases $b = 1$ 
or $\alpha + 2\beta \neq 2\alpha + \beta$, and $b \neq 1$ with
$\alpha + 2\beta = 2\alpha + \beta$.

\begin{figure}
\centerline{\input{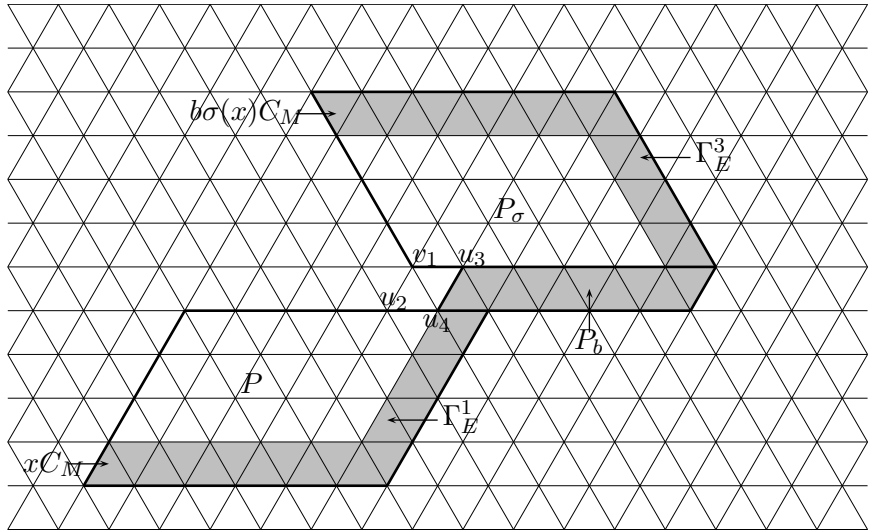}}
\caption{An example of $P$, $P_b$ and $P_\sigma$}
\label{PPbPsigma}
\end{figure}

\begin{figure}
\centerline{\input{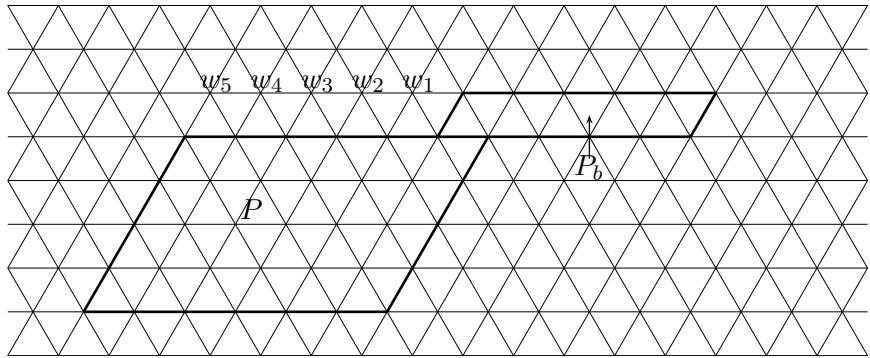}}
\caption{$w_2$, $w_3$, $w_4$ and $w_5$}
\label{w2345}
\end{figure}

As we just mentioned, we start by assuming $(t,e) \in I_1$ and either
$b = 1$ or $\alpha + 2\beta \neq 0$. Recall that the associated 
composite gallery $\Gamma_x$ is $\Gamma_E^1 \cup \Gamma_E^2 \cup
\Gamma_E^3$, where each $\Gamma_E^i$ is minimal (if $b = 1$ then
$\Gamma_E^2 = \emptyset$).  Let $P = P(\Gamma_E^1)$,
let $P_b = P(\Gamma_E^2)$, and let $P_\sigma = P(\Gamma_E^3)$.
It is worth noting that although we never made a specific choice of 
$\Gamma_E^2$, $P(\Gamma_E^2)$ is well defined (this relies on the 
fact that we have either $b = 1$ or $b$ with $\alpha + 2\beta \neq 0$).  
Figure~\ref{PPbPsigma} has a picture of $P$, $P_b$, and 
$P_\sigma$ for $$b = \left(
\begin{matrix}
	\pi^3 & 0 & 0 \\
	0 & \pi^{-1} & 0 \\
	0 & 0 & \pi^{-2}
\end{matrix}
\right),$$
and the example type-edge pair pictured in Figure~\ref{GammaD1}.
By Lemma~\ref{UsefulLemmaStatement}, we can choose a vertex $w_1$ adjacent to 
$u_2$, $u_3$, and $u_4$.  This vertex is not necessarily $v_1$, although
it may be.  By
repeated application of Lemma~\ref{UsefulLemmaStatement}, 
we can choose the vertices 
$w_2$, $w_3$, $w_4$ and $w_5$, labeled on Figure~\ref{w2345}.

\begin{figure}
\centerline{\input{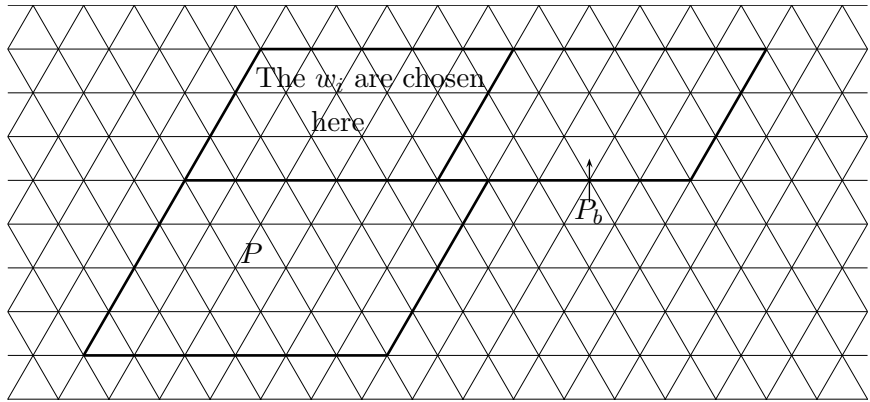}}
\caption{Multiple rows of the $w_i$}
\label{MultipleRowsofwi}
\end{figure}

\begin{figure}
\centerline{\input{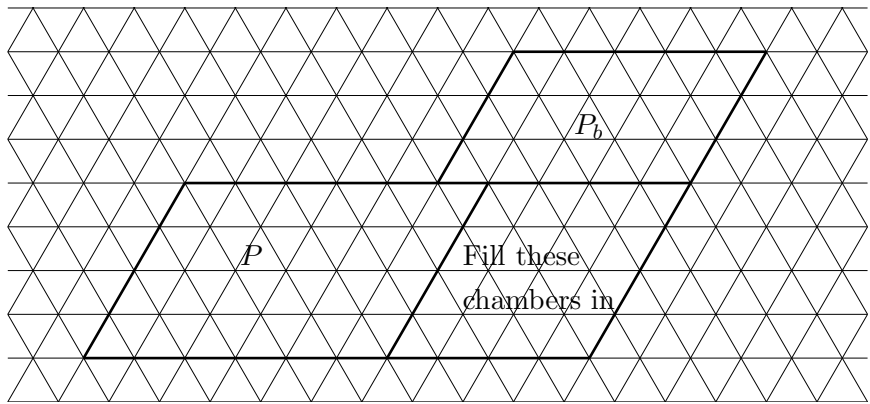}}
\caption{Filling in the rest of $\overline{P}$}
\label{RestofPBar}
\end{figure}

For an arbitrary $b \neq 1$ with $\alpha + 2\beta \neq 0$, 
and for any type-edge pair in $I_1$, we can choose
vertices analogous to the $w_i$ by repeated application of 
Lemma~\ref{UsefulLemmaStatement}.  Note that if $P_b$ has more than one row of
chambers, then we must choose correspondingly many rows of the $w_i$, one
row at a time, from the bottom up.  See Figure~\ref{MultipleRowsofwi}.
We can also use Lemma~\ref{UsefulLemmaStatement} repeatedly to
fill in chambers as labeled in Figure~\ref{RestofPBar}.  We refer to the
resulting unique construction consisting of $P$, $P_b$, and the chambers
added in the two above processes as $\overline{P}$.  Note that 
$\overline{P}$ is the intersection of all apartments containing
$P$ and $P_b$, and is also the union of all minimal galleries from the 
first chamber of $\Gamma_E^1$ to the last chamber in $\Gamma_E^2$. 

We now define two invariants which describe the way in which $\overline{P}$
is connected to $P_\sigma$.
\begin{defn}
Let $\gamma_1$ denote the number of edges that $\overline{P}$ and
$P_\sigma$ have in common.
\end{defn}
\noindent So $\gamma_1 \geq 1$, and in the 
case pictured in Figure~\ref{PPbPsigma},
$1 \leq \gamma_1 \leq 6$.

\begin{figure}
\centerline{\input{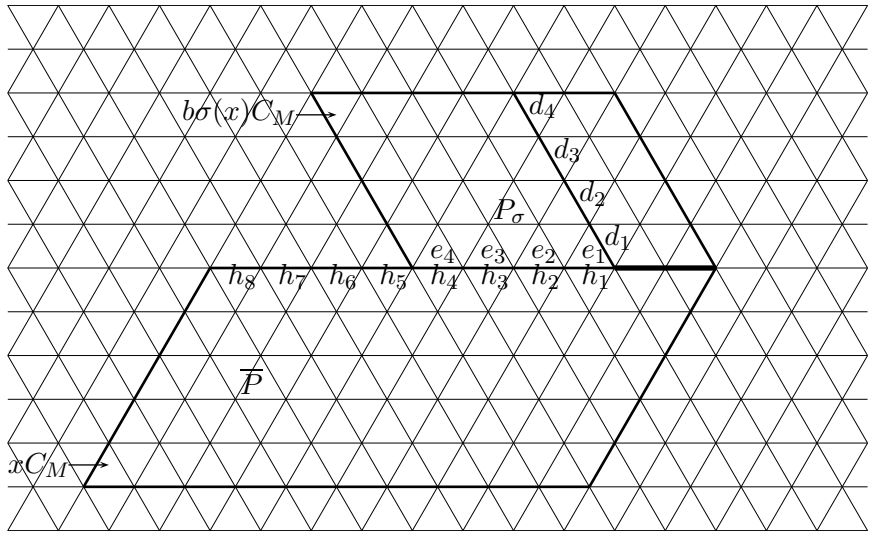}}
\caption{The definition of the invariant $\gamma_2$, part 1}
\label{Definitiongamma2}
\end{figure}

\begin{figure}
\centerline{\input{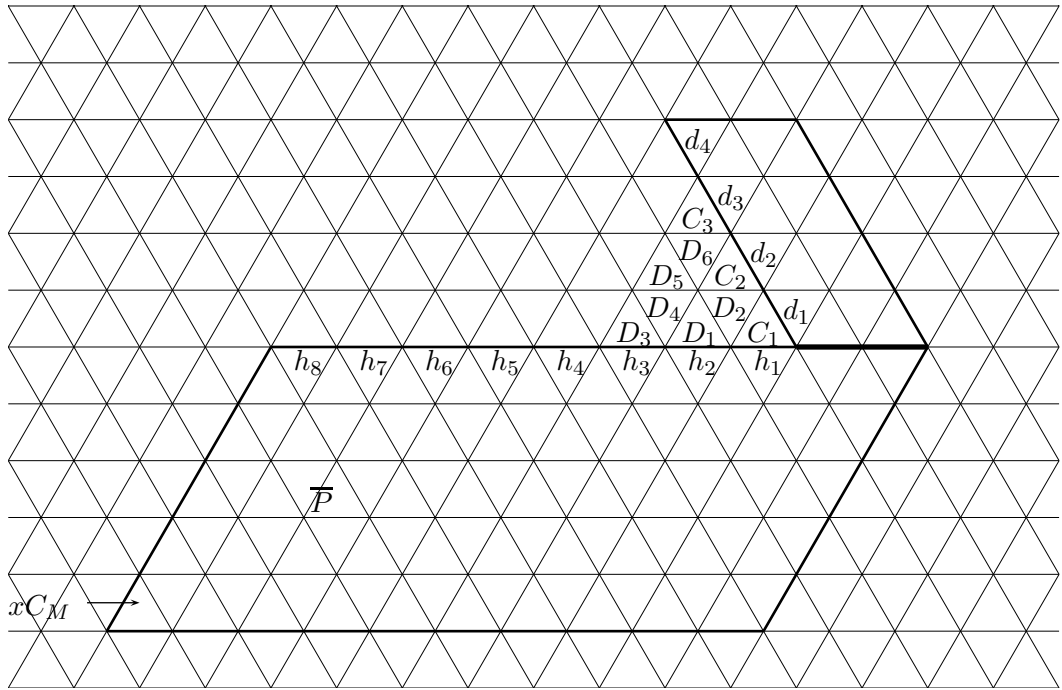}}
\caption{The definition of the invariant $\gamma_2$, part 3}
\label{Definitiongamma23}
\end{figure}

We define the second invariant $\gamma_2$ using 
Figure~\ref{Definitiongamma2}, which shows a situation with $\gamma_1 = 2$.
Note that for each $i = 1$,$2$,$3$,$4$, the edges $e_i$ of $P_{\sigma}$ 
and $h_i$ of $\overline{P}$ are
different, although they look the same in the figure.  
We pass to Figure~\ref{Definitiongamma22}, which is the same as 
Figure~\ref{Definitiongamma2}, only part of $P_\sigma$ is not shown for
clarity.  We now define $\gamma_2$ for our example situation.  It will be
clear from this example how $\gamma_2$ is defined in general.  If there is
a chamber $C_1$ (pictured in Figure~\ref{Definitiongamma23}) having
$h_1$ and $d_1$ as edges, then it is unique, and not contained 
in $P_\sigma$.  In
this case we require $\gamma_2 \geq 1$, and we get $D_1$ and $D_2$ as pictured
(in Figure~\ref{Definitiongamma23}) by using Lemma~\ref{UsefulLemmaStatement}.
If there is no $C_1$ as described then we say $\gamma_2 = 0$.  Given $C_1$, 
if there
is a chamber $C_2$ connecting $D_2$ to $d_2$, then it is unique.  In this
case we require $\gamma_2 \geq 2$, and we get $D_3$, $D_4$, and $D_5$,
$D_6$ by using Lemma~\ref{UsefulLemmaStatement} twice.  If there is no $C_2$
as described then we say $\gamma_2 = 1$.  Given $C_2$, we look
for $C_3$ as labeled.  We continue in this manner, with the restriction
that $\gamma_2$ cannot be bigger than either the number of the $h_i$ or the
number of the $d_i$.  If $\gamma_2$ is equal to the lesser of these numbers,
then we say $\gamma_2$ is maximal.
\begin{figure}
\centerline{\input{fig55.tex}}
\caption{The definition of the invariant $\gamma_2$, part 2}
\label{Definitiongamma22}
\end{figure}

For any parallelogram $P$, let
$T(P)$ be the number of edges along the top of $P$.  Let $S(P)$
be the number of edges along the right side of $P$.  Recall that
this could be different from the number of edges along the left side
of $P$.  Similarly, $T(P)$ could be different from the number of edges
along the bottom of $P$. This is because the ``parallelograms'' 
with which we are working are
not true parallelograms in the traditional sense. They are parallelograms that
may have one chamber removed from either or both of the acute angle
corners.  Although note that the upper right acute angle of $P$ (and
therefore $\overline{P}$) never has
a chamber removed; $P_b$ will always be a parallelogram in the traditional
sense of the word; and the lower right acute angle of $P_{\sigma}$ 
never has a chamber removed.

Fixing $b$ and a type-edge pair $(t,e) \in I_1$, 
let $\gamma_1^{\textrm{max}} = T(P)$.  Choose $\gamma_1$
such that $1 \leq \gamma_1 \leq \gamma_1^{\textrm{max}}$, and choose $\gamma_2$
such that $0 \leq \gamma_2 \leq \min \{ S(P) , T(\overline{P}) - \gamma_1 \}
= \gamma_2^{\textrm{max}}(\gamma_1)$.  We will show
that one can find $x$ such that $xC_M$ gives the type-edge pair
$(t,e)$, and the invariants $\gamma_1$ and $\gamma_2$.  We begin by making 
a simplification.  

If $b$ is still fixed; 
$(t_1 ,e_1 )$ and $(t_2 , e_2 )$ are two type-edge pairs in $I_1$;
$P_1$ and $P_2$ are the corresponding parallelograms $P(\Gamma_{E_1}^1)$
and $P( \Gamma_{E_2}^2 )$; and if $T(P_1) \leq T(P_2)$ and 
$S(P_1) \leq S(P_2)$, then $\gamma_1^{\textrm{max}}(P_1) \leq 
\gamma_1^{\textrm{max}}(P_2)$.  Given 
$\gamma_1 \leq \gamma_1^{\textrm{max}} (P_1)$, we also have 
$\gamma_2^{\textrm{max}}(P_1,\gamma_1) 
\leq \gamma_2^{\textrm{max}} (P_2,\gamma_1)$.
To construct $P_1$ with invariants
$\gamma_1 (P_1) \leq \gamma_1^{\textrm{max}} (P_1)$ and 
$\gamma_2 (P_1) \leq \gamma_2^{\textrm{max}} (P_1, \gamma_1(P_1))$, 
it would suffice
to construct $P_2$ with the same invariants.  One could then chop off
part of $P_2$ to get the desired $P_1$.  Therefore it suffices to
demonstrate the construction of an ``infinite parallelogram''
$P$ with $T(P) = S(P) = \infty$ with any fixed (and still finite)
invariants $\gamma_1$ and $\gamma_2$.  The equation $T(P) = S(P) = \infty$
is meant to imply a parallelogram that is half-infinite along the top and
along the right side as in Figure~\ref{HalfInfiniteParallelogram}
\begin{figure}
\centerline{\input{fig57.tex}}
\caption{$T(P) = S(P) = \infty$}
\label{HalfInfiniteParallelogram}
\end{figure}

\begin{figure}
\centerline{\input{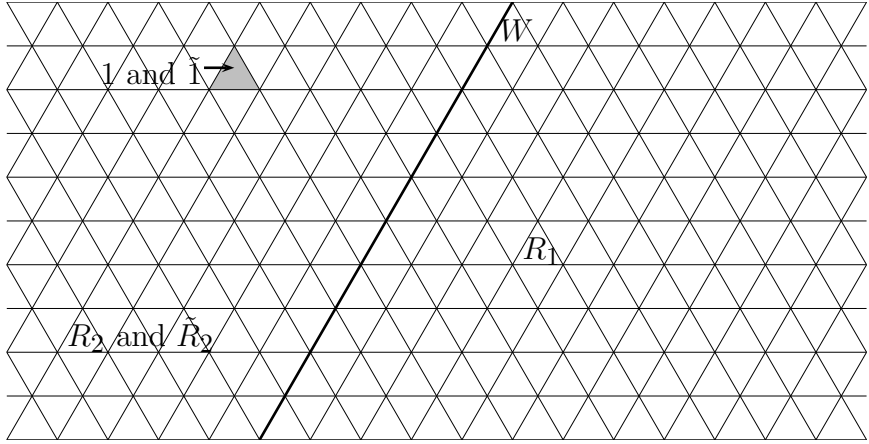}}
\caption{The meaning of a tilde}
\label{TwidleNotation}
\end{figure}

We now discuss notation to be used in constructions to come.  Regions of
an apartment of the building will always be denoted using Roman letters
with subscript.  Single chambers will be labeled using Roman letters
with subscript, or just using numbers.  At several points we will have to
graphically represent structures in the building that cannot be embedded into
a Euclidean plane.  For instance, we will discuss the main apartment $A_M$,
together with a half apartment $A_1^{\frac{1}{2}}$ coming out of a wall
$K$ in $A_M$.  We would picture this example as in Figure~\ref{TwidleNotation},
where $R_1$ is a region in $A_M$, $R_2$ is also in $A_M$, and $\tilde{R}_2 
= A_1^{\frac{1}{2}}$.  Thus $R_2$ and $\tilde{R}_2$ are the same
region graphically, but different regions in the construction we are
describing.  The chamber shaded in Figure~\ref{TwidleNotation} could
be denoted $1$ if in $R_2$ or $\tilde{1}$ if in $\tilde{R}_2$.  The tilde
will be used generally in this way, with constructions marked with a
tilde understood to not be in $\overline{P}$, $P_b$, or $P_{\sigma}$,
unless specifically included.

We now proceed to demonstrate that one can find $x$ such that $xC_M$ 
gives $(t,e)$, $\gamma_1$, and $\gamma_2$.  We work with two cases.

\begin{figure}
\centerline{\input{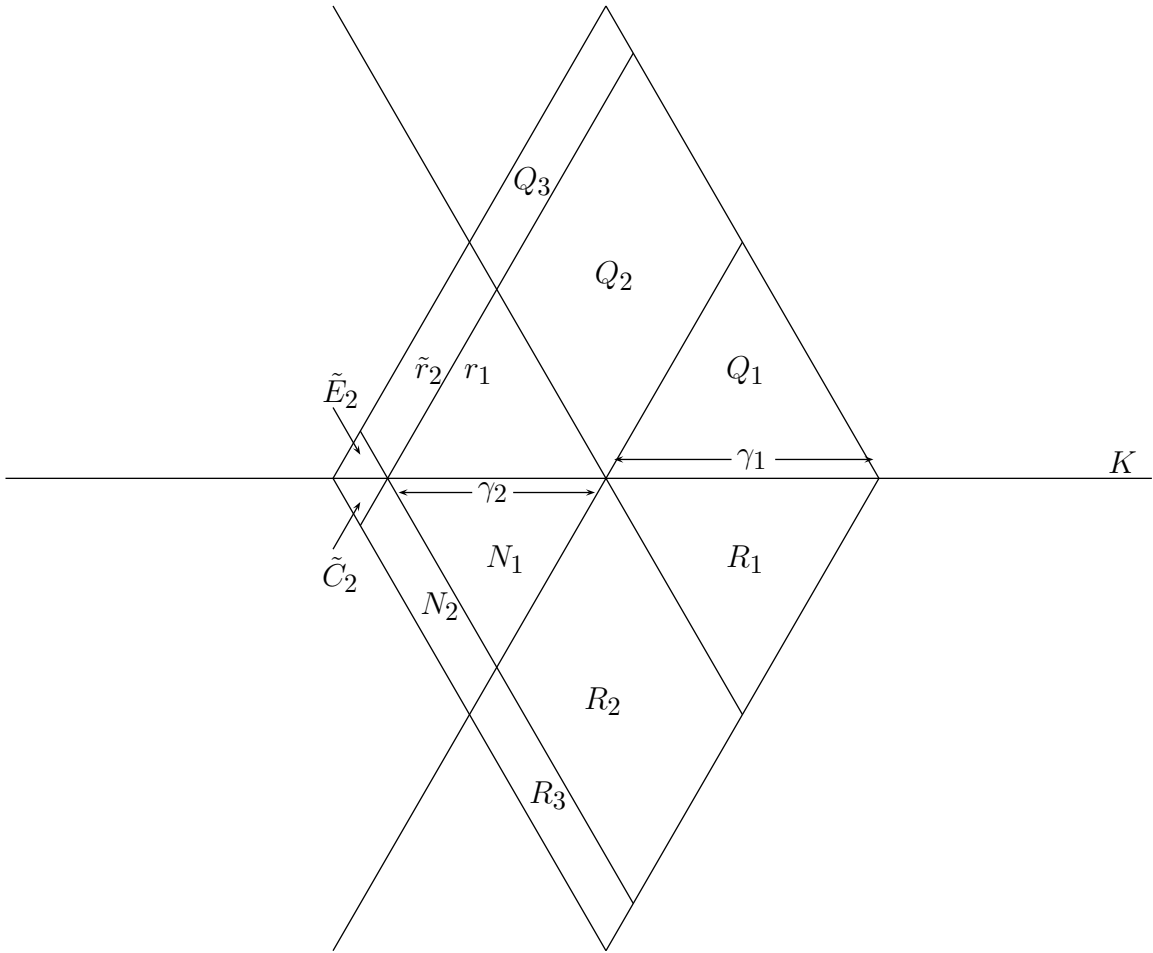}}
\caption{Existence of all values of $\gamma_1$ and $\gamma_2$, case $1$}
\label{InvariantsCase1}
\end{figure}

\noindent {\em Case 1:} $b = 1$.  Let $A_1$ be an apartment in $\mathcal{B}_1$,
and let $K$ be a wall in $A_1$ (see Figure~\ref{InvariantsCase1}).
Let $A_1^{\frac{1}{2}}$ be a half-apartment coming out of $K$ that lies in
$\mathcal{B}_2$. If $\gamma_1$ is as labeled in Figure~\ref{InvariantsCase1},
choose region $R_1$ in $A_1^{\frac{1}{2}}$.  Choose region $R_2$ from
$\mathcal{B}_6 \setminus \mathcal{B}_2$.  Let $Q_1 = \sigma (R_1)$ and
$Q_2 = \sigma (R_2 )$. This determines the regions $N_1$ and
$\tilde{r}_1$ by repeated application of Lemma~\ref{UsefulLemmaStatement},
and therefore determines $r_1 := \sigma (N_1)$. So $\tilde{r}_1$ is attached to
$N_1$.  We must check that the bottom boundaries of 
$r_1$ and $\tilde{r}_1$ have no edges
in common to ensure that $\gamma_1$ is no bigger than intended.  
But such commonality
would imply that for some edge $e$ along the upper boundary of $N_1$, 
we have $\sigma (e) = e$, i.e., $e \in \mathcal{B}_1$.  But let $G$ be
a minimal gallery through $Q_2$ from $e$ to some chamber $D$ in 
$Q_1 \subseteq \mathcal{B}_2$.  The galleries $G$ and $\sigma^2 (G)$
are both galleries of the same type from $e$ to $D$, so therefore
$G = \sigma^2 (G)$.  But $G$ passes through $Q_2 \subseteq 
\mathcal{B}_6 \setminus \mathcal{B}_2$, so $G \neq \sigma^2 (G)$, 
and we have a contradiction.

We now choose $R_3$ arbitrarily.  This determines $Q_3 = \sigma (R_3)$,
and $\tilde{r}_2$, $N_2$, $\tilde{C}_2$, $\tilde{E}_2$ by 
Lemma~\ref{UsefulLemmaStatement}.  The region $\tilde{r}_2$ is attached on 
its right hand side to $\tilde{r}_1$.  Choose $C_2 \neq \tilde{C}_2$ to ensure
that $\gamma_2$ is no bigger than intended.  Fill out the rest of $P$ 
arbitrarily.

\begin{figure}
\centerline{\input{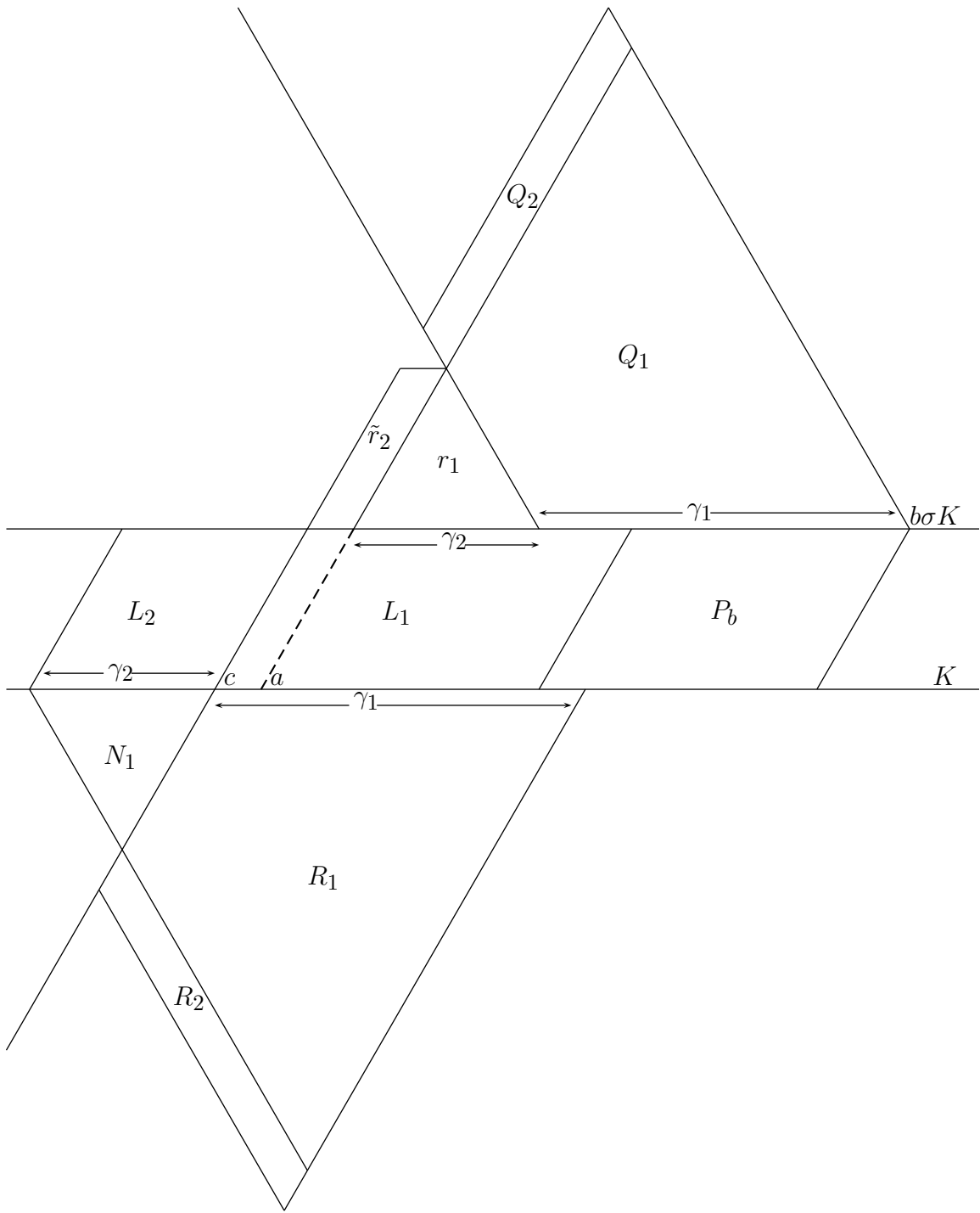}}
\caption{Existence of all values of $\gamma_1$ and $\gamma_2$, 
case 2 subcase 1}
\label{InvariantsCase2Subcase1}
\end{figure}

\noindent {\em Case 2:} $b \neq 1$.
We break this case into subcases:

{\em Subcase 1:} $\gamma_1 + T(P_b) - 1 \geq \gamma_1 + \gamma_2$, and if
equality holds then $\gamma_2 \neq 0$.  This situation is pictured in
Figure~\ref{InvariantsCase2Subcase1}.  The inequality implies that $c$ is
at least as far left as $a$ in that figure, and if $a = c$ then 
$\gamma_2 \neq 0$.  Let $K$ be a wall in $A_M$,
and let $A_1^{\frac{1}{2}}$ be a half apartment coming out of $K$.
Then $b \sigma (A_1^{\frac{1}{2}})$ is a half apartment coming out of
$b \sigma (K) = b K$.  Choose $Q_1$ and $\tilde{r}_1$ from
$b \sigma (A_1^{\frac{1}{2}})$.  Let $R_1 = \sigma^{-1} b^{-1} Q_1$.
Then $R_1$ is in $A_1^{\frac{1}{2}}$, so $L_1 \subseteq \overline{P}$ 
is in $A_M$, and $Q_1$
and $\tilde{r}_1$ are attached to $\overline{P}$ along $b \sigma (K)$.
Choose $r_1$ the same shape as $\tilde{r}_1$, and attached to $Q_1$ along
the line between $Q_1$ and $\tilde{r}_1$, but sharing no chambers in common
with $\tilde{r}_1$.   This ensures that none of the chambers in $r_1$
is attached to $L_1$, and therefore
that $\gamma_1$ is no bigger than intended.  
Let $N_1 = \sigma^{-1} b^{-1} (r_1)$. We can fill in $L_2$ by repeated
application of Lemma~\ref{UsefulLemmaStatement}.  Note that $L_2$ may be empty
if $\gamma_2 = 0$. In this case, $r_1$, $\tilde{r}_1$ and $\tilde{r}_2$
would also be empty.  If $\gamma_2 > 0$,
choose $\tilde{r}_2$ attached to $L_1$, or attached to $L_2$ if 
$c = a$ (and in this case we have $L_2 \neq \emptyset$, 
since $\gamma_2 \neq 0$). This is done by repeated application 
of Lemma~\ref{UsefulLemmaStatement}.
Choose $Q_2$ such that there is no chamber attaching it to $\tilde{r}_2$
(we know this can be done by the uniqueness statement in 
Lemma~\ref{UsefulLemmaStatement}).  
Now let $R_2 = \sigma^{-1} b^{-1} Q_2$.  This 
ensures that $\gamma_2$ is no bigger than desired.  We now choose the rest
of $P$ arbitrarily.

{\em Subcase 2:} $\gamma_1 + T(P_b) - 1 = \gamma_1 + \gamma_2$ and 
$\gamma_2 = 0$.  In this case, $T(P_b) = 1$.  The situation is pictured
in Figure~\ref{InvariantsCase2Subcase2}.  Let $A_1^{\frac{1}{2}}$
be a half apartment coming out of $K$, as before.  Choose $R_1$ in 
$A_1^{\frac{1}{2}}$, and let $Q_1 = b \sigma(R_1)$. We see that
$L_1 \subseteq A_M$.  Choose $R_2$ attached to
$R_1$.  This determines $Q_2= b \sigma(R_2)$, which in turn determines
$\tilde{D}_2$, $\tilde{L}_2$, $\tilde{E}_2$, and $\tilde{C}_2$ by 
repeated application of 
Lemma~\ref{UsefulLemmaStatement}.  For $C_2$, choose 
any chamber attached to $R_2$ other than $\tilde{C}_2$.  
This ensures that $\gamma_1$ is no bigger than intended, and that 
$\gamma_2 = 0$.
Choose the rest of $P$ arbitrarily.
\begin{figure}
\centerline{\input{fig61.tex}}
\caption{Existence of all values of $\gamma_1$ and $\gamma_2$, 
case 2 subcase 2}
\label{InvariantsCase2Subcase2}
\end{figure}

\begin{figure}
\centerline{\input{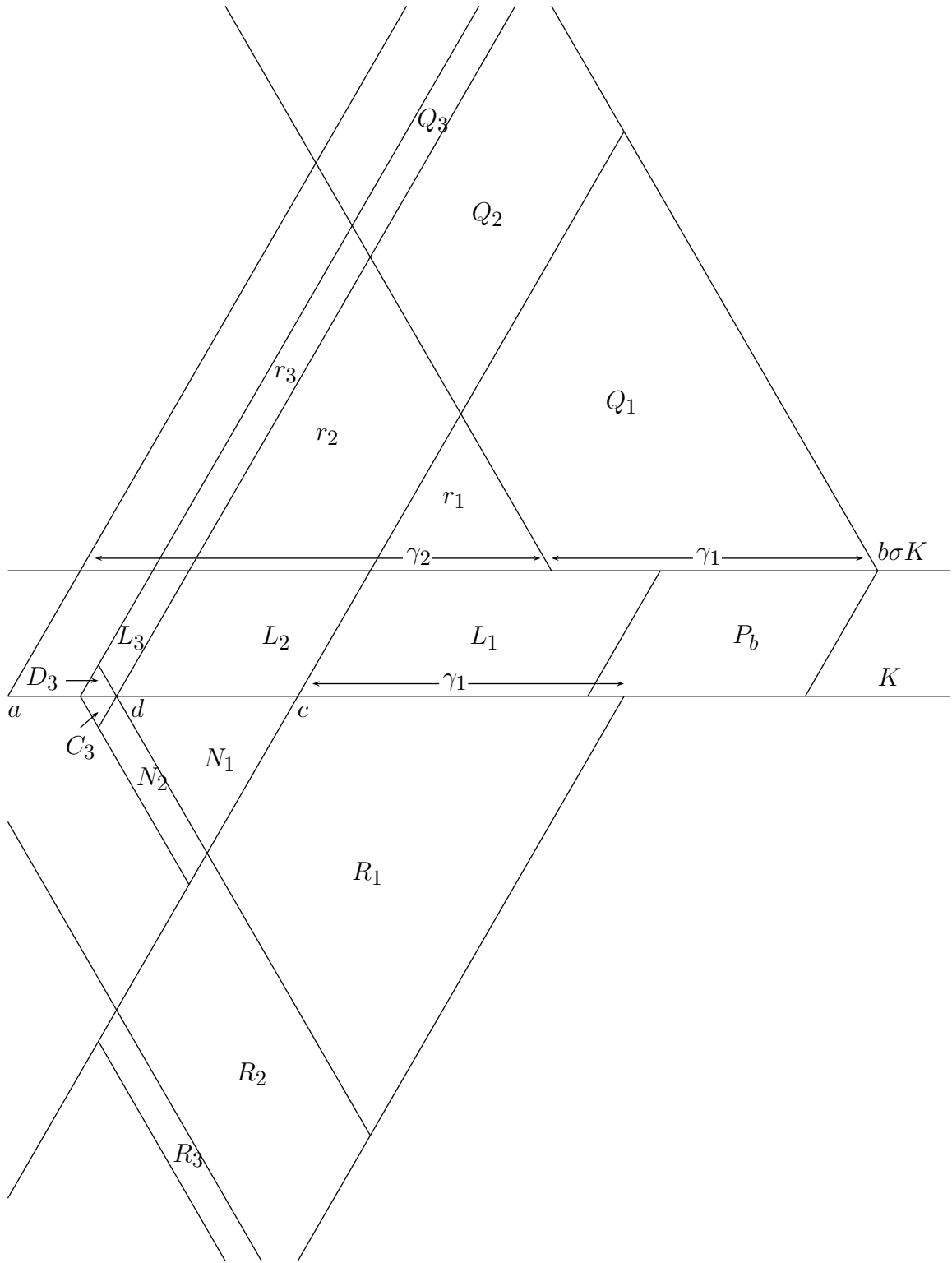}}
\caption{Existence of all values of $\gamma_1$ and $\gamma_2$, 
case 2 subcase 4}
\label{InvariantsCase2Subcase4}
\end{figure}

{\em Subcase 3:} $\gamma_1 + T(P_b) - 1 < \gamma_1 + \gamma_2$, 
and point $a$ in Figure~\ref{InvariantsCase2Subcase3} is at least
as far right as point $d$ (the alternative is 
Figure~\ref{InvariantsCase2Subcase4}, which we will consider in subcase 4).
Point $a$ is a distance of $\gamma_1 + \gamma_2 - T(P_b) + 1$
edge lengths from the right boundary of $P$, and point $d$ is a distance
of $\gamma_1 + T(P_b) - 1$ edge lengths from the same boundary.
So our condition on $a$ and $d$ is the same as saying 
$\gamma_1 + \gamma_2 - T(P_b) + 1 \leq \gamma_1 + T(P_b) - 1$,
i.e., $\gamma_2 + 2 \leq 2T(P_b)$.
Choose $Q_1$ and $\tilde{r}_1$ from $b \sigma (A_1^{\frac{1}{2}})$.  Let
$R_1 = \sigma^{-1} b^{-1} Q_1$.  Then $R_1 \subseteq A_1^{\frac{1}{2}}$,
so $L_1$ is in $A_M$, and $Q_1$ and $\tilde{r}_1$ are attached to 
$\overline{P}$ along $b \sigma (K)$.  Now consider the apartment made up
of $A_1^{\frac{1}{2}}$, $b \sigma (A_1^{\frac{1}{2}})$, and the strip
in $A_M$ between $K$ and $b \sigma (K)$.  Extend the line $l$ in
Figure~\ref{InvariantsCase2Subcase3} to a wall in this apartment, and take 
$A_2^{\frac{1}{2}}$ to be a half apartment coming out of this wall.
Choose $Q_2$, $\tilde{r}_2$, $L_2$ and $N_1$ from $A_2^{\frac{1}{2}}$.
This determines $R_2$.  Choose $Q_3$ arbitrarily.  This determines
$\tilde{r}_3$, $\tilde{L}_3$ and $\tilde{N}_3$, none of 
which will be part of the final 
construction.  Choose $N_2$ using Lemma~\ref{UsefulLemmaStatement}.  
Let $\tilde{C}$ and 
$\tilde{D}$ be the unique chambers attached to $\tilde{N}_3$ and 
$N_2$, given
by Lemma~\ref{UsefulLemmaStatement}.  
Choose $C \neq \tilde{C}$. This assures that $\gamma_2$ 
is no bigger than desired.
Construct the rest of $P$ arbitrarily.
\begin{figure}
\centerline{\input{fig62.tex}}
\caption{Existence of all values of $\gamma_1$ and $\gamma_2$, 
case 2 subcase 3}
\label{InvariantsCase2Subcase3}
\end{figure}

{\em Subcase 4:} $\gamma_1 + T(P_b) - 1 < \gamma_1 + \gamma_2$,
and point $a$ in Figure~\ref{InvariantsCase2Subcase4} is strictly
to the left of point $d$ (so $\gamma_2 + 2 > 2T(P_b)$ and the condition
$\gamma_1 + T(P_b) - 1 < \gamma_1 + \gamma_2$ become redundant).  
Construct $R_1$, $L_1$, $Q_1$ and $\tilde{r}_1$
exactly as in subcases $3$ and $1$.  Construct $N_1$, $L_2$, $\tilde{r}_2$,
and $Q_2$ exactly as in subcase 3.  This determines $R_2$ and $r_1$.
Choose $Q_3$ arbitrarily. This determines $R_3$.  Use 
Lemma~\ref{UsefulLemmaStatement}
repeatedly to fill in $\tilde{r}_3$, $L_3$ and $N_2$.  Choose the
unique possible $C_3$ and $D_3$ to connect $N_2$ and $L_3$.  Repeat
this process to point $a$, then, at the next step choose for $C_i$ any chamber
other than the unique one that connects to $L_i$.  
This ensures $\gamma_2$ is no bigger than desired.
Choose the rest of $P$ arbitrarily.

\begin{figure}
\centerline{\input{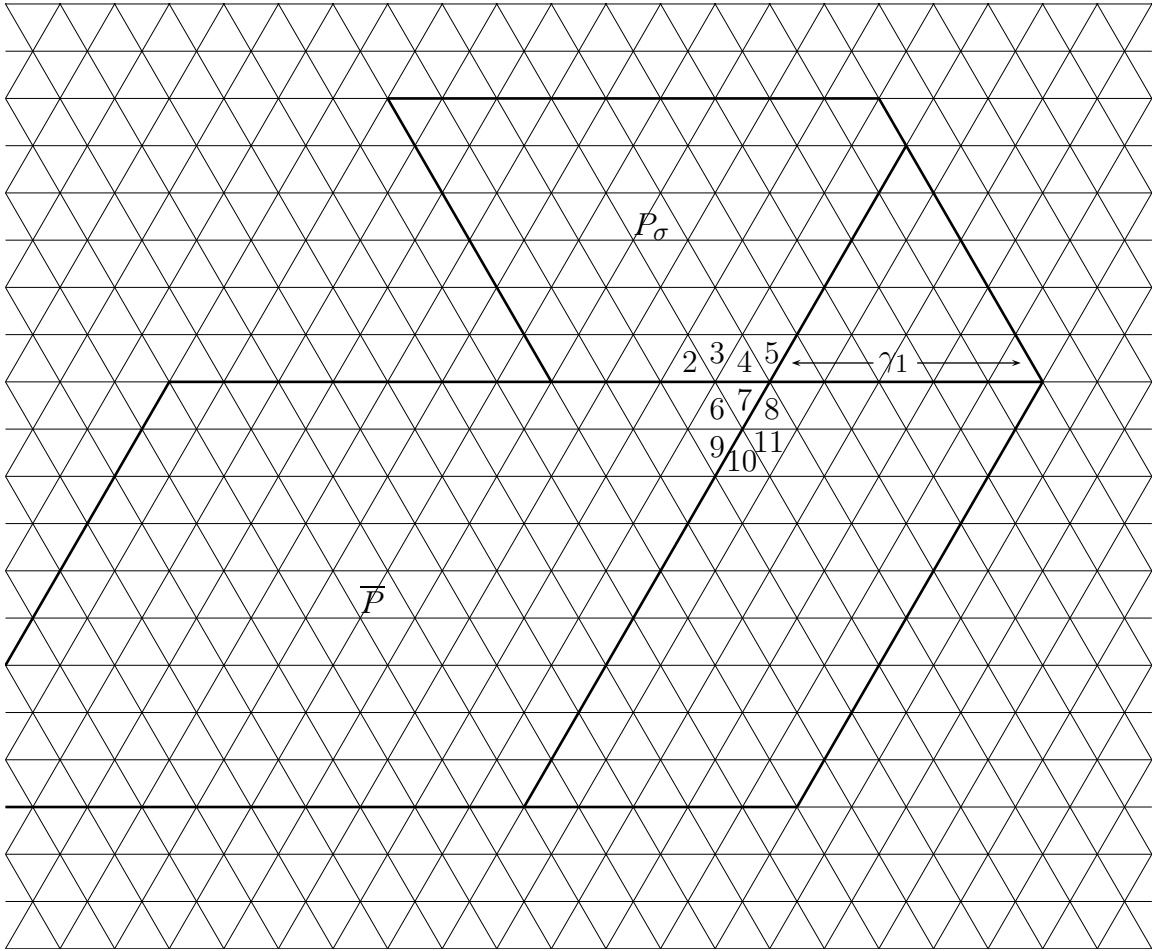}}
\caption{Definition of negative values of $\gamma_2$}
\label{Broadergamma2}
\end{figure}

Up to this point, we have defined $\overline{P}$, $P_{\sigma}$, and
the invariants $\gamma_1$ and $\gamma_2$ for $b = 1$ or $b$ with
$\alpha + 2\beta \neq 0$, and $(t,e) \in I_1$.  
We have also shown that if $1 \leq \gamma_1 \leq \gamma_1^{\textrm{max}}$,
$0 \leq \gamma_2 \leq \gamma_2^{\textrm{max}}(\gamma_1)$, then there
is some $x \in SL_3(L)$ such that $xC_M$ gives rise to
$(t,e)$, $\gamma_1$ and $\gamma_2$.
We now broaden the definition of $\gamma_2$ to include negative values.
We will then have to re-address the issue of whether all possible
$\gamma_2$ actually occur.  
Consider Figure~\ref{Broadergamma2}. If there is no chamber connecting
chamber $7$ to chamber $5$, we say $\gamma_2 \leq 0$ 
(previously we just said $\gamma_2 = 0$).
In this case we ask if there is a chamber connecting $4$ to $8$. If not,
then $\gamma_2 = 0$.  If so, then $\gamma_2 \leq -1$, and we call
the connecting chamber $\tilde{7}$.  In this case, we fill in chambers
$\tilde{6}$ and $\tilde{9}$ using chambers 
$\tilde{7},8,11,10$ and Lemma~\ref{UsefulLemmaStatement}.
We now ask if there is a chamber connecting $\tilde{6}$ to $2$. If not, then
$\gamma_2 = -1$, and if so then $\gamma_2 \leq -2$.  We proceed in this way
to determine the value of $\gamma_2$. Note that if 
$\gamma_2^{\textrm{min}}(\gamma_1)
= - \min \{ S(\overline{P}), T(P) - \gamma_1 \}$, then 
$\gamma_2 \geq \gamma_2^{\textrm{min}}$ automatically.

We now show that given any $b$ such that $b = 1$ or $\alpha + 2\beta \neq 0$,
given $(t,e) \in I_1$, and given $\gamma_1$ and $\gamma_2$
such that $1 \leq \gamma_1 \leq \gamma_1^{\textrm{max}}$ and 
$\gamma_2^{\textrm{min}}(\gamma_1) 
\leq \gamma_2 \leq \gamma_2^{\textrm{max}}(\gamma_1)$, 
there exists an $x$ such that $x$ and $b$ give rise to $(t,e)$,
$\gamma_1$ and $\gamma_2$.  Given what we have already done, we may assume
that $\gamma_2 \leq 0$. As usual we proceed in cases.

\noindent {\em Case 1:} $b = 1$. See Figure~\ref{Negativegamma2case1}.
Let $L$ be a wall in $A_M$, and let $A_1^{\frac{1}{2}}$ be a half apartment
in $\mathcal{B}_2$ coming out of $L$.  Choose $Q_1$ in $A_1^{\frac{1}{2}}$.
Choose chamber $4$ in 
$\mathcal{B}_6 \setminus \mathcal{B}_2$.  This determines $8$, and therefore
by Lemma~\ref{UsefulLemmaStatement} also determines 
$\tilde{3}$ and $\tilde{7}$
connecting $4$ and $8$.  To arrange $\gamma_2 = 0$, choose $7$ such that
$7 \neq \tilde{7}$ and $\sigma(7) \neq \tilde{3}$, and construct the
rest of $P$ arbitrarily.  To arrange $\gamma_2 \leq -1$, choose 
$7 = \sigma^{-1}(\tilde{3})$.  Since $4, 8 \in \mathcal{B}_6 \setminus
\mathcal{B}_2$, $\sigma(\tilde{7}) \neq \tilde{3}$, so $\tilde{7} \neq 7$,
so $\gamma_2 \leq 0$ still holds.  On the other hand $3 = \sigma(7) 
= \tilde{3}$, so $3$ can be connected to $8$ (via $\tilde{7}$), ensuring
$\gamma_2 \leq -1$.  Choose $\tilde{6}$ and $\tilde{9}$ to fit between
$\tilde{7}$ and $10$ using Lemma~\ref{UsefulLemmaStatement}.  
Choose $\tilde{1}$
and $\tilde{5}$ to fit between $2$ and $\tilde{6}$.  To arrange
$\gamma_2 = -1$ choose $1 \neq \tilde{1}$ and construct the rest
of $P$ arbitrarily.  To arrange $\gamma_2 \leq -2$, choose 
$1 = \tilde{1}$.  Continue in this manner.
\begin{figure}
\centerline{\input{fig65.tex}}
\caption{Existence of all negative values of $\gamma_2$, case $1$}
\label{Negativegamma2case1}
\end{figure}

\noindent {\em Case 2:} $b \neq 1$.  We break this case into subcases.

\begin{figure}
\centerline{\input{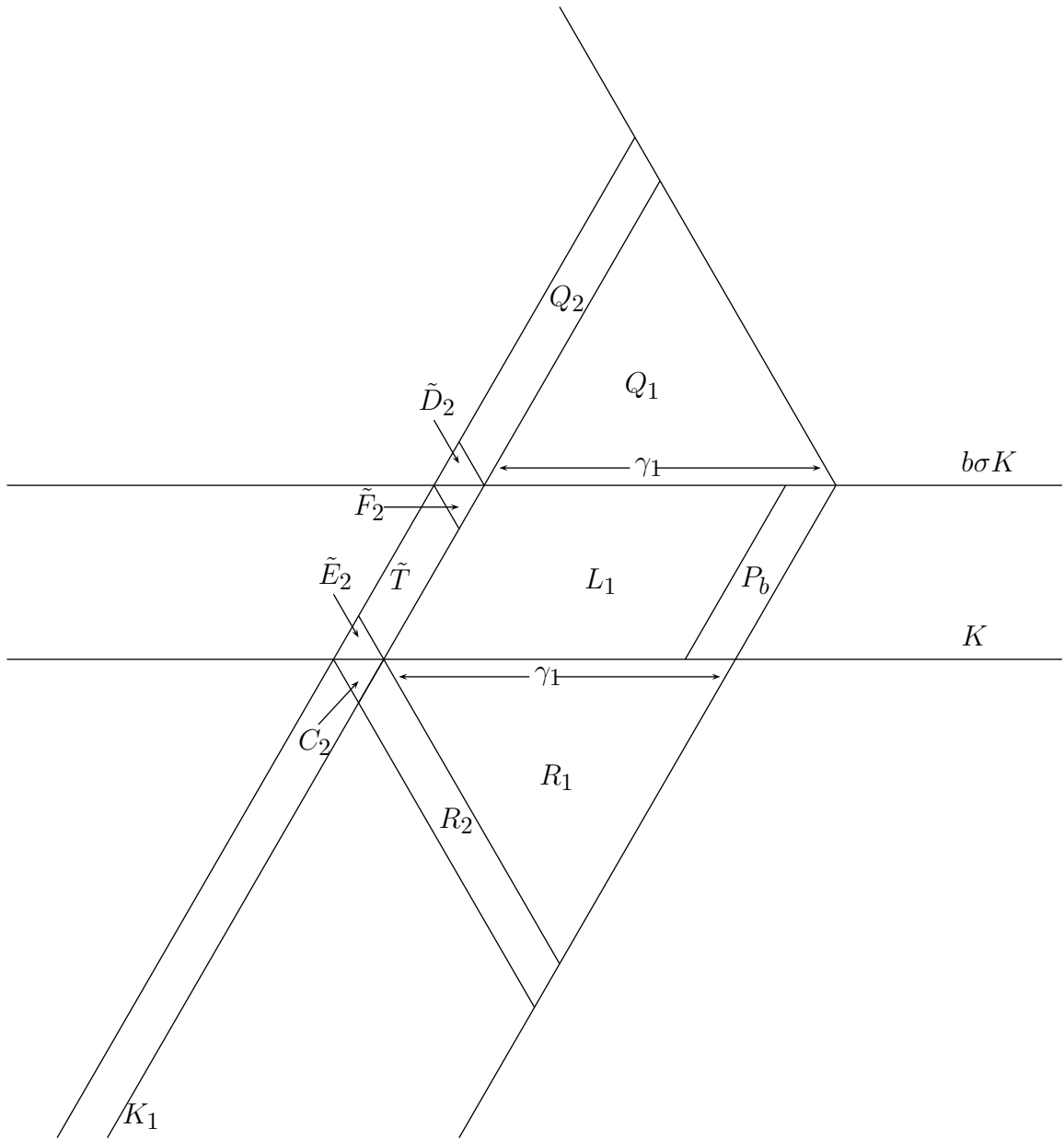}}
\caption{Existence of all negative values of $\gamma_2$, case $2$ subcase $1$}
\label{Negativegamma2case2subcase1}
\end{figure}

{\em Subcase 1:} $T(P_b) = 1$.  See Figure~\ref{Negativegamma2case2subcase1}.
We let $K$ and $K_1$ be walls in $A_M$, $K$ a horizontal wall, and $K_1$ 
a wall such that $bK_1 = K_1$.  Let $\tilde{T}$ be the first row of chambers
in a half apartment $A_1^{\frac{1}{2}}$ coming out of $K_1$ to the left.
Choose $A_1^{\frac{1}{2}}$ such that $b\sigma(\tilde{T}) \neq \tilde{T}$.
This implies that $b\sigma(\tilde{T})$ has no chambers in common with
$\tilde{T}$.  Let $A$ be the apartment made 
up of $A_M$ on the right of $K_1$ and 
$A_1^{\frac{1}{2}}$ on the left.  Then half of $K$ is in $A$, and we can
extend this half to $\tilde{K}$, a wall in $A$.  Let $A_2^{\frac{1}{2}}$
be a half apartment coming out of $\tilde{K}$, and choose $R_1$,
$R_2$ and $\tilde{C}_2$ from $A_2^{\frac{1}{2}}$.  We can now throw 
away all of $A_1^{\frac{1}{2}}$ except $\tilde{T}$.  We get 
$Q_1 = b\sigma(R_1)$, $Q_2 = b\sigma(R_2)$, $P_b$, and $L_1 \subseteq A_M$.
Note that $b\sigma(\tilde{T})$ is another row of chambers coming out of 
$L$, and note that we have $\tilde{D}_2 = b\sigma(\tilde{C}_2)$ attaching
$Q_2$ to $b\sigma(\tilde{T})$, as drawn in the figure.  We have 
$\tilde{F}_2 \subseteq b\sigma(\tilde{T})$ which connects $\tilde{D}_2$
to $L_1$, and by Lemma~\ref{UsefulLemmaStatement}, we have $\tilde{E}_2$
and $\tilde{\tilde{C}}_2$ attaching $b\sigma(\tilde{T})$ to $R_2$ as drawn.
Note that $\tilde{\tilde{C}}_2 \neq \tilde{C}_2$.  If these chambers were 
equal, we could create a contradiction to the uniqueness part of 
Lemma~\ref{UsefulLemmaStatement} by producing two galleries of length
two from the edge between $\tilde{\tilde{C}}_2 = \tilde{C}_2$
and $\tilde{E}_2$, down through $\tilde{T}$ and $b\sigma(\tilde{T})$
respectively, to $L_1$.  To make sure $\gamma_2 \leq 0$, choose
$C_2 \neq \tilde{\tilde{C}}_2$.  To make sure $\gamma_2 \leq -1$,
choose $C_2 = \tilde{C}_2$.  For $\gamma_2 = 0$, choose any other
$C_2$. If we chose $\gamma_2 = 0$, fill out the rest of $P$ arbitrarily.
If we chose $\gamma_2 \leq -1$, then what we have done becomes the base
case of an induction.

\begin{figure}
\centerline{\input{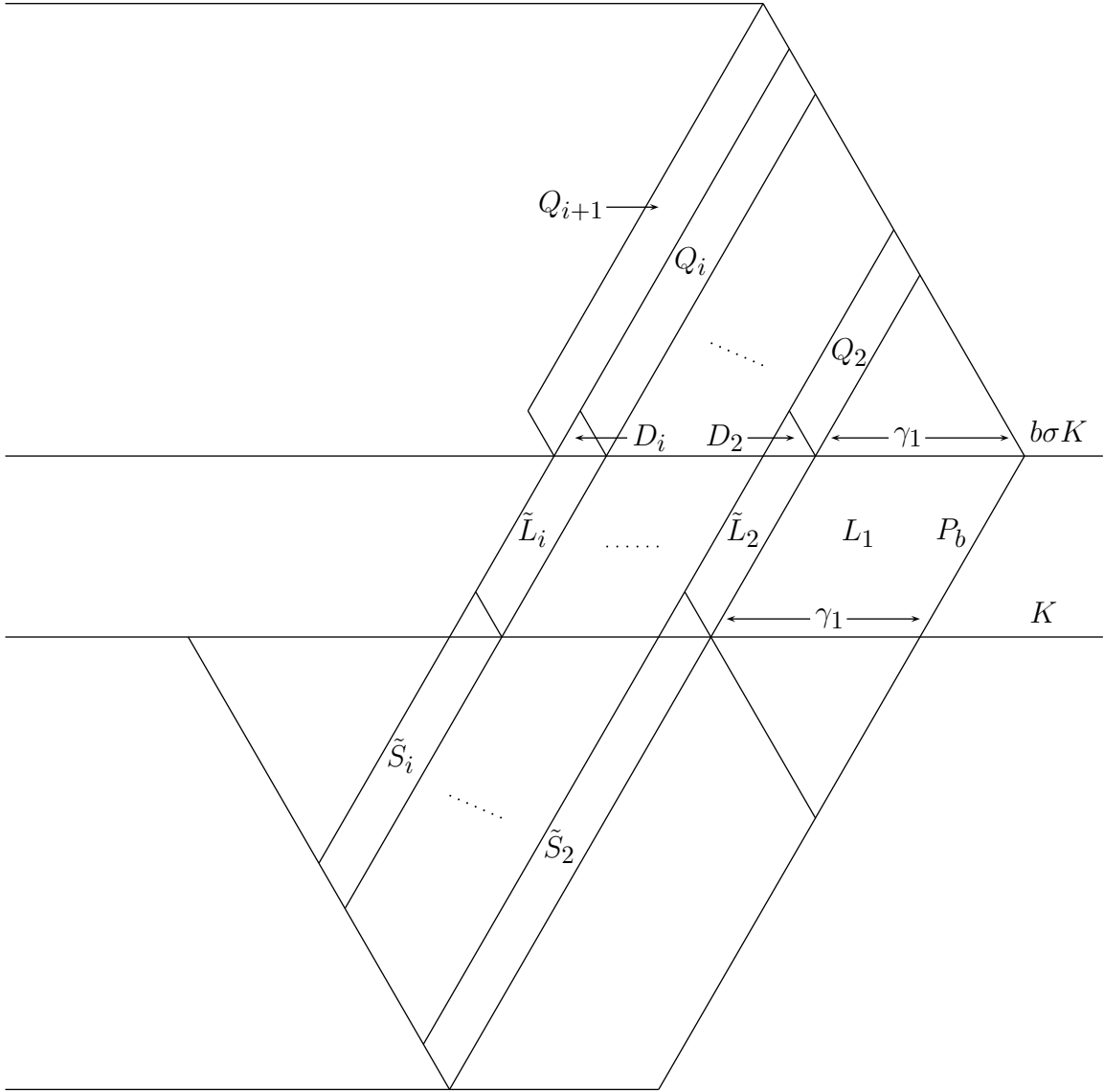}}
\caption{Existence of all negative values of $\gamma_2$, case $2$ subcase $1$
diagram $2$}
\label{Negativegamma2Case2Subcase1}
\end{figure}

Let $D_2 = b\sigma(C_2)$, and see 
Figure~\ref{Negativegamma2Case2Subcase1}.  The region made up of 
$\tilde{L}_2$, $\tilde{S}_2$ is taken from $b\sigma(\tilde{T})$.
We assume by induction
that the $Q_j$, $D_j$, $\tilde{L}_j$ and $\tilde{S}_j$ for $2 \leq j \leq i$
are all contained in a single half apartment, and that $\gamma_2 \leq -i +1$.
We choose $Q_{i+1}$ arbitrarily.  If we want $\gamma_2 = -i + 1$, we choose
$D_{i+1}$ such that there is no chamber attaching it to $\tilde{L}_i$. 
If we want
$\gamma_2 \leq -i$, we choose a 
chamber $D_{i+1}$ that can be attached to $\tilde{L}_i$.
This enables us to fill in an $\tilde{L}_{i+1}$ and an $\tilde{S}_{i+1}$ 
uniquely
using Lemma~\ref{UsefulLemmaStatement}.  We can continue thus until we either
reach the desired value for $\gamma_2$, or until we are stopped by 
$\gamma_2^{\textrm{min}} (\gamma_1)$.

{\em Subcase 2:} $T(P_b) > 1$.  See Figure~\ref{Negativegamma2Case2Subcase2}.
Construct $R_1$ and $Q_1$ in the standard way.
Choose $Q_2$ such that there is no chamber connecting it to $L_1$ (to ensure
$\gamma_2 \leq 0$).  If we want $\gamma_2 = 0$, we choose $D_2$ such that
there is no chamber connecting it to the line $l$ (which is in $L_1$ and $P$).
If we want $\gamma_2 \leq -1$, we choose $D_2$ such that there is a chamber
$\tilde{E}_2$ connecting it to the line $l$.  Note that $\tilde{E}_2 
\nsubseteq L_1$.
Construct $\tilde{S}_2$ using Lemma~\ref{UsefulLemmaStatement}. 
The region $\tilde{S}_2$ does
not have any chambers is common with $L_1$ or with $P$.  We then choose
$Q_3$ arbitrarily, and if we want $\gamma_2 = -1$, we choose $D_3$ such that
there is no chamber connecting it to $S_2$. If we want $\gamma_2 \leq -2$,
we choose $D_3$ such that there is such a chamber.  Proceed in this way.
\begin{figure}
\centerline{\input{fig68.tex}}
\caption{Existence of all negative values of $\gamma_2$, case $2$ subcase $2$}
\label{Negativegamma2Case2Subcase2}
\end{figure}

We have now shown that given a $\sigma$-conjugacy class $b$ such that
either $b = 1$ or $\alpha + 2\beta \neq 0$, and given 
a type edge pair $(t,e) \in I_1$, one can find 
$x \in SL_3 (L)$ such that $xC_M$ has SMG of type $t$ and departure edge
$e$, and such that the invariants $\gamma_1$, $\gamma_2$ associated
with the composite gallery $\Gamma_x$ are as desired (within the 
limitations $1 \leq \gamma_1 \leq \gamma_1^{\textrm{max}}$ and
$\gamma_2^{\textrm{min}}(\gamma_1) 
\leq \gamma_2 \leq \gamma_2^{\textrm{max}}(\gamma_1)$).
This is useful because we will now prove that $\rho (x^{-1}b \sigma (x) C_M )$
is determined by $b$, $t$, $e$, $\gamma_1$, and $\gamma_2$.  Again, we
consider cases.

\begin{figure}
\centerline{\input{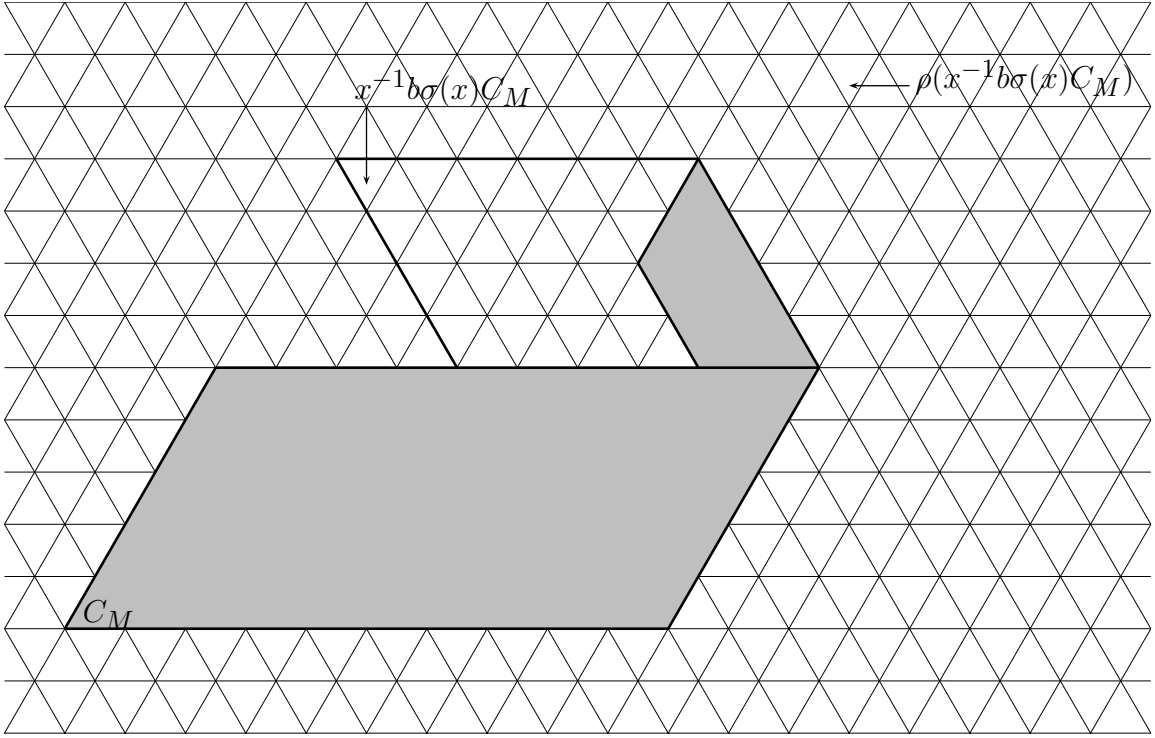}}
\caption{Example of how $\gamma_1$ and $\gamma_2$ determine 
$\rho (x^{-1} b \sigma (x) C_M )$}
\label{BackTransportExample}
\end{figure}

\noindent {\em Case 1:} $\gamma_2 > 0$.  
See Figure~\ref{BackTransportgamma2gt0}, which represents the structure
$\overline{P} \cup P_{\sigma}$ transported back to $C_M$ (in other words,
$x^{-1}(\overline{P} \cup P_{\sigma}))$. We can find $g \in I$ that sends 
the shaded part of this picture to $A_M$.  This is because one can find 
an apartment containing the shaded part of the figure but not containing
the unshaded part.  Therefore
$g$ does not send the unshaded part of the diagram to $A_M$. So we assume
from this point forward that the shaded part of 
Figure~\ref{BackTransportgamma2gt0} is in $A_M$ and the unshaded part is not.
It is easy to see now that $\rho$ has the effect of reflecting the chambers
in $R_1$ across $l_1$, reflecting the chambers in $R_2$ across
$l_2$, reflecting the chambers in $R_3$ first across $l_1$ and then 
across $l_3$,
and reflecting the chambers in $R_4$ across $l_2$.  So one can determine the
value of $\rho (x^{-1} b \sigma (x) C_M )$.  A specific example is given 
in Figure~\ref{BackTransportExample}. Note that the existence of a corner
at point $x$ in the figure is necessary for the considerations that arrive
at these results. 
\begin{figure}
\centerline{\input{fig69.tex}}
\caption{$x^{-1}(\overline{P} \cup P_{\sigma})$ for $\gamma_2 > 0$}
\label{BackTransportgamma2gt0}
\end{figure}

\begin{figure}
\centerline{\input{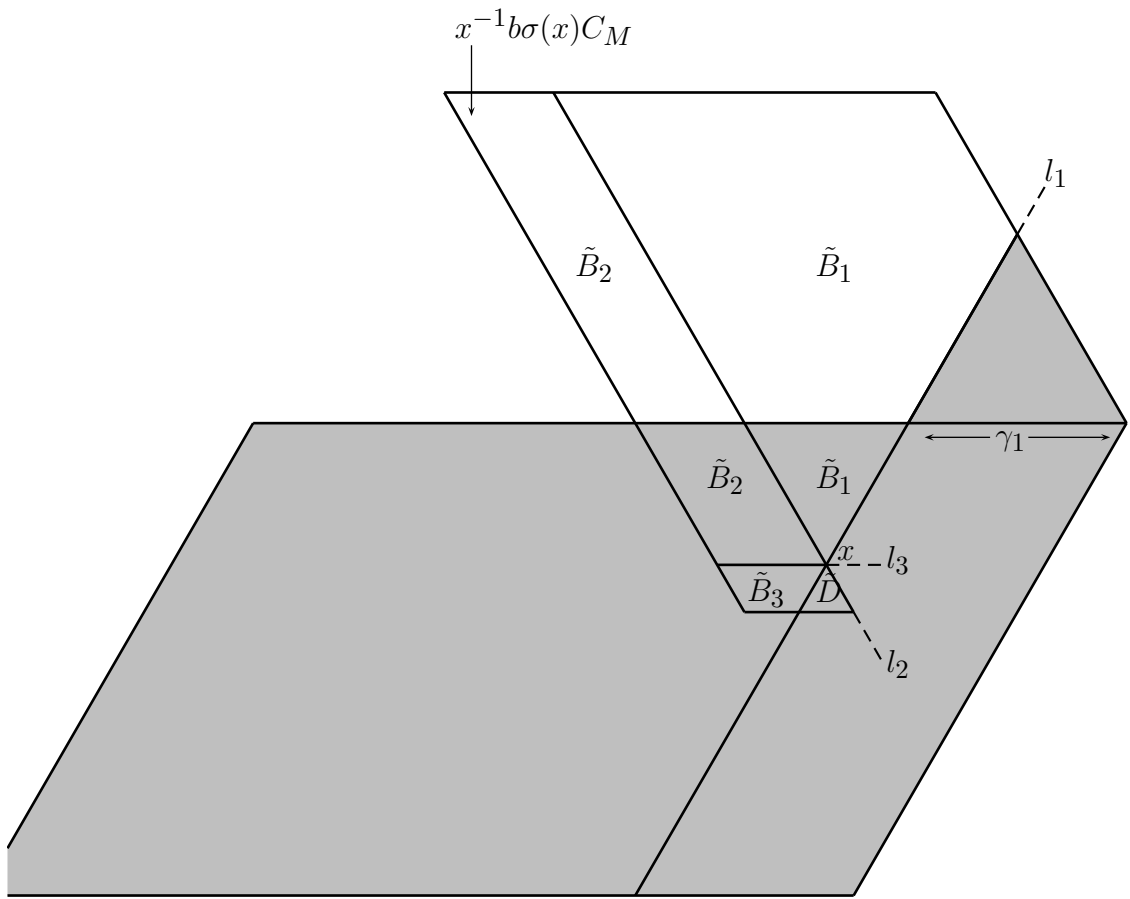}}
\caption{$x^{-1}(\overline{P} \cup P_{\sigma})$ for $\gamma_2 \leq 0$, folding
explained further}
\label{BackTransportgamma2ltet0_2}
\end{figure}

\noindent {\em Case 2:} $\gamma_2 \leq 0$. 
See Figure~\ref{BackTransportgamma2ltet0}.  Once again,
we can arrange for the shaded part of the figure to be in $A_M$, and for the
unshaded part to not be in $A_M$. According to the value of $\gamma_2$,
we can connect $P_{\sigma}$ to $l$ using region $\tilde{R}_1$ in 
Figure~\ref{BackTransportgamma2ltet0_1}.  Here, $l$ is in 
$\overline{P}$, and $\tilde{R}_1$ shares no chambers in common
with $\overline{P}$.  We can now fill in $\tilde{R}_2$ uniquely
using Lemma~\ref{UsefulLemmaStatement} so that it is connected above to 
$P_{\sigma}$ and on the right to $\tilde{R}_1$.  
There now exist chambers $\tilde{C}$ and $\tilde{D}$
such that neither is in $\overline{P}$, but such that $\tilde{C}$ is
connected to $\tilde{R}_2$, and $\tilde{D}$ has its right 
edge in $\overline{P}$. Also, $\tilde{C}$ and $\tilde{D}$ are connected to
each other.  We fill in $\tilde{R}_3$ using 
Lemma~\ref{UsefulLemmaStatement}.
With respect to folding, the corner labeled with an $x$ is similar in 
function to the corner labeled with an $x$ in 
Figure~\ref{BackTransportgamma2gt0}.  As such, referring to 
Figure~\ref{BackTransportgamma2ltet0_2}, chamber $\tilde{D}$ is
reflected across $l_2$ by $\rho$, the chambers in region $\tilde{B}_1$ 
(some of which are in $P_{\sigma}$ and some of which are not in
$\overline{P} \cup P_{\sigma}$) are
reflected across $l_1$, the chambers in $\tilde{B}_2$ (some of which are
in $P_{\sigma}$ and some of which are not in $\overline{P} \cup P_{\sigma}$) 
are reflected first across
$l_1$, then across $l_3$, and the chambers in $\tilde{B}_3$ 
are reflected across 
$l_2$.  So one can determine the value of 
$\rho ( x^{-1} b \sigma (x) C_M )$ for any specific example.
\begin{figure}
\centerline{\input{fig71.tex}}
\caption{$x^{-1}(\overline{P} \cup P_{\sigma})$ for $\gamma_2 \leq 0$}
\label{BackTransportgamma2ltet0}
\end{figure}
\begin{figure}
\centerline{\input{fig72.tex}}
\caption{$x^{-1}(\overline{P} \cup P_{\sigma})$ for $\gamma_2 \leq 0$, folding
explained}
\label{BackTransportgamma2ltet0_1}
\end{figure}

In Section~\ref{Superset}, we computed chambers to include in the superset
$S_1$ by computing all possible foldings of the composite gallery
$\Gamma_x$ for $x$ giving rise to any fixed type-edge pair $(t,e)$. 
The set $S_1$ consisted of the union of these results for all $(t,e)$
in $I_1 \cup I_2$.  We are now in a position to show that for all
$(t,e) \in I_1$, the collection of possible foldings of $\Gamma_x$
gives rise only to chambers in $S$.  This, combined with similar
results for $I_2$ proves that $S = S_1$ for $b = 1$ or $b$ with 
$\alpha + 2\beta \neq 0$.

Let $S^f_{(t,e)}$ be the collection of chambers obtained by computing
the possible foldings of $\Gamma_x$ for $x$ giving rise to the type-edge
pair $(t,e)$.  So 
$$S_1 = \left[ \cup_{(t,e) \in I_1} S^f_{(t,e)} \right] \cup
\left[ \cup_{(t,e) \in I_2} S^f_{(t,e)} \right].$$
Let $S^{\gamma}_{(t,e)}$ be the collection of the 
$\rho (x^{-1} b \sigma (x) C_M )$ (which can now be computed for 
$(t,e) \in I_1$ and either $b = 1$ or $\alpha + 2\beta \neq 0$) coming from
$x$ giving rise to $(t,e)$ and any invariants $\gamma_1$ and $\gamma_2$ such
that $1 \leq \gamma_1 \leq \gamma_1^{\textrm{max}}(t,e)$ and
$\gamma_2^{\textrm{min}}(t,e,\gamma_1) 
\leq \gamma_2 \leq \gamma_2^{\textrm{max}}(t,e,\gamma_1)$.
We will show that $S^{\gamma}_{(t,e)} = S^f_{(t,e)}$ (again, only for
$(t,e) \in I_1$ and either $b = 1$ or $\alpha + 2\beta \neq 0$).  To do this,
we begin by dissecting $S^f_{(t,e)}$. 
\begin{figure}
\centerline{\input{fig74.tex}}
\caption{Choice Points}
\label{ChoiceEdges}
\end{figure}

When computing the possible values of $\rho (x^{-1} b \sigma (x) C_M )$,
one can work through the chambers of $x^{-1}\Gamma_x$ 
starting with $C_M$, making
a choice of folding direction at each choice point.  At 
choice point $e$ between
chambers $i$ and $i+1$, we define the {\em status quo choice} to be the one
in which $i$ and $i+1$ go to different chambers. The {\em change choice}
is the choice in which $i$ and $i+1$ go to the same chamber.  Note that
if a change choice is made at an edge $e$, then choice points subsequent
to $e$ in $x^{-1}\Gamma_x$ may become non-choice points, and vice versa. The 
possible values of $\rho (x^{-1} b \sigma (x) C_M )$ can be 
enumerated by exploring the binary tree described by the choice point
structure of $x^{-1}\Gamma_x$.

It is easy to see that if one makes the status quo choice at all edges 
before some edge $u_i$ in Figure~\ref{ChoiceEdges}, and then one makes
a change choice at $u_i$ itself, then all subsequent edges become 
non-choice points and $\rho (x^{-1} b \sigma (x) C_M )$ is determined.
This value for $\rho ( x^{-1} b \sigma (x) C_M )$ is in $S^{\gamma}_{(t,e)}$
because we may choose $\gamma_1 = i$ and $\gamma_2$ maximal. The value
of $\rho (x^{-1} b \sigma (x) C_M )$ obtained by making only status quo
choices is in $S^{\gamma}_{(t,e)}$ because we may choose $\gamma_1$ and
$\gamma_2$ maximal.

It is easy to see that if one makes the first change choice at some 
$d_i$, then one may make at most one subsequent change choice, and 
such a subsequent change choice can only be made an odd number $s$
of edges after $d_i$.  The resulting value of 
$\rho (x^{-1} b \sigma (x) C_M )$ is in $S^{\gamma}_{(t,e)}$ because
we may choose $\gamma_1$ and $\gamma_2$ such that
$\gamma_1 + \max \{ \gamma_2 , 0 \} = S(P_\sigma ) + i - 1 $ 
(to ensure that the first change choice occurs at $d_i$ )
and such that $\gamma_1 - \min \{ 0 , \gamma_2 \} = \frac{s+1}{2}$ 
(to ensure that the second change choice occurs $s$ edges after $d_i$).
If no second change choice is desired, then choose $\gamma_2 = 
\gamma_2^{\textrm{min}}(\gamma_1)$.

If one makes the first change choice at some $e_i$, then only one subsequent
change choice is possible.  If such a change choice is made, it must
be done at one of the $u_j$.  The resulting 
$\rho (x^{-1} b \sigma (x) C_M )$ is in $S^{\gamma}_{(t,e)}$
because we may choose $\gamma_1$ and $\gamma_2$ such that
$\gamma_1 + \max \{ \gamma_2 , 0 \} = i$ and such that
$\gamma_1 - \min \{ \gamma_2 , 0 \} = j$.  If no second change choice 
is desired, choose $\gamma_2 = \gamma_2^{\textrm{min}}(\gamma_1)$.

In essence, the choice of $\gamma_1$ and $\gamma_2$ above corresponds
to choosing the point $x$ in Figures~\ref{BackTransportgamma2gt0} 
and~\ref{BackTransportgamma2ltet0_2} such that $l_1$ 
passes through the first change choice and $l_2$ passes through
the second.

We now address the situation where $(t,e) \in I_1$ still, but now $b \neq 1$
and $\alpha + 2\beta = 0$.  As before, let $p$ be the number of chambers
in $\Gamma_E^3$ between the edge of departure 
and the turning edge.  So $p \geq 1$ is odd.
We first consider the case where $p = 1$.  So, fixing $t,e,b$ with the 
specified properties, and choosing $y \in SL_3(L)$ such that $y$ gives
rise to $(t,e)$, and given some possible folding of a gallery of the same
type as $y^{-1}\Gamma_y$, we want to show that there is $x \in SL_3(L)$
such that $\rho(x^{-1}b\sigma(x)C_M)$ is the last chamber of this folding.

\begin{figure}
\centerline{\input{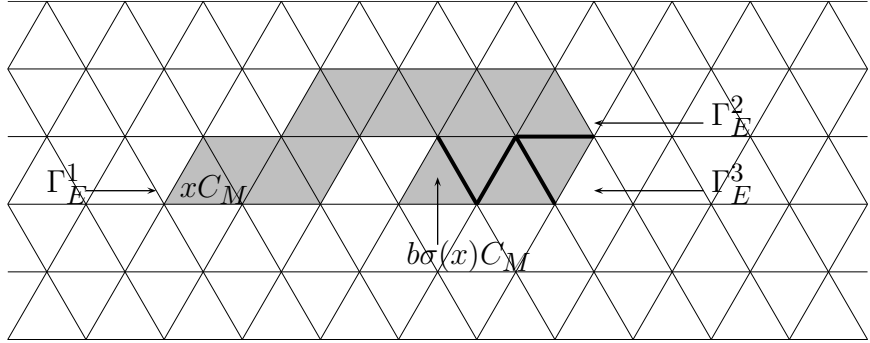}}
\caption{A composite gallery for $p = 1$, $\alpha = 2$, $\beta = -1$}
\label{alpa2betaiszerop1example}
\end{figure}

For $p = 1$, a sample composite gallery for $b$ with $\alpha = 2$,
$\beta = -1$ is pictured in Figure~\ref{alpa2betaiszerop1example}.
The choice points are marked.  For general $p$, $b$, we get a similar
structure, with every edge in $\Gamma_E^3$ a choice point, and the edge
between $\Gamma_E^2$ and $\Gamma_E^3$ another choice point.  Note that
for $p = 1$, if one makes a change choice at any choice point other than the 
horizontal one between $\Gamma_E^2$ and $\Gamma_E^3$, no subsequent choices
are available.  A change choice at the horizontal change point allows
at most two more choices.

Let $K$ be a wall in $A_M$ such that $bK = K$.  Let $T$ be the first
row of chambers in some half apartment coming out of $K$.  Choose
$T \subseteq \mathcal{B}_1$, so $\sigma$ does not affect $T$.  Then
$b\sigma(T) = T$.  Choose the first $i$ chambers of $\Gamma_E^3$
in $T$, and the rest not in $T$. So $\Gamma_E^1 
= \sigma^{-1}b^{-1}(\Gamma_E^3)$ and $\Gamma_E^2 \subseteq A_M$ are 
determined.  This composite gallery construction folds in such a way
that the first (and only) change choice occurs in $\Gamma_E^3$
between the $i$th and $(i+1)$st chambers.  To arrange no change choices at
all, choose $\Gamma_E^3 \subseteq T$.

\begin{figure}
\centerline{\input{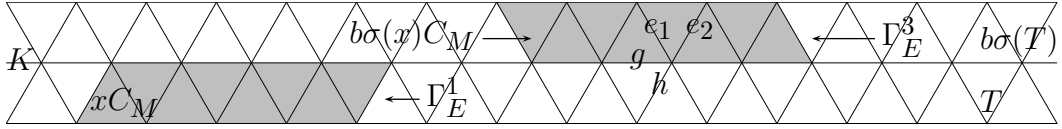}}
\caption{All foldings of composite galleries with $p=1$ occur, where $b \neq 1$
and $\alpha + 2\beta = 0$}
\label{allfoldingsoccur}
\end{figure}

To arrange that the edge between $\Gamma_E^2$ and $\Gamma_E^3$ be a change 
choice, take $T \subseteq \mathcal{B}_2 \setminus \mathcal{B}_1$.  Then 
$b\sigma(T) \cap T = K$.  If we want no subsequent change choice, then 
choose $\Gamma_E^3 \subseteq b\sigma(T)$. 
If we want the second change choice to
be between the $i$th and $(i+1)$st chambers in $\Gamma_E^3$, then
choose the first $i$ chambers of $\Gamma_E^3$ in $b\sigma(T)$, and the rest
outside of $b\sigma(T)$.  We now have a situation for which 
Figure~\ref{allfoldingsoccur} is an example
(here $b$ has $\alpha = 4$ and $\beta = -2$).  The galleries $\Gamma_E^1$
and $\Gamma_E^3$ are drawn as part of $T$ and $b\sigma(T)$, respectively,
in this diagram, even though only part of each need be.

Let $e$ be the edge in $\Gamma_E^3$ at which the second change choice
is to occur.  Then the chambers of $\Gamma_E^3$ that are part of $b\sigma(T)$
are exactly those that occur before $e$.  If $e$ is an edge of the same
angle as $e_1$ in Figure~\ref{allfoldingsoccur}, then no subsequent
change choice is possible.  If $e$ is of the same angle as $e_2$,
then $e_1$ could be a third (and final) change choice.  To arrange this,
choose $\Gamma_E^3$ in such a way that there is no chamber connecting
edge $g$ to edge $h$.  To arrange no third change choice, choose
$\Gamma_E^3$ such that there is a chamber connecting $g$ to $h$.

\begin{figure}
\centerline{\input{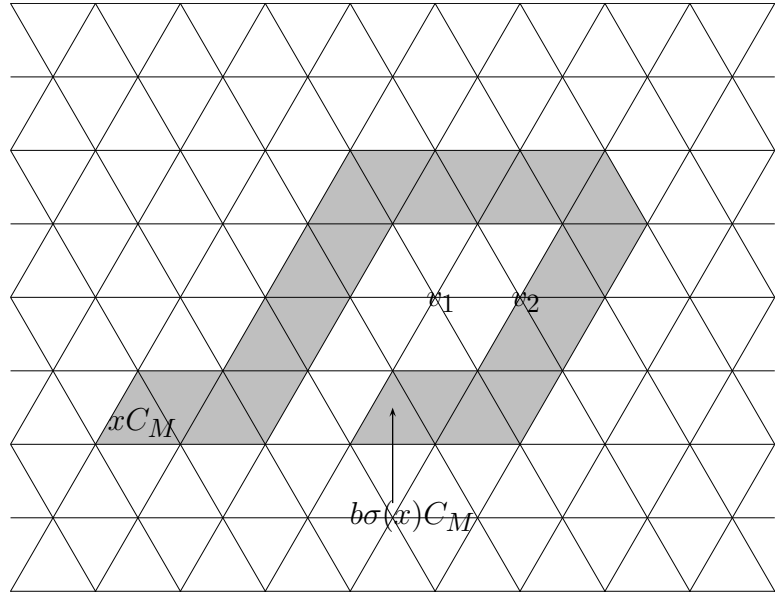}}
\caption{Definition of $\gamma_3$}
\label{gamma3}
\end{figure}

We have now shown for $p = 1$ that any possible folding can occur.
We now assume inductively that any possible folding can occur for
$p \leq n$.  Suppose that $x$ gives rise to a type-edge pair $(t,e)$ with
$p = n+1$, and a composite gallery $\Gamma_x$.  We begin by defining
an invariant $\gamma_3$ of $x$.  Consider Figure~\ref{gamma3}, which shows an
example composite gallery for $b$ with $\alpha = 2$, $\beta = -1$.
We use Lemma~\ref{UsefulLemmaStatement} to add vertices $v_1$ and
$v_2$ to the structure. Note that $v_2$ may or may not be part of 
$\Gamma_x$.  If $v_2$ is in $\Gamma_x$ then we say $\gamma_3 \geq 1$.
If $v_2$ is not in $\Gamma_x$, then we say $\gamma_3 = 0$. Although we
have only given the definition of $\gamma_3$ for the specific example
of Figure~\ref{gamma3}, it is clear how one would proceed in general.
Note that $\gamma_3$ does not necessarily have a specific value.  We
either have $\gamma_3 = 0$ or $\gamma_3 \geq 1$.

In case $\gamma_3 \geq 1$, we can find a gallery from $xC_M$
to $b\sigma(x)C_M$ that is of the same shape as a composite
gallery of some type-edge pair with $p = n$. Therefore, by 
induction, we need not worry about this case.

If $\gamma_3 = 0$, then the edge between $\Gamma_E^2$ and $\Gamma_E^3$
is forced to be a change choice.  For instance, the gallery
in Figure~\ref{gamma3} becomes the gallery in 
Figure~\ref{gamma3iszero} upon addition
of the vertex $v_3$ (which can be done using Lemma~\ref{UsefulLemmaStatement}).
We can use techniques similar to those used in the first 
part of this section to show that one can arrange for such a gallery
to fold in any desired way.
\begin{figure}
\centerline{\input{fig112.tex}}
\caption{What happens to Figure~\ref{gamma3} when $\gamma_3 = 0$}
\label{gamma3iszero}
\end{figure}

To finish the proof that $S_1 = S$, it would be necessary to consider
$(t,e) \in I_2$.  One should consider the case $b = 1$ and the cases
$2\alpha + \beta \neq \alpha + 2\beta$ separately from the cases
in which $b \neq 1$ but $2\alpha + \beta = \alpha + 2\beta$.  The 
process is similar to what we have already done, and is omitted.

\subsection{Symmetry Under a $\mathbb{Z}/3$-action}\label{Symmetry3}

In Section~\ref{Symmetry} we noticed from results developed in 
Sections~\ref{Identity} and~\ref{Non-identity} that for $SL_2$, the set 
of $w$ with non-empty $X_w (b\sigma)$ has $\mathbb{Z}/2$-symmetry
about the center of $C_M$.  We gave an {\em a priori} proof of this
fact.  The results developed for $SL_3$ in the previous sections 
of this chapter indicate that in this case the set of $w$ with
non-empty $X_w (b \sigma )$ has $\mathbb{Z} /3$-rotational symmetry
about the center of $C_M$.  We adapt the proof from 
Section~\ref{Symmetry} to give an {\em a priori} proof of the
$SL_3$ symmetry.

Let $$q = \left(
\begin{matrix}
	0 & 0 & \pi \\
	1 & 0 & 0 \\
	0 & 1 & 0
\end{matrix}
\right) \in GL_3 (F).$$
This matrix acts on $A_M$ by rotating it $120^{\circ}$ counterclockwise
about the center of $C_M$.  Let $$g = \left(
\begin{matrix}
	\pi & 0 & 0 \\
	0 & 1 & 0 \\
	0 & 0 & 1
\end{matrix}
\right) \in GL_3(F).$$
\begin{prop}\label{Z3SymmetryProp}
$\inv (gxC_M , b \sigma (gx) C_M ) = R_{120^{\circ}} (\inv (xC_M , 
b \sigma (x) C_M ))$, where $R_{120^{\circ}}$ denotes rotation
counterclockwise by $120^{\circ}$ around the center of $C_M$.  
So $\inv (g^2 xC_M , 
b \sigma (g^2 x) C_M ) = R_{240^{\circ}} (\inv (xC_M , 
b \sigma (x) C_M ))$.
\end{prop}
\begin{proof}
We know that $\inv(gxC_M , b\sigma (gx) C_M) = \rho(s^{-1} b \sigma (gx) C_M)$
where $s \in SL_3 (L)$ is chosen such that $s^{-1}gxC_M = C_M$.  We may choose
$s$ such that $s^{-1} = qx^{-1}g^{-1}$.  So $\rho (s^{-1}b \sigma (gx) C_M)
= \rho (qx^{-1}g^{-1}b \sigma (gx) C_M) = \rho(qx^{-1}b \sigma (x)C_M)$, and
this is in turn equal to
$R_{120^{\circ}}(\rho(x^{-1}b \sigma (x)C_M)) = 
R_{120^{\circ}}(\inv (xC_M, b \sigma (x) C_M ))$.
\end{proof}
This is relevant because if $y \in SL_3 (L)$ is such that $yC_M = gxC_M$,
then we have $\inv (y,b \sigma (y)) = \inv (gx C_M , b \sigma (gx) C_M)$.

One could prove an analogous result for $SL_n$, only it would
be an invariance under a $\mathbb{Z}/n$-action.  The proof is the same,
only one replaces $g$ and $q$ with the analogous $n \times n$ matrices.

\section{$GL_3$ and $PGL_3$}\label{GL3}

We now have $\inv : G(L) \times G(L) \rightarrow I \backslash G(L) / I
\simeq \tilde{W} \simeq W_a \rtimes M$, where if $G = GL_3$ then
$M = \mathbb{Z}$, and if $G = PGL_3$ then $M = \mathbb{Z}/3$.
In either case, $\inv(x,y) = (\rho (x^{-1}yC_M) , v(\det (x^{-1} y)))$,
but if $G = GL_3$ then $v(\det (x^{-1} y)) \in \mathbb{Z}$ and if
$G = PGL_3$ then $v(\det (x^{-1} y)) \in \mathbb{Z}/3$. The matrices 
$b$ listed at the beginning of Section~\ref{SL3} as representing distinct
$\sigma$-conjugacy classes in $SL_3$ still represent distinct classes
when considered as elements of $GL_3$ or $PGL_3$.  Further, any two
$\sigma$-conjugacy classes of $SL_3$ remain distinct when we pass
to $\sigma$-conjugacy classes of $GL_3$ or $PGL_3$.  It is also true that
$GL_3$ and $PGL_3$ have more $\sigma$-conjugacy classes than $SL_3$,
just as was the case for the analogous rank one groups.  Unlike in the
rank one case, we will not address these additional $\sigma$-conjugacy
classes, since for one thing we have not even addressed all the 
$\sigma$-conjugacy classes of $SL_3$ itself.

If $b$ is one of the matrices listed at the beginning of Section~\ref{SL3},
this time considered as an element of $GL_3$ or $PGL_3$, then we ask
for a description of $\{ \inv(x , b \sigma (x) ) : x \in G(L) \}$, i.e.,
we want to describe $\{ ( \rho (x^{-1} b \sigma (x) C_M ) , 
v( \det (x^{-1} b \sigma (x)))) \in W_a \ltimes M : x \in G(L) \}$.
We note that $v( \det (x^{-1} b \sigma (x))) = v( \det(b))$ is fixed for
fixed $b$ (and equal to $0$ for the $b$ we have chosen).  So we need to
examine $\{ \rho (x^{-1} b \sigma (x) C_M ) \in W_a : x \in G(L) \}$,
which, {\em a priori}, may be bigger than 
$\{ \rho (x^{-1} b \sigma (x) C_M ) \in W_a : x \in SL_3(L) \}$,
since $G(L)$ acts on $\mathcal{B}_{\infty}$ in ways that $SL_3$ does
not.  However, if $q$ is as defined in Section~\ref{Symmetry3},
and if $y = q^{-v(\det(x))}$ then $v(\det(xy)) = 0$.  Note that if
$G = PGL_3$, then $v(\det(x))$ is only determined mod $3$.
But since $q^3$ is a scalar matrix, this does not cause problems
with the definition of $y$.  Now $\det(xy) \in \mathcal{O}_L^{\times}$,
and since $(\mathcal{O}_L^{\times})^3 = \mathcal{O}_L^{\times}$,
we may further modify $xy$ on the right by a scalar matrix of determinant
$\det(xy)^{-1}$.  Therefore we may assume without loss of generality that
$\det(xy) = 1$, i.e., $xy \in SL_3(L)$.  But $xyC_M = xC_M$, so
$\sigma(xy)C_M = \sigma(x)C_M$, so $b\sigma(xy)C_M = b\sigma(x)C_M$,
so $x^{-1}b\sigma(xy)C_M = x^{-1}b\sigma(x)C_M$, so
$\rho(x^{-1}b\sigma(xy)C_M) = \rho(x^{-1}b\sigma(x)C_M)$.
But it is easy to see that $\rho(y^{-1}x^{-1}b\sigma(xy)C_M)$ is just a
rotation of $\rho(x^{-1}b\sigma(xy)C_M)$ about the center of 
$C_M$ by $(120^{\circ})(v(\det(x)))$ counterclockwise.  Since we already
proved that the set $\{ \rho(x^{-1}b\sigma(x)C_M) : x \in SL_3(L) \}$
is rotation invariant by $120^{\circ}$ and $240^{\circ}$, we see that
$\{ \rho(x^{-1}b\sigma(x)C_M) : x \in SL_3(L) \} 
= \{ \rho(x^{-1}b\sigma(x)C_M) : x \in G(L) \}$.

\section{$Sp_4$}\label{SP4}

We have applied some of the methods of Section~\ref{SL3} 
to $Sp_4$.  This section contains outlines of the application
of these methods, and the results.  The structure of the section
is the same as that of Section~\ref{SL3}.

The version of $Sp_4(E)$ that we will use for $E/F$ a field extension is
the fixed points in $GL_4(E)$ of the involution $\Theta : GL_4(E) \rightarrow
GL_4(E)$ where $\Theta(g) = M^{-1}(g^t)^{-1}M$ for 
$$M = \left(
\begin{matrix}
	0 & 0 & 0 & 1 \\
	0 & 0 & -1 & 0 \\
	0 & 1 & 0 & 0 \\
	-1 & 0 & 0 & 0 
\end{matrix} \right).$$
Note that if $T$ is the standard
maximal torus in $GL_4$, then the fixed points $T^{\Theta}$ form a split
maximal torus in $Sp_4$, and if $B$ is the standard Borel subgroup of
$GL_4$, then the fixed points $B^{\Theta}$ form a Borel for $Sp_4$.

The letter $b$ in this section will always stand for an element of 
$Sp_4$ of the form
$$ \left(
\begin{matrix}
	\pi^{\alpha} & 0 & 0 & 0 \\
	0 & \pi^{\beta} & 0 & 0 \\
	0 & 0 & \pi^{- \beta} & 0 \\
	0 & 0 & 0 & \pi^{- \alpha}
\end{matrix} \right), $$
where $\alpha \geq \beta \geq 0$.  These represent distinct 
$\sigma$-conjugacy classes in $Sp_4(L)$, but they do not constitute a
complete collection of $\sigma$-conjugacy classes. \cite{K1} \cite{K2}

\begin{figure}
\centerline{\input{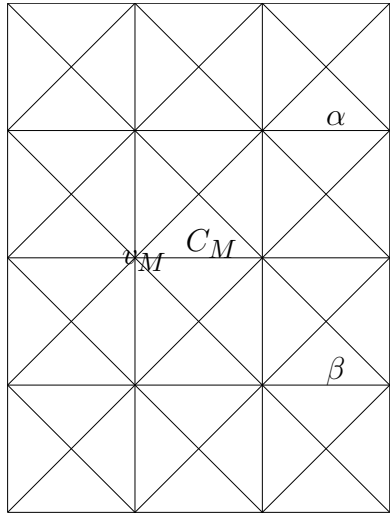}}
\caption{How the torus of $Sp_4$ acts on $A_M$}
\label{Sp4torusaction}
\end{figure}

To understand something about how $Sp_4$ acts on the main apartment of
its building, consider Figure~\ref{Sp4torusaction}.  
The main chamber and the main vertex 
are labeled in this figure, and the value of $bC_M$ with $\alpha = 1$ and
$\beta = 0$ is labeled with an $\alpha$.  The value of $bC_M$ with 
$\alpha = 0$ and $\beta = 1$ is labeled with a $\beta$.

\subsection{Standard Minimal Galleries and Composite Galleries}

\begin{figure}
\centerline{\input{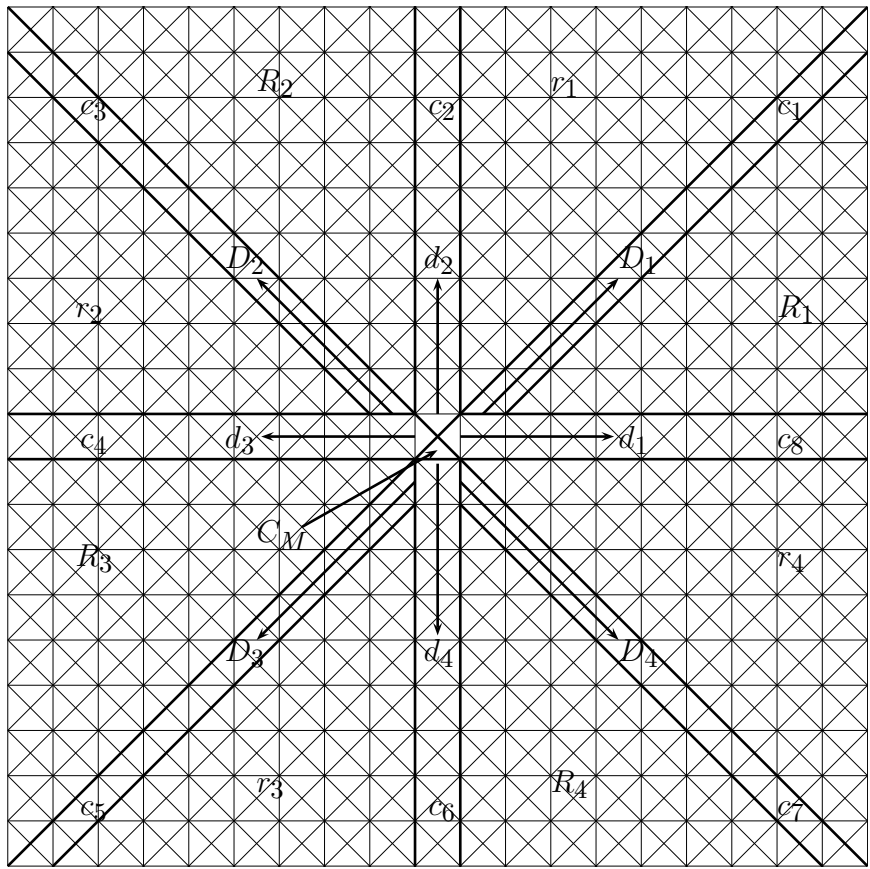}}
\caption{Primary and Secondary Directions in $Sp_4$}
\label{DirectionsSp4}
\end{figure}

We define the four {\em primary directions} $D_1$, $D_2$, $D_3$, $D_4$
and the four {\em secondary directions} $d_1$, $d_2$, $d_3$, $d_4$
in $A_M$ for $Sp_4$ as marked in Figure~\ref{DirectionsSp4}.  Given a 
chamber $E \subseteq A_M$, we define the {\em standard minimal gallery}
(SMG) between $C_M$ and $E$ as follows.  If $E$ is in one of the 
corridors marked $c_1 , \ldots , c_8$ on Figure~\ref{DirectionsSp4},
then the SMG is the unique minimal gallery from $C_M$ to $E$.  If $E$ 
is in region $R_i$, proceed first in direction $D_i$, then in 
direction $d_i$.  If $E$ is in region $r_i$, proceed first in
direction $D_i$, then in direction $d_j$, where $j \equiv i+1 \mod 4$.

As was the case for $SL_3$, one can prove that if $E$ is now allowed to
be a chamber anywhere in the building $\mathcal{B}_{\infty}$ of $Sp_4(L)$,
then there exists a unique gallery $G_E$ between $C_M$ and $E$ which is 
minimal and such that $\rho (G_E)$ is the SMG from $C_M$ to $\rho (E)$.
We define the SMG from $C_M$ to $E$ to be $G_E$ in this case.

We construct the {\em composite gallery} $\Gamma_E$ associated to
$E$ as in the $SL_3$ case.  We also define the {\em edge of departure}
in the same way as for $SL_3$

\subsection{A Conjectural Superset of the Solution Set}\label{SupersetSp4}

\begin{figure}
\setlength{\unitlength}{1in}
\begin{picture}(6,7.5)(0,0)
\centerline{\includegraphics[height=7.5in]{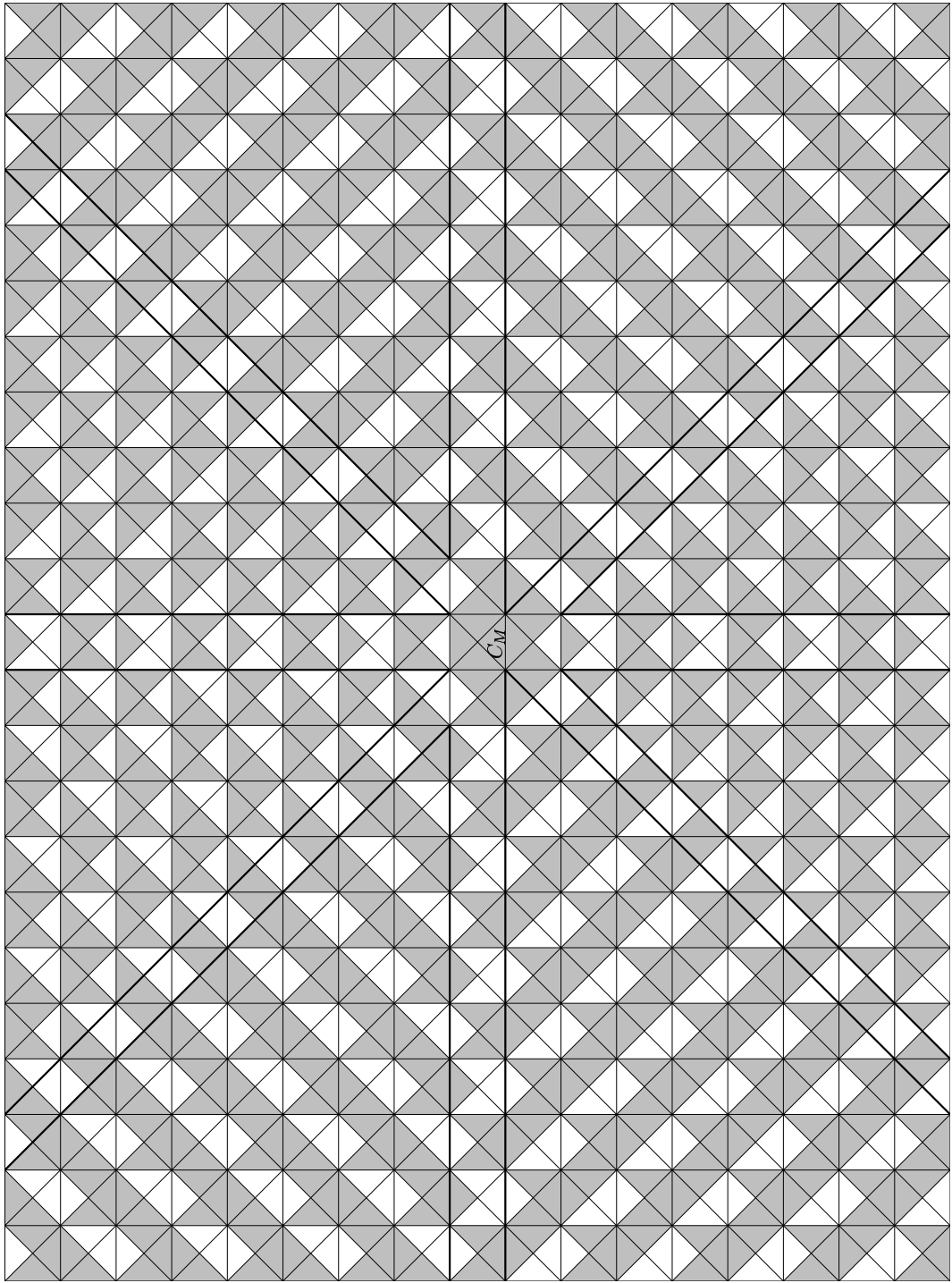}}
\end{picture}
\caption{Result for $b = 1$, $Sp_4$}
\label{SupersetResultalpha0beta0Sp4}
\end{figure}

Just as for $SL_3$, we can arrive at a superset $S_1$ of the solution
set $S$ by computing $S_1 = \cup_{(t,e)} S_{(t,e)}$, where the union 
is over all type-edge pairs $(t,e)$, and $S_{(t,e)}$ is the collection
of final chambers of possible foldings of composite galleries arising 
from SMGs of type $t$ and edge of departure $e$.  This computation
is even more prohibitively lengthy than was the case for $SL_3$, and
has only been carried out to completion for $b = 1$, for which the 
results are in Figure~\ref{SupersetResultalpha0beta0Sp4}.

\begin{figure}
\setlength{\unitlength}{1in}
\begin{picture}(6,7.5)(0,0)
\centerline{\includegraphics[height=7.5in]{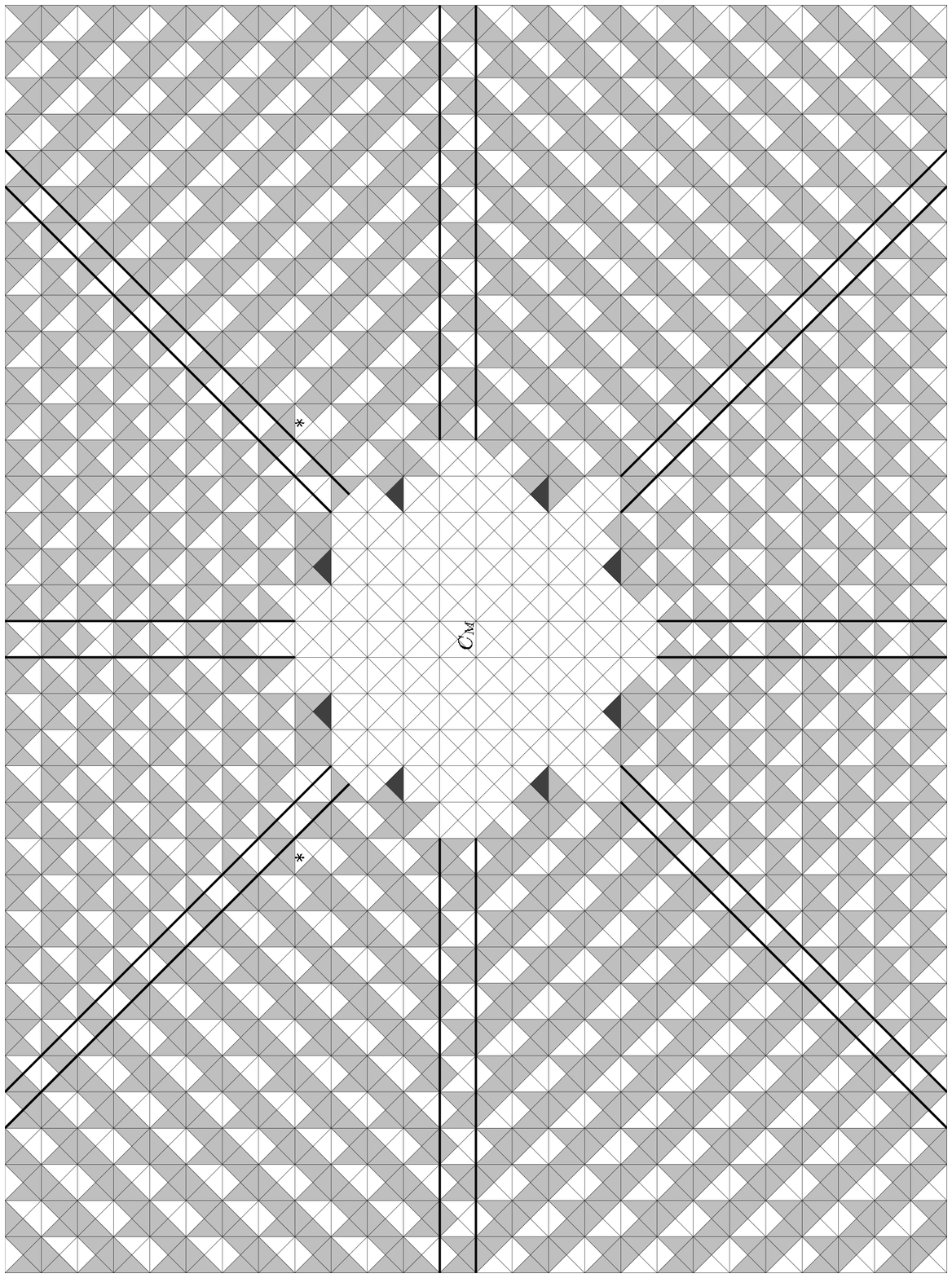}}
\end{picture}
\caption{Result for $\alpha = 3$, $\beta = 1$, $Sp_4$}
\label{SupersetResultalpha3beta1Sp4}
\end{figure}

\begin{figure}
\setlength{\unitlength}{1in}
\begin{picture}(6,7.5)(0,0)
\centerline{\includegraphics[height=7.5in]{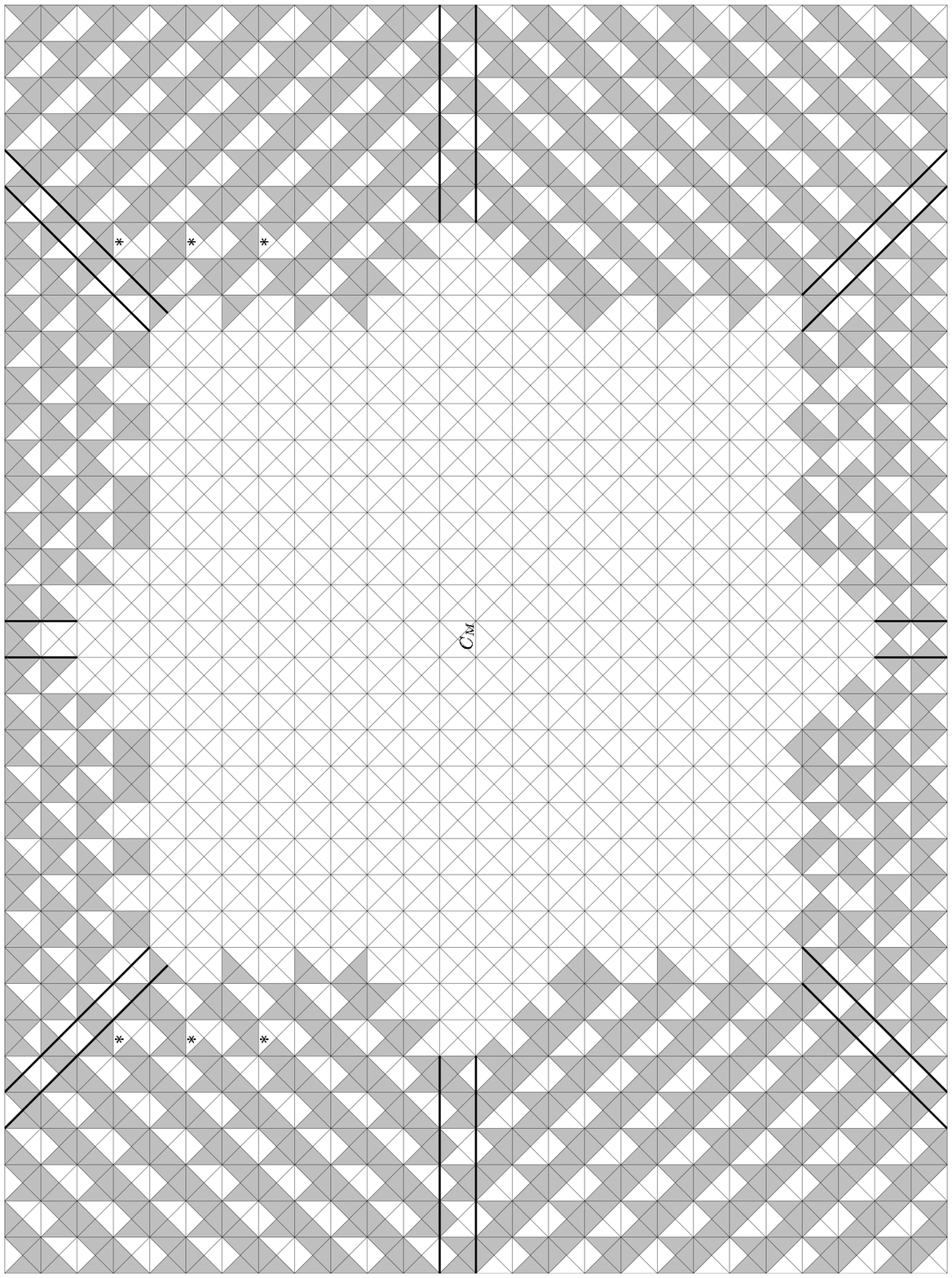}}
\end{picture}
\caption{Result for $\alpha = 6$, $\beta = 3$, $Sp_4$}
\label{SupersetResultalpha6beta3Sp4}
\end{figure}

\begin{figure}
\setlength{\unitlength}{1in}
\begin{picture}(6,7.5)(0,0)
\centerline{\includegraphics[height=7.5in]{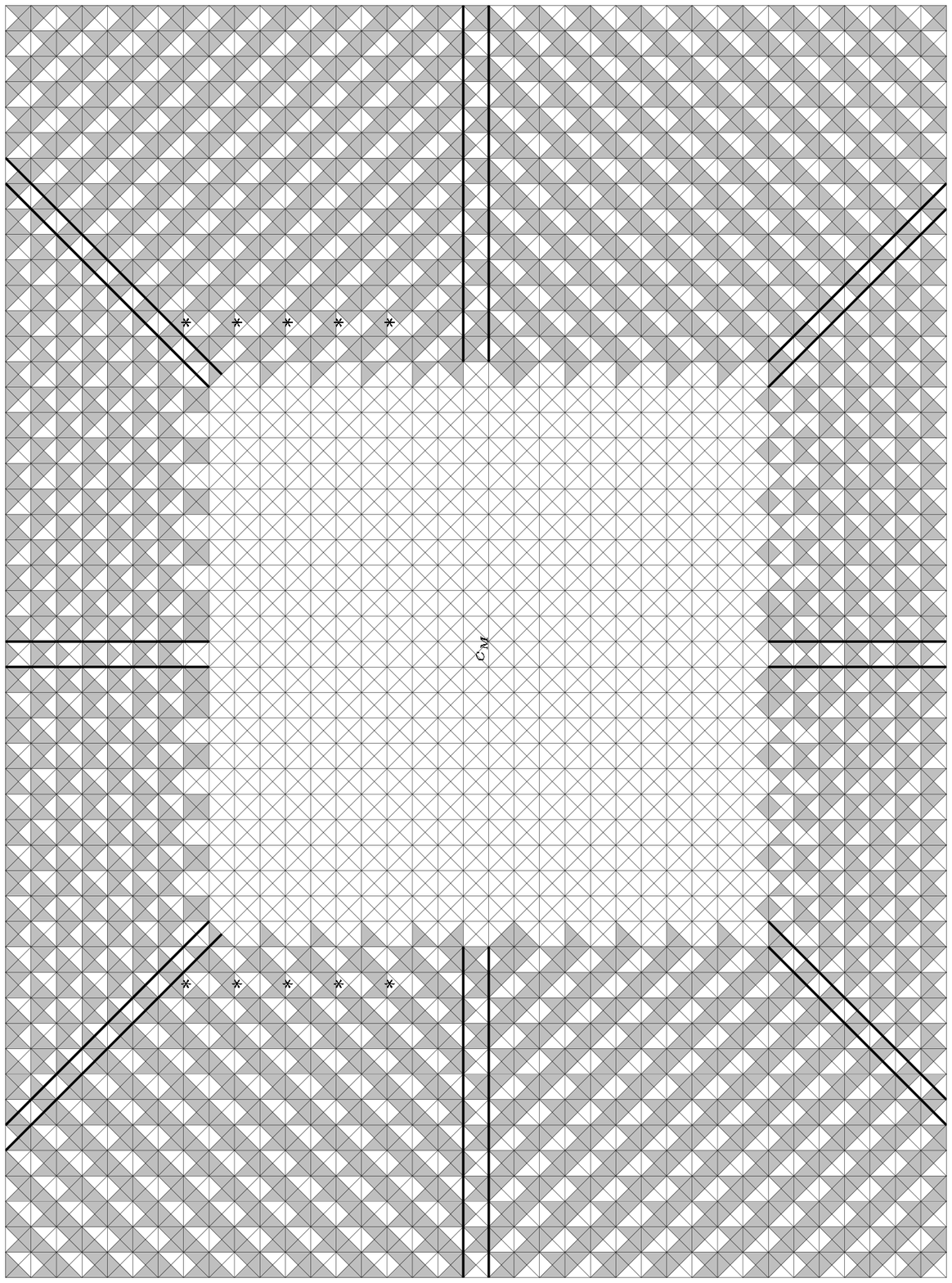}}
\end{picture}
\caption{Result for $\alpha = 6$, $\beta = 5$, $Sp_4$}
\label{SupersetResultalpha6beta5Sp4}
\end{figure}

\begin{figure}
\setlength{\unitlength}{1in}
\begin{picture}(6,7.5)(0,0)
\centerline{\includegraphics[height=7.5in]{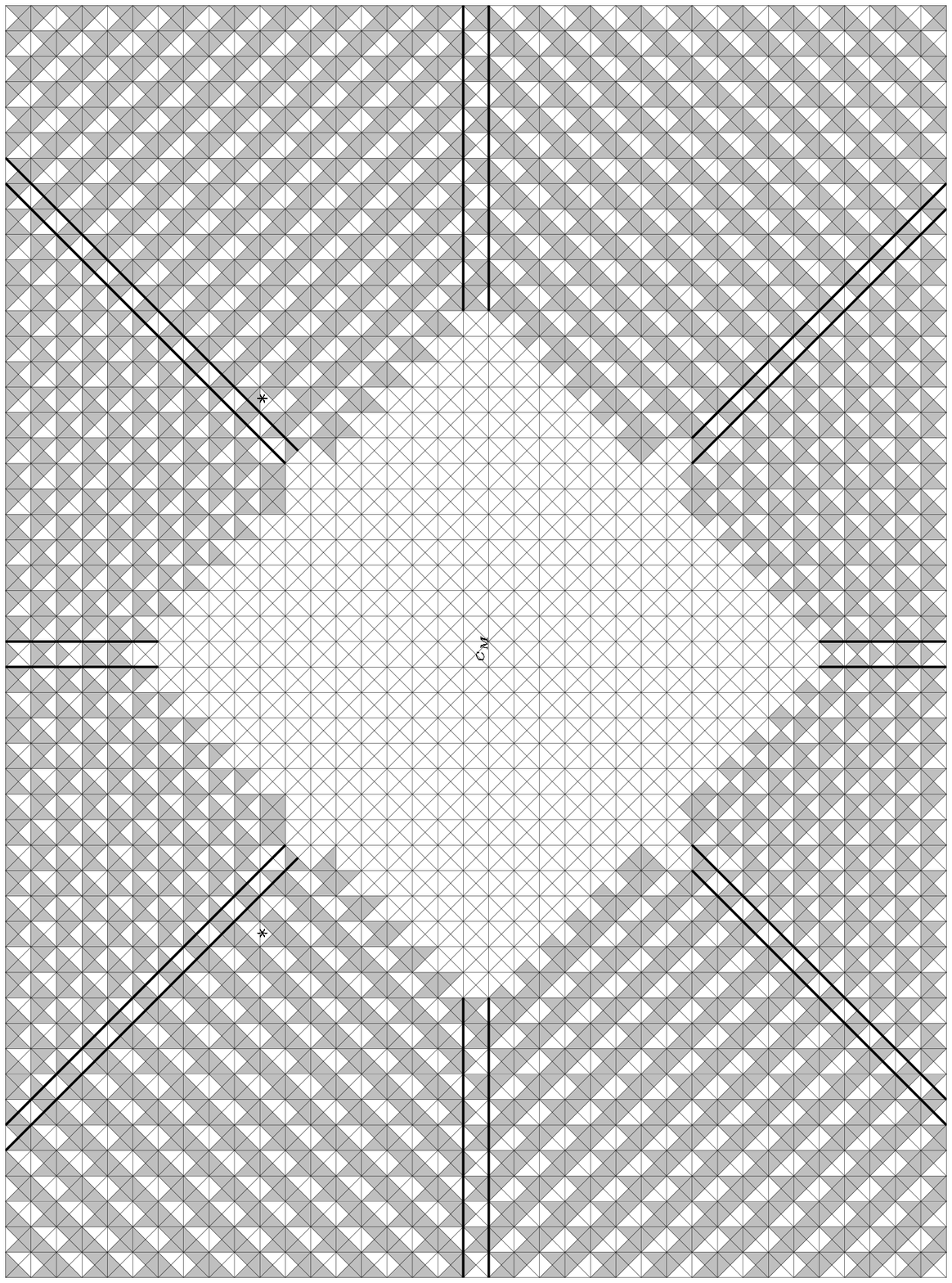}}
\end{picture}
\caption{Result for $\alpha = 7$, $\beta = 1$, $Sp_4$}
\label{SupersetResultalpha7beta1Sp4}
\end{figure}

Let $I_1$ be the infinite class of type-edge pairs $(t,e)$ with $t$ the
type of the SMG of some chamber in $r_2$, and $e$ some horizontal
departure edge.  Let $I_2$ be the infinite class of type-edge
pairs $(t,e)$ with $t$ the type of the SMG of some chamber in $c_4$,
and $e$ an arbitrary departure edge.
These $I_1$ and $I_2$ are analogous to the infinite classes of the same
name for $SL_3$.
\begin{conj}\label{conjecture}
$S_1 = \cup_{(t,e)} S_{(t,e)}$, where the union is over all type-edge
pairs in $I_1 \cup I_2$.
\end{conj}
\noindent Computation of $\cup_{(t,e)} S_{(t,e)}$ as $(t,e)$
ranges over $I_1 \cup I_2$ is a lengthy but reasonable computation.
The details of this computation for any $b$ are discussed in 
Section~\ref{RelationshipSp4}.  For $\alpha = 3$ and $\beta = 1$, the 
results are given in Figure~\ref{SupersetResultalpha3beta1Sp4}.  For
$\alpha = 6$, $\beta = 3$, the results are given in 
Figure~\ref{SupersetResultalpha6beta3Sp4}, for $\alpha = 6$, $\beta = 5$ the
results are given in Figure~\ref{SupersetResultalpha6beta5Sp4}, and for
$\alpha = 7$, $\beta = 1$ the results are given in 
Figure~\ref{SupersetResultalpha7beta1Sp4}.  The lines on these figures
and on Figure~\ref{SupersetResultalpha0beta0Sp4} are only there
to make the figures easier to look at. The chambers in 
Figure~\ref{SupersetResultalpha3beta1Sp4}
that are shaded more darkly are the chambers $w^{-1}bwC_M$ for $w \in W$.

The first piece of evidence supporting Conjecture~\ref{conjecture} is that
for $\alpha = \beta = 0$, $I_1$ and $I_2$ yield 
Figure~\ref{SupersetResultalpha0beta0Sp4}, which is the complete
computation $\cup_{(t,e)} S_{(t,e)}$, as $(t,e)$ ranges over all
type-edge pairs.

The chambers marked with a $*$ in 
Figures~\ref{SupersetResultalpha3beta1Sp4},~\ref{SupersetResultalpha6beta3Sp4},~\ref{SupersetResultalpha6beta5Sp4}, and ~\ref{SupersetResultalpha7beta1Sp4}
provide evidence against Conjecture~\ref{conjecture}, as they seem to
be holes in the pattern.  However, it could be the case that these chambers
are actually missing from the true $S_1 = \cup_{(t,e)}S_{(t,e)}$, where
the union is over all $(t,e)$.  Preliminary computations
not presented seem to provide some evidence that this is the case,
in support of Conjecture~\ref{conjecture}.  However, this evidence is 
not a strong indication.

\subsection{Relationship Between $X_w(1 \sigma)$ and $X_w (b \sigma)$ for 
$Sp_4$, and an Efficient Way of Computing Supersets}\label{RelationshipSp4}

In Section~\ref{Relationship3} we gave a more efficient method of computing
the results of the folding of $I_1 \cup I_2$.  In that case, we knew these
results gave the entire superset $S_1$ of the solution set.  One can apply
the methodology of Section~\ref{Relationship3} to $Sp_4$ to produce
an efficient way of computing the folding of $I_1 \cup I_2$ for $Sp_4$.
However, we have only conjectured that the folding of $I_1 \cup I_2$
for $Sp_4$ gives the entire superset $S_1$.

Since the methods of this section are very similar to those of 
Section~\ref{Relationship3}, we present only results.
\begin{prop}
For $b$ with $\alpha > \beta$, computing the folding results 
of $I_1 \cup I_2$ is the same as computing
the folding results of some collection of half-infinite galleries 
analogous to those pictured in 
Figure~\ref{HalfInfinitenondegSp4}, and combining those results with
their own reflection across the vertical line of symmetry of $C_M$. 
For $b$ with $\alpha = \beta$, computing the folding results of 
$I_1 \cup I_2$ is the same as computing the folding results of some 
collection of half-infinite galleries analogous to those pictured
in Figure~\ref{HalfInfinitedegSp4}, and combining these results with
their reflection across the vertical line of symmetry of $C_M$.
Although only one specific value of $b$ is pictured in each figure,
it is clear what the situation would be for arbitrary $b$.  In other
words, it is clear what the analogous collection of half-infinite
galleries would be for any $b$
\begin{figure}
\centerline{\input{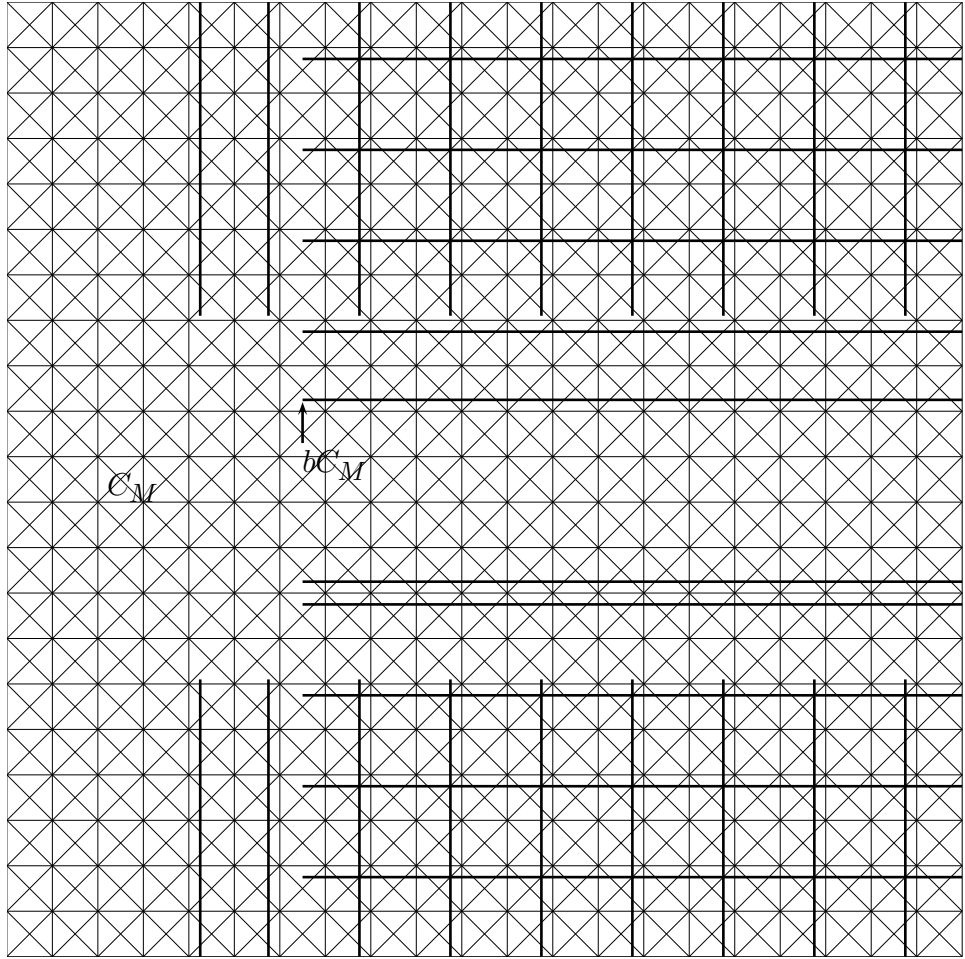}}
\caption{Half infinite galleries for $Sp_4$, $\alpha > \beta$ 
(specifically, $\alpha = 3$, $\beta = 1$)}
\label{HalfInfinitenondegSp4}
\end{figure}
\begin{figure}
\centerline{\input{fig100.tex}}
\caption{Half infinite galleries for $Sp_4$, $\alpha = \beta$ (specifically,
$\alpha = \beta = 3$)}
\label{HalfInfinitedegSp4}
\end{figure}
\end{prop}

\subsection{The Method of Kottwitz and Rapoport Applied to $Sp_4$}\label{RapoportSp4}

As in the $SL_3$ case, we produce a subset $S_2 \subseteq S$ for $b = 1$.  
Let $a \in Sp_4(F)$ be of the form
$$a = \left(
\begin{matrix}
	\pi^m & 0 & 0 & 0 \\
	0 & \pi^n & 0 & 0 \\
	0 & 0 & \pi^{-n} & 0 \\
	0 & 0 & 0 & \pi^{-m}
\end{matrix} \right)$$
where there are no conditions on $m$ and $n$.  Let $w$ be one of the 
following matrices:
$$r = \left(
\begin{matrix}
	0 & 0 & 1 & 0 \\
	1 & 0 & 0 & 0 \\
	0 & 0 & 0 & -1 \\
	0 & 1 & 0 & 0 
\end{matrix} \right) \hspace{.2in} ; \hspace{.2in}
r^2 = \left(
\begin{matrix}
	0 & 0 & 0 & -1 \\
	0 & 0 & 1 & 0 \\
	0 & -1 & 0 & 0 \\
	1 & 0 & 0 & 0
\end{matrix} \right) \hspace{.2in} ; \hspace{.2in}
r^3 = \left(
\begin{matrix}
	0 & -1 & 0 & 0 \\
	0 & 0 & 0 & -1 \\
	-1 & 0 & 0 & 0 \\
	0 & 0 & 1 & 0
\end{matrix} \right).$$
Note that $r^4 = -1$. One can
check that in fact $r \in Sp_4(F)$ by computing $\Theta(r) = r$.
The $r^i$ are representatives of Weyl group elements for $Sp_4$.
In fact, $W = D_4$, and $r$ is a representative of a generator
of the order $4$ cyclic subgroup of $W$.

\begin{figure}
\setlength{\unitlength}{1in}
\begin{picture}(6,5)(0,0)
\centerline{\includegraphics[height=5in]{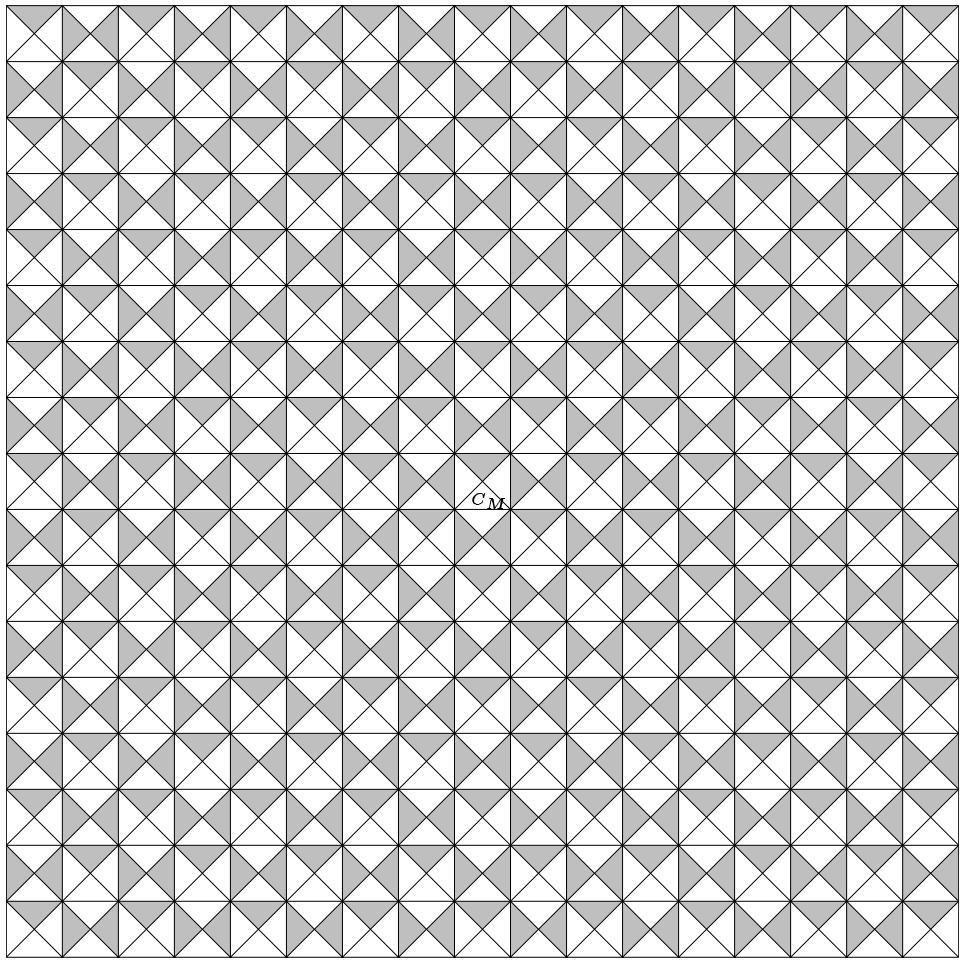}}
\end{picture}
\caption{Main results of the method suggested by Rapoport and Kottwitz, $Sp_4$}
\label{RISp4}
\end{figure}

The matrix $aw$ belongs to a basic $\sigma$-conjugacy class if there
is some $l$ such that $aw \sigma(aw) \sigma^2 (aw) \cdots \sigma^{l-1}(aw)$
is central in $Sp_4$.  Since $Sp_4$ is simply connected, $b = 1$ is 
the only basic $\sigma$-conjugacy class. \cite{K1} \cite{K2}
\begin{lem}
There exists an $l$ such that $aw \sigma(aw) \sigma^2 (aw) 
\cdots \sigma^{l-1}(aw)$ is central for any choice of $m$ and $n$.
\end{lem}
\begin{proof}
This is similar to Lemma~\ref{CentralityLemma}. 
The only difference is that $l = 4$.
\end{proof} 
\noindent Therefore, the double-$I$-cosets pictured in Figure~\ref{RISp4}
all meet the $\sigma$-conjugacy class of $1$ non-trivially.

\begin{figure}
\centerline{\input{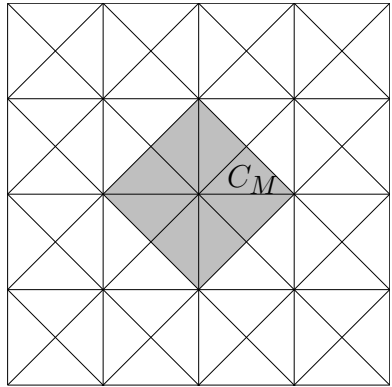}}
\caption{The $a = 1$ results of the method suggested 
by Rapoport and Kottwitz, $Sp_4$}
\label{RIIa1Sp4}
\end{figure}

\begin{figure}
\centerline{\input{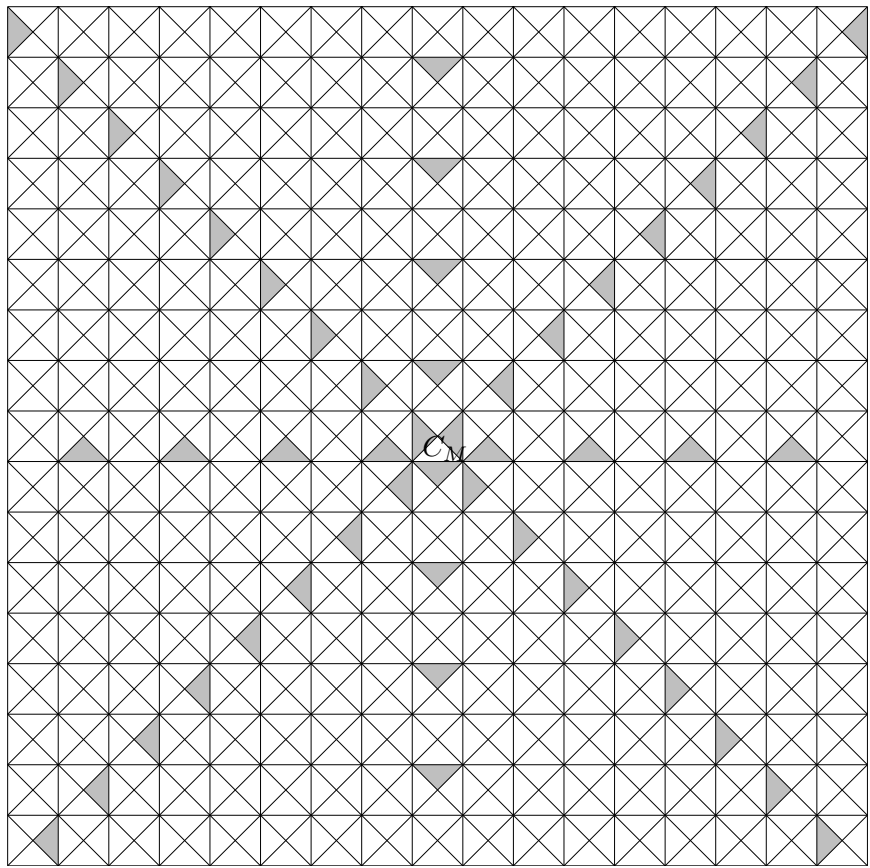}}
\caption{Results of the method suggested by Rapoport and Kottwitz
under the circumstances of Lemma~\ref{zipzang}, $Sp_4$}
\label{RIIan1Sp4}
\end{figure}

If $w$ is a representative of an element of $W$ other than $r^i$ for 
$i = 1,2,3$, and if certain conditions are placed on $m$ and $n$, then it
is still possible for $(aw)^l$ to be central for some $l$.
\begin{lem}
If $a = 1$, then for $l = 4$, $(aw)^l$ is central for any $w \in W$.
\end{lem}
\noindent This means the double-$I$-cosets corresponding to chambers
shaded in Figure~\ref{RIIa1Sp4} all intersect the $\sigma$-conjugacy class
of $b = 1$ non-trivially.  Let 
$$f = \left(
\begin{matrix}
	0 & -1 & 0 & 0 \\
	1 & 0 & 0 & 0 \\
	0 & 0 & 0 & -1 \\
	0 & 0 & 1 & 0
\end{matrix} \right).$$
\begin{lem}\label{zipzang}
If $n = 0$ then $aa^wa^{w^2}a^{w^3}w^4 = 1$ for $w = rf$.  
If $m = 0$, then $aa^wa^{w^2}a^{w^3}w^4 = 1$
for $w = r^3 f$.
If $m = -n$, then $aa^wa^{w^2}a^{w^3}w^4 = 1$
for $w = f$.
If $m = n$, then $aa^wa^{w^2}a^{w^3}w^4 = 1$
for $w = r^2 f$.
\end{lem}
\begin{proof}
Compute the relevant matrix products.
\end{proof}
\noindent The implication of this lemma is that the double-$I$-cosets
corresponding to the chambers in Figure~\ref{RIIan1Sp4} all intersect 
the $\sigma$-conjugacy class of $1$ non-trivially.

Just as the methods of Section~\ref{Rapoport} generalize to $SL_n$,
the results of this section generalize to $Sp_{2n}$.

\subsection{A Subset of the Solution Set for $b=1$ for $Sp_4$}\label{SubsetSp4}

This section is the $Sp_4$ analogue of Section~\ref{Subset}, and the setup is
the same.  In slightly different language, we choose a gallery $G$ of length
$2$ or $3$ starting at $C_M$ and having only its first chamber $C_M$ in $A_M$.
We denote the chambers of $G$ after $C_M$ by 
$G_1, \ldots, G_i$ ($i = 1$ or $2$), and we let $e$ be the edge by which 
$G_1$ is adjacent to $C_M$.  We let $\Gamma^2$ be a minimal gallery
from $e$ to $\tilde{b}e$, where $\tilde{b}C_M$ is one of the chambers in
$A_M$ obtained in the previous section. We form a composite gallery
$\Gamma$ from $G_1, \ldots, G_i$, $\Gamma^2$, and 
$\tilde{b}\sigma(G_1), \ldots, \tilde{b}\sigma(G_i)$.  The collection
of possible foldings of $x^{-1}\Gamma$ gives candidates for an 
addition to $S_2$. If there is only one possible folding then it gives
us a new element of $S_2$.

\begin{figure}
\centerline{\input{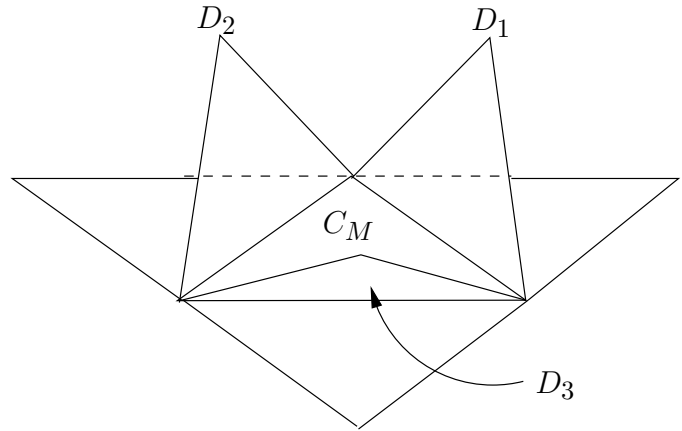}}
\caption{The choices $D_1$, $D_2$ and $D_3$ of $G_1$}
\label{ChoicesSp4}
\end{figure}

\begin{figure}
\setlength{\unitlength}{1in}
\begin{picture}(4,4.25)(0,0)
\centerline{\includegraphics[height=4.25in]{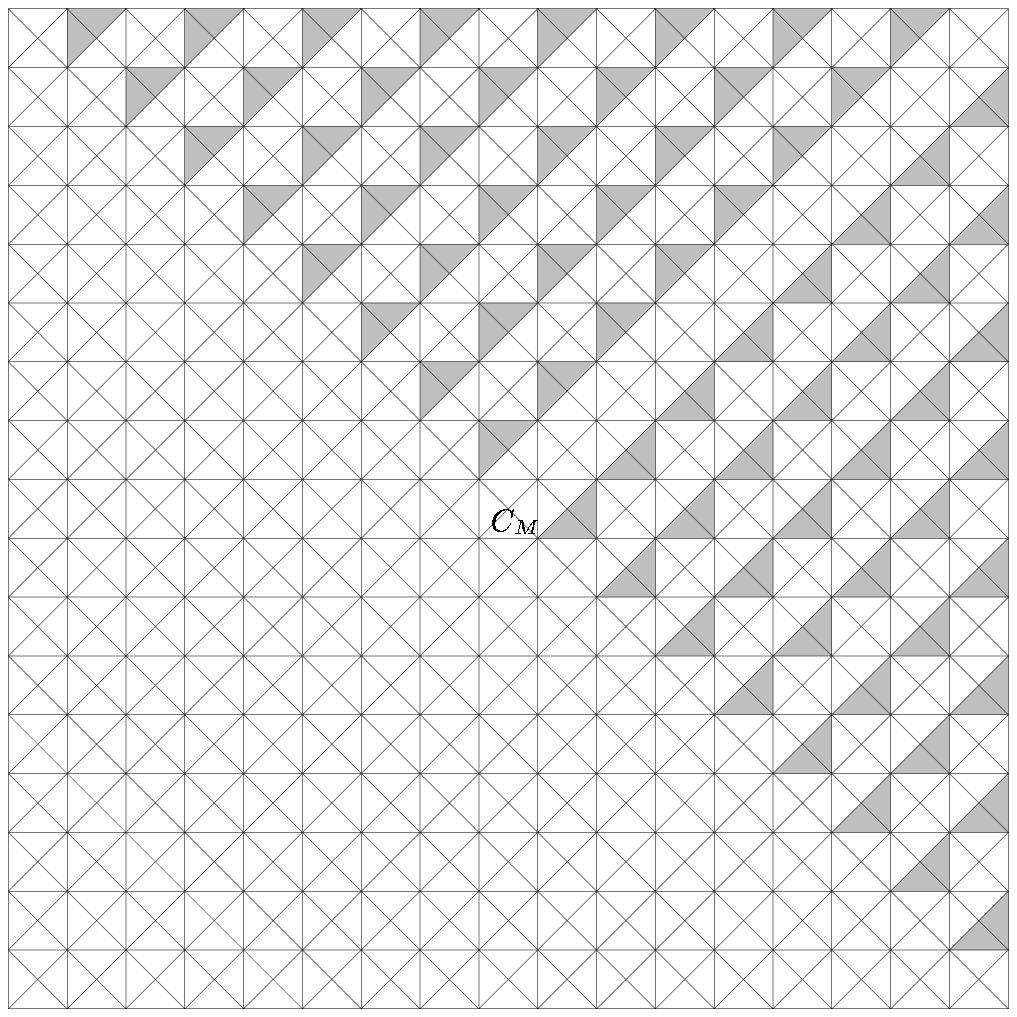}}
\end{picture}
\caption{Results for $xC_M = D_1$ and $w = r$ or for 
$xC_M = D_1$ and $w = r^3$}
\label{Sp411}
\end{figure}

\begin{figure}
\setlength{\unitlength}{1in}
\begin{picture}(4,4.25)(0,0)
\centerline{\includegraphics[height=4.25in]{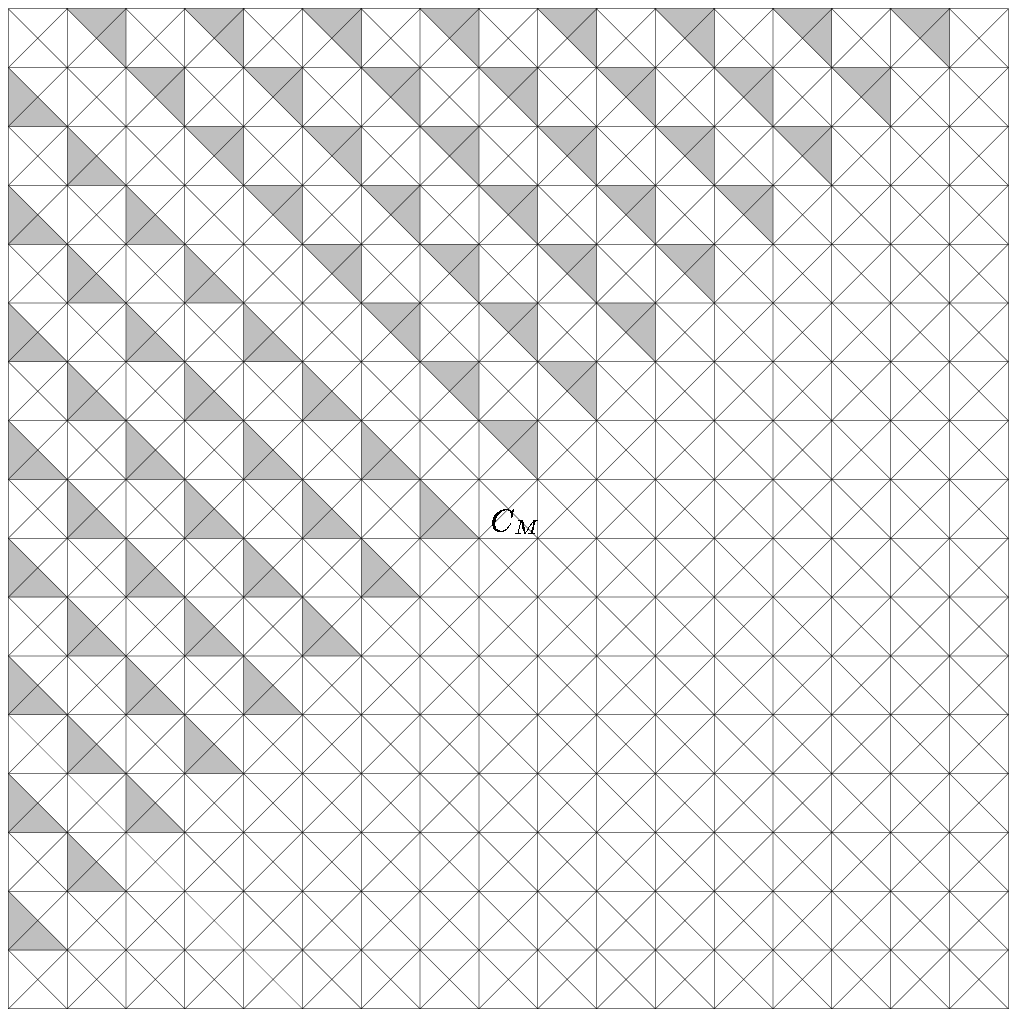}}
\end{picture}
\caption{Results for $xC_M = D_2$ and $w = r$ or for
$xC_M = D_2$ and $w = r^3$}
\label{Sp421}
\end{figure}

\begin{figure}
\setlength{\unitlength}{1in}
\begin{picture}(4,4.25)(0,0)
\centerline{\includegraphics[height=4.25in]{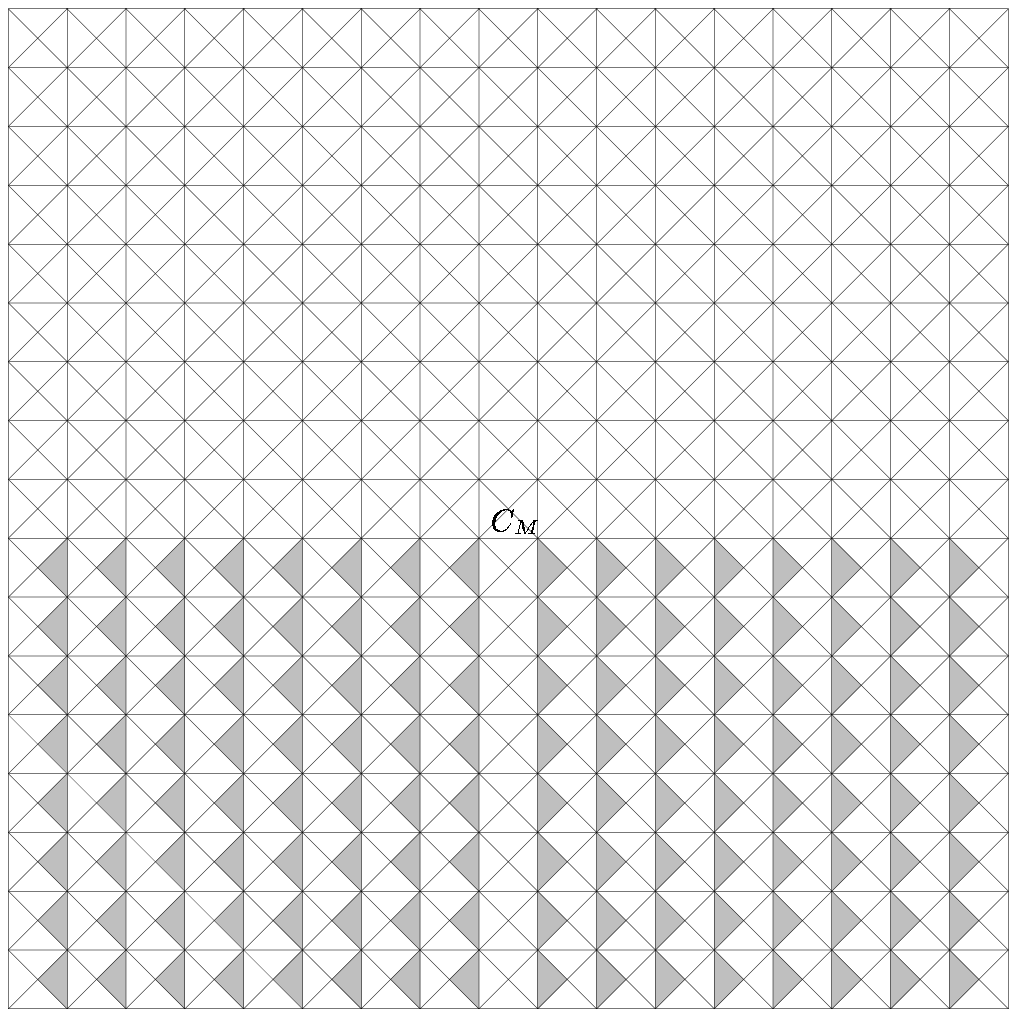}}
\end{picture}
\caption{Results for $xC_M = D_3$ and $w = r$ or for
$xC_M = D_3$ and $w = r^3$}
\label{Sp431}
\end{figure}

Starting with $i=1$ and $xC_M = G_1$ equal to one of the three chambers 
$D_1$, $D_2$, $D_3$ pictured in Figure~\ref{ChoicesSp4}, 
we get results as pictured
in Figures~\ref{Sp411},~\ref{Sp421}, and~\ref{Sp431} for 
$\tilde{b} = aw$ with $w = r$, 
Figures~\ref{Sp412},~\ref{Sp422}, and~\ref{Sp432} for $w = r^2$,
and Figures~\ref{Sp411},~\ref{Sp421}, and~\ref{Sp431} again 
for $w = r^3$.
Occasionally during these computations, the composite gallery
$\Gamma$ is not minimal, and so $\rho(x^{-1}\Gamma)$ cannot be
determined (although in all such cases $\rho(x^{-1}\Gamma)$ has
one of two possible values).  In these cases, neither of the possible
values of $\rho(x^{-1}\Gamma)$ is included in the result figures.  So
all chambers in figures~\ref{Sp411} through~\ref{Sp432} are actually in $S_2$.
\begin{figure}
\centerline{\input{fig88.tex}}
\caption{Results for $xC_M = D_1$ and $w = r^2$}
\label{Sp412}
\end{figure}
\begin{figure}
\centerline{\input{fig89.tex}}
\caption{Results for $xC_M = D_2$ and $w = r^2$}
\label{Sp422}
\end{figure}
\begin{figure}
\centerline{\input{fig90.tex}}
\caption{Results for $xC_M = D_3$ and $w = r^2$}
\label{Sp432}
\end{figure}

Combining the results obtained so far in this section with those of the
previous section gives Figure~\ref{SupersetResultalpha0beta0Sp4}, which was
already known to be the superset $S_1 \supseteq S$.  Therefore
there is no need to consider the $i=2$ case, and we see that 
$S_2 = S_1$, and so $S_2 = S = S_1$.  

\subsection{Subsets for Other $b$, and Symmetry under a $\mathbb{Z}/2$-action}\label{SymmetrySp4}

So far we have established a superset $S_1$ and a subset $S_2$ for $b=1$,
and we have noted that $S_1 = S_2$, so $S_2 = S = S_1$.  For $b \neq 1$
we have produced a conjectural superset.  
We also know a very lengthy computation
which would verify or refute the conjecture by producing the actual
superset.  We have not produced a subset $S_2$ for any $b \neq 1$.  The
process of producing such an $S_2$ for $SL_3$ was a very lengthy and 
involved proof,
and was done in Section~\ref{SubsetGeometric}.  It seems likely that a 
similar process could be carried out for $Sp_4$.  We therefore make 
the following conjecture.
\begin{conj}\label{conjecture2}
The set $\cup_{(t,e)} S_{(t,e)}$ as $(t,e)$ ranges over all type-edge
pairs is in fact equal to $S$.  In other words, $S_1 = S$.
\end{conj}
\noindent If this conjecture and Conjecture~\ref{conjecture} are both
true, the Figures~\ref{SupersetResultalpha3beta1Sp4} 
through~\ref{SupersetResultalpha7beta1Sp4} represent actual solution
sets for their respective $\sigma$-conjugacy classes.

In Section~\ref{Symmetry3}, we gave an {\em a priori} proof of 
$\mathbb{Z}/3$-rotational symmetry of the $SL_3$ solution set $S$.
This proof was based on the existence of a matrix $q \in GL_3(F)$
that acted on $A_M$ by rotating it by $120^{\circ}$ around the center
of $C_M$.  We have a matrix $q \in GSp_4(F)$ that acts on the main
apartment of the building of $Sp_4$ by flipping it 
across the (vertical) line of 
symmetry of $C_M$.  One can prove a proposition exactly
analogous to Proposition~\ref{Z3SymmetryProp} in Section~\ref{Symmetry3}
to show that the solution set $S$ for $Sp_4$ will have 
symmetry about the vertical line through the middle of $C_M$.  As the 
rotational symmetry result for $SL_3$ generalizes easily to a 
$\mathbb{Z}/n$-symmetry result for $SL_n$, so this symmetry result for $Sp_4$
generalizes to a $\mathbb{Z}/n$-symmetry result for $Sp_{2n}$.  

\section{Comments on $GSp_4$, $PSp_4$, and $G_2$}\label{GSP4}

If $G$ is $GSp_4$ or $PSp_4$, we have 
$\inv : G(L) \times G(L) \rightarrow \tilde{W} \simeq W_a \ltimes M$,
where $M = \mathbb{Z}$ for $G = GSp_4$ and $M = \mathbb{Z}/2$ for
$G = PSp_4$.  We have $\inv(x,y) = (\rho(x^{-1}yC_M), 
\frac{1}{2}v(\det(x^{-1}y)))$.  One can check that if $G = GSp_4$
and $x \in G(L)$ then $v(\det(x)) \in 2\mathbb{Z}$.  If $G = PSp_4$
and $x \in G(L)$ then $v(\det(x)) \in 2\mathbb{Z}/4\mathbb{Z} \simeq
\mathbb{Z}/2\mathbb{Z}$.  The representatives $b$ listed at the beginning
of Section~\ref{SP4} still represent distinct $\sigma$-conjugacy classes
in $G(L)$, but there are additional $\sigma$-conjugacy classes we will 
not consider.  So if $b$ is one of the matrices at the beginning of 
Section~\ref{SP4}, this time considered as an element of $G(L)$, then
we ask for a description of the set $\{(\rho(x^{-1}b\sigma(x)C_M),
\frac{1}{2}v(\det(x^{-1}b\sigma(x)))) : x \in G(L) \}$.
Using methods very similar to those of Section~\ref{GL3}, we can prove that
$\{\rho(x^{-1}b\sigma(x)C_M) : x \in G(L)\} 
= \{\rho(x^{-1}b\sigma(x)C_M) : x \in Sp_4(L)\}$.  It is already clear
that $v(\det(x^{-1}b\sigma(x))) = v(\det(b))$ is fixed for fixed $b$
(and equal to $0$ for the $b$ we have chosen).

\begin{figure}
\centerline{\input{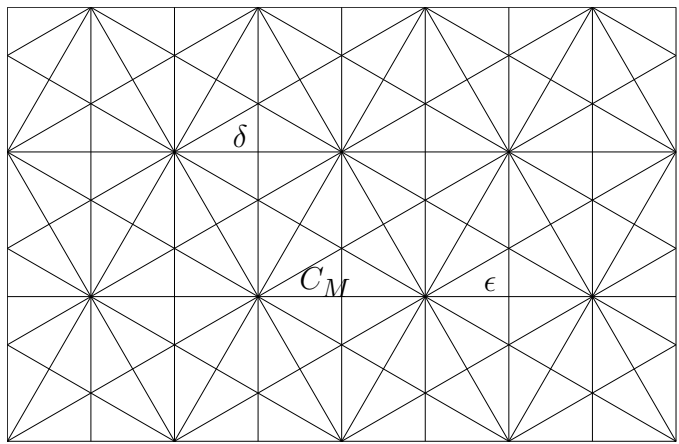}}
\caption{\hspace{0.05in}The translations $\epsilon$ and $\delta$ of $A_M$ for $G_2$}
\label{ed}
\end{figure}

We would also like to make some comments about $G_2$.  We first note that
$W = D_6$, the group of symmetries of a hexagon, and that 
$W_a = \mathbb{Z}^2 \rtimes W$. We ask for a description of 
$\{\inv(x,\sigma(x)) : x \in G_2(L) \}$.  Let $\epsilon$ and 
$\delta$ be the two standard generators of $\mathbb{Z}^2 \subseteq
\mathbb{Z}^2 \ltimes W \simeq W_a$, where $\epsilon$ and $\delta$
act on $A_M$ by translating it so that $C_M$ goes, respectively, 
to the chambers marked with an $\epsilon$ and a $\delta$ in
Figure~\ref{ed}.

Let $r \in W$ be the element that rotates $A_M$ counterclockwise
by $60^{\circ}$ around $v_M$, and let $f \in W$ be the element 
that flips $A_M$ about the horizontal line through $v_M$.  We know
that $\epsilon$, $\delta$, $r$ and $f$ generate $W_a$, 
$\epsilon$ and $\delta$ commute, and $rf = fr^5$.  One can check that
$r\epsilon = \epsilon \delta r$, $r \delta = \epsilon^{-1} r$,
$f\epsilon = \epsilon f$, and $f\delta = \epsilon^{-1}\delta^{-1}f$.

Using a similar process to that used in Sections~\ref{Rapoport} 
and~\ref{RapoportSp4}, if $a = \epsilon^s \delta^t$ and 
$w \in W$, and if $(aw)^l = 1$ in
$W_a$ for some $l$, then $aw \in \{ \inv (x, \sigma(x)) : x \in G_2(L) \}$.
We get:
\begin{lem}
For any integers $s$, $t$ and any $1 \leq n \leq 5$, $\epsilon^s
\delta^t r^n \in \{ \inv (x,\sigma(x)) : x \in G_2(L) \}$.
\end{lem}
\begin{proof}
Just simplify $(\epsilon^s\delta^t r^n)^6$ for each $n = 1,2,3,4,5$
using the rules listed in the previous paragraph.  For instance,
if $n = 1$, then $(\epsilon^s\delta^t r)^6 = 
(\epsilon^s\delta^t r)(\epsilon^s\delta^t r)(\epsilon^s\delta^t r)^4
= (\epsilon^s\delta^t)(\epsilon^{s-t}\delta^s r^2)(\epsilon^s\delta^t r)^4
= (\epsilon^s\delta^t)(\epsilon^{s-t}\delta^s)(\epsilon^{-t}\delta^{s-t})
(\epsilon^{-s}\delta^{-t})(\epsilon^{-s+t}\delta^{-s})(\epsilon^t\delta^{-s+t})
r^6 = 1$.
\end{proof}

One could also create lemmas for $G_2$ analogous to Lemma~\ref{balllemma} and 
Lemma~\ref{RapRes} in a similar way.

These results give a starting point for the subset $S_2$ of the
solution set $S$ for $G_2$ with $b = 1$.  One could enlarge
$S_2$ using methods similar to those in Section~\ref{Subset} and
Section~\ref{SubsetSp4}.  One could produce a superset
$S_1 \supseteq S$ by folding galleries in an analogous way to 
Section~\ref{Superset} and Section~\ref{SupersetSp4}.  This has not been done,
so we cannot say for sure if the resulting $S_1$ and $S_2$ are
equal, although we guess that they would be.

We also note that $G_2$ has no analogue to the rotational invariance
developed for $SL_3$ in Section~\ref{Symmetry3}, or the flip invariance
developed for $Sp_4$ in Section~\ref{SymmetrySp4}, although the general 
results that will be developed in Section~\ref{Invariance} will hold for $G_2$.

\section{Invariance Properties of Solution Sets}\label{Invariance}

In this section we prove some results that give invariance properties 
of the solution set $S$. These results hold true if $S$ is the solution set 
$\{ \inv (x,b \sigma (x)) : x \in G (L) \}$ for
any simply-connected group $G$ and any $\sigma$-conjugacy class $b$.
We apply the invariance results of this section to $SL_3$, $b = 1$
under the assumption that one knows whether $X_w (1 \sigma)$ is empty
or non-empty only for any $w$ in a certain subset of $W_a = \tilde{W}$.
We will see that knowledge on whether $X_w (b \sigma )$ is empty
can be obtained on a much larger class of $w$.

Let $G$ be a simple simply-connected group, let $W_a$ be 
its affine Weyl group, and
let $A_M$ be the main apartment of its building.  Suppose $G$ has rank $n$.
Let $C_M$ be the main chamber and let $v_M$ be the main vertex.  Let 
$L_1, \ldots, L_n$ be the hyperplanes in $A_M$ that contain $v_M$
and that intersect $C_M$ in an $n-1$ dimensional simplex.  Let $L_{n+1}$
be the hyperplane in $A_M$ that intersects $C_M$ in an $n-1$ dimensional
simplex, but that does not contain $v_M$.  Let $s_i$ be reflection of 
$A_M$ about $L_i$.  Then the $s_i$ generate $W_a$ as a Coxeter group.

If $D$ is a chamber in $A_M$, define $p_i D$ to be the chamber
obtained by reflecting $D$ about the wall $\tilde{L}_i$, where
$\tilde{L}_i$ is parallel to $L_i$, and intersects $D$ in an $n-1$ 
dimensional simplex.  So $p_i$ is a map from the set of chambers in 
$A_M$ to itself.  The elements $s_i \in W_a$ are also maps of this
kind, and as such we can state the following lemma:
\begin{lem}\label{shisto}
If $i \neq j$ then $p_i s_j = s_j p_k$ for some $k$.  We also
have $p_i s_i = s_i p_i$.
\end{lem}
This can be easily checked, and is valuable because a gallery in $A_M$
starting at $C_M$ is just a sequence of the $p_i$.  We can now prove:
\begin{lem}
If $w \in W_a$ then there is a one-to-one correspondence between
Coxeter expansions of $w$ (using the $s_i$) and galleries from 
$C_M$ to $wC_M$.  A Coxeter expansion of length $m$ corresponds to
a length $m+1$ gallery, so minimality is preserved under the 
correspondence.
\end{lem}
\begin{proof}
If $w = s_{i_1} \cdots s_{i_m}$ then $wC_M = s_{i_1} \cdots s_{i_m} C_M
= s_{i_1} \cdots s_{i_{m-1}} p_{i_m} C_M = 
p_{k_m}s_{i_1}\cdots s_{i_{m-1}} C_M$, where the last equality is achieved
by using Lemma~\ref{shisto} $m-1$ times.  We then proceed: 
$p_{k_m} s_{i_1} \cdots s_{i_{m-1}} C_M 
= p_{k_m} s_{i_1} \cdots s_{i_{m-2}} p_{i_{m-1}} C_M 
= p_{k_m} p_{k_{m-1}} s_{i_1} \cdots s_{i_{m-2}}C_M = \cdots
=p_{k_m} p_{k_{m-1}} \cdots p_{k_1} C_M$.  Note that length is preserved,
and the process is reversible.  
\end{proof}

Let $l(w)$ denote the length of $w \in W_a$ as an element of the Coxeter
group $W_a$.  Then we have the following results.  Note that 
Proposition~\ref{bij} is the affine analogue of some parts of the proof of 
Theorem $1.6$ in \cite{DL1}.
\begin{prop}\label{bij}
If $s = s_i$ for some $i$ and if $l(sws) = l(w)$ then there exists
a bijective map from $X_w (b \sigma )$ to $X_{sws} (b \sigma)$.
\end{prop}
\begin{proof}
We define $P$ to be the parallelogram spanned by $C_M$ and $wC_M$ (so
it is the intersection of all apartments containing $C_M$ and $wC_M$,
or, alternatively, the union of all minimal galleries from $C_M$
to $wC_M$).  Because of the fact that $l(sws) = l(w)$, one can show
that exactly one of $wsC_M$ and $sC_M$ is in $P$.  If 
$x \in X_w (b \sigma )$, then let $g \in I$ be such that 
$gx^{-1}b \sigma (x) C_M = wC_M$.  Consider $xg^{-1}P$, which
contains $xg^{-1}C_M = xC_M$ and $xg^{-1}wC_M = b \sigma(x) C_M$.

\noindent {\em Case 1:} $sC_M \subseteq P$. Then let $y = xg^{-1}s$,
so $yC_M = xg^{-1}sC_M$.  Since $yC_M$ is adjacent to $xC_M$,
$b \sigma(y)C_M$ is adjacent to $b \sigma(x) C_M$, and not a part of
$xg^{-1}P$.  One can now see that $\inv (y , b \sigma (y) ) = sws C_M$.

\noindent {\em Case 2:} $wsC_M \subseteq P$.  Then let 
$yC_M = \sigma^{-1} b^{-1} x g^{-1} wsC_M$, so we have that
$b \sigma (y)C_M = xg^{-1}wsC_M$.  Since $wsC_M$ is adjacent to $wC_M$,
$xg^{-1}wsC_M$ is adjacent to $b \sigma (x) C_M$, so $y C_M$ is 
adjacent to $xC_M$ along the edge $xg^{-1}e$, where $e$ is the
edge of adjacency of $C_M$ and $sC_M$.  One can now see that
$\inv (y , b \sigma (y)) = sws C_M$.

Now let $\Gamma_1 : X_w (b \sigma ) \rightarrow X_{sws} (b \sigma)$
be defined so that $\Gamma_1 (x C_M) = yC_M$.  I claim that this
is a bijection. To check this, note that the roles of $w$ and $sws$
in the construction of $\Gamma_1$ were symmetric, so by replacing
each with the other we get a map 
$\Gamma_2 : X_{sws} (b \sigma ) \rightarrow X_w (b \sigma)$.  One
can check that $\Gamma_1$ and $\Gamma_2$ are inverses.
\end{proof}

\begin{prop}
If $s = s_i$ for some $i$ and if $l(w) > l(sws)$ (and so 
$l(sws)= l(w) - 2$) then there is a surjective map from $X_w(b \sigma)$
to $X_{sws} (b \sigma ) \cup X_{sw}(b \sigma)$.
\end{prop}
\begin{proof}
Again, let $P$ be the parallelogram spanned by $C_M$ and $wC_M$.  The fact
that $l(w) > l(sws)$ means that both $wsC_M$ and $sC_M$ are in $P$.
As in the previous proposition, if $x \in X_w (b \sigma)$, let 
$g \in I$ be such that $gx^{-1}b \sigma (x) C_M = wC_M$.  Let $e$ be 
the edge of adjacency between $wC_M$ and $wsC_M$, and consider 
$xg^{-1}P$.  One can see that $b \sigma (xg^{-1} s) C_M$ is 
adjacent to $b \sigma (x)C_M = xg^{-1} wC_M$ via $xg^{-1}e$, but we may
have $b \sigma (xg^{-1} s)C_M = xg^{-1}wsC_M$ or we may have
$b \sigma (xg^{-1}s)C_M \neq xg^{-1}wsC_M$ (at least {\em a priori},
either the equality or the inequality could possibly hold, and it will
turn out that both actually do arise).

\noindent {\em Case 1:} $b \sigma (xg^{-1} s)C_M = xg^{-1}wsC_M$.
Then if we let $y = xg^{-1}s$, we have that
$\inv (y , b \sigma (y)) = \inv (xg^{-1}s, xg^{-1}ws)
= \inv(s, ws) = \inv(1,sws) = sws$, so $y \in X_{sws} (b \sigma)$.

\noindent {\em Case 2:} $b \sigma (xg^{-1}s) C_M \neq xg^{-1}wsC_M$.
Then if $y = xg^{-1}s$, we have $\inv(y, b \sigma (y)) = \inv (s,w)
= \inv(1,sw) = sw$, so $y \in X_{sw}(b \sigma)$.

So we have a map 
$\Gamma: X_w (b \sigma) \rightarrow X_{sws}(b \sigma)\cup X_{sw} (b \sigma)$
where $\Gamma(xC_M) = yC_M$.  We must show that $\Gamma$ is surjective.
Take $z \in X_{sw}(b \sigma)$, so $\inv (z, b \sigma (z)) = sw$.
So $\rho_{sC_M} (sz^{-1} b \sigma (z) C_M) = wC_M$, where $\rho_{sC_M}$
is the retraction of $\mathcal{B}_{\infty}$ onto $A_M$ centered at $sC_M$.
So there exists $g \in sIs$ such that $gsz^{-1}b \sigma (z) C_M = wC_M$.
Consider $zsg^{-1}P$, which contains $zC_M$ and 
$b \sigma (z) C_M = zsg^{-1}wC_M$.  Note that if $e$ is the edge of
adjacency of $C_M$ and $sC_M$, then $zC_M$ contains $zsg^{-1}e$.
Choose a chamber $xC_M$ containing $e$, but not equal to $zC_M$.
Then $b\sigma(x)C_M$ is adjacent to $b\sigma(z)C_M$.  We may assume
without loss of generality that $b\sigma(x)C_M \nsubseteq zsg^{-1}P$
(we can arrange this by choosing an appropriate $xC_M$).  Then
$x \in X_w(b\sigma)$ and $\Gamma(xC_M) = zC_M$.

If $z \in X_{sws}(b \sigma)$, then $\rho_{sC_M} (sz^{-1}b \sigma(z) C_M)
= wsC_M$.  So there exists $g \in sIs$ such that 
$gsz^{-1}b \sigma(z)C_M = wsC_M$.  Consider $zsg^{-1}P$, which contains $zC_M$
and $b \sigma(z)C_M = zsg^{-1}wsC_M$.  Let $e$ be the edge between
$C_M$ and $sC_M$.  We know $zsg^{-1}e \subseteq zC_M$, and 
$b \sigma (e) = zsg^{-1}we$.  So choose $xC_M$ containing $e$ and not in 
$P$.  Then $b \sigma(x) C_M$ contains $b \sigma(e)$.  We can require
without loss of generality that $b \sigma(x)C_M \nsubseteq P$.
Then $\Gamma(xC_M) = zC_M$.  This proves surjectivity.
\end{proof}

\begin{cor}
If $l(sws) > l(w)$ (so $l(sws) = l(w) + 2$) then there is a surjective
map from $X_{sws}(b \sigma )$ to $X_w (b \sigma ) \cup X_{ws} (b \sigma)$.
\end{cor}
\begin{proof}
Just apply the previous proposition with $sws$ in place of $w$.
\end{proof}

\begin{cor}\label{lengthoneapp}
If $l(sw) > l(w)$, $l(ws) < l(w)$ and $X_w (b \sigma) \neq \emptyset$
then $X_{sw} (b \sigma ) \neq \emptyset$.  If $l(sw) < l(w)$, $l(ws) > l(w)$,
and $X_w(b \sigma ) \neq \emptyset$, then $X_{ws} (b \sigma) \neq \emptyset$.
\end{cor}
\begin{proof}
This follows from the two propositions.
\end{proof}

Note that we used this last result repeatedly in Section~\ref{Subset} and 
Section~\ref{SubsetSp4} for $SL_3$ and $Sp_4$.

\begin{figure}
\setlength{\unitlength}{1in}
\begin{picture}(6,5)(0,0)
\centerline{\includegraphics[width=6in]{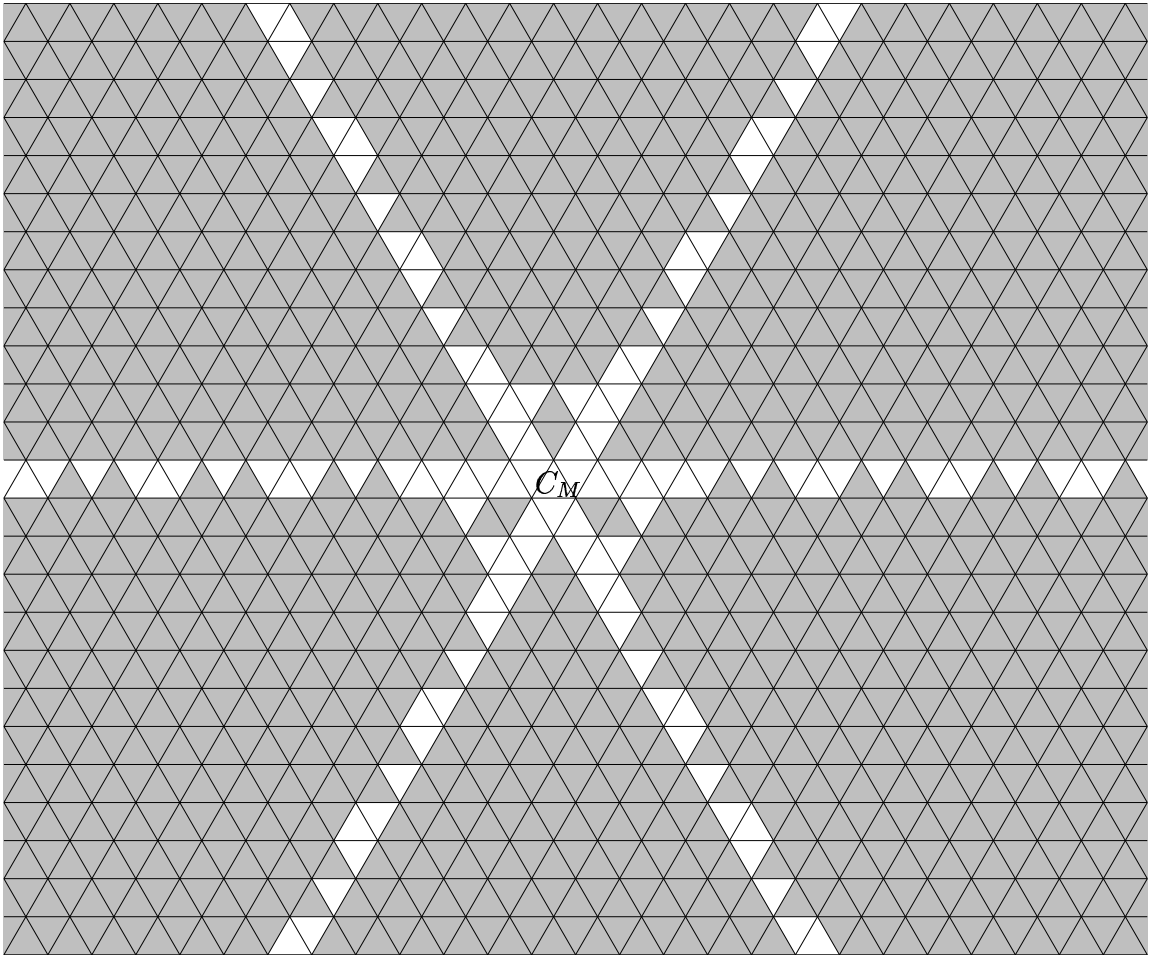}}
\end{picture}
\caption{\hspace{0.05in}Understood region after applying invariance properties to region $R$}
\label{RegionRInvariantResults}
\end{figure}

These results are valuable because they can be used on a wide variety
of groups $G$ and $\sigma$-conjugacy classes $b$ to increase
partial information on the nature of the solution set 
$\{ \inv(x,b \sigma (x)) : x \in G(L) \}$.  To illustrate the possible utility
of the above propositions, we consider the hypothetical case that 
for $SL_3$, $b = 1$, one knows for any $w$ in region $R$ on 
Figure~\ref{RegionR} whether $X_w (1\sigma)$ is or is not empty.
So we assume that part of Figure~\ref{TotalResults0_0}, the solution set for 
$SL_3$, $b=1$, is given.  In this case, one can repeatedly apply the results
of this section, together with the rotational invariance result of 
Section~\ref{Symmetry3}, to the point where it is known whether $X_w(b \sigma)$
is empty or non-empty for each $w$ corresponding to a chamber shaded
in Figure~\ref{RegionRInvariantResults}.  We note that this represents
nearly complete information about the solution set of $SL_3$, $b = 1$, 
since nearly every chamber in that figure is shaded.  The results
of this section could be applied with the same level of effectiveness
if $b \neq 1$, or for $Sp_4$.  Our results could even be applied to
other groups, including higher rank groups.
\begin{figure}
\centerline{\input{fig91.tex}}
\caption{\hspace{0.05in}Region $R$}
\label{RegionR}
\end{figure}

Also note that Sections~\ref{Subset} and~\ref{SubsetSp4} were essentially just 
Corollary~\ref{lengthoneapp} applied to the partial information of 
Sections~\ref{Rapoport} and~\ref{RapoportSp4}.

\singlespacing
\pagebreak

\end{document}